\documentclass[11pt]{ucthesis}
\usepackage{epsf}
\usepackage{graphicx}
\usepackage{url}
\usepackage{amsmath}
\usepackage{amssymb}
\usepackage{float}
\newcommand{\urlBiBTeX}[1]{\url{#1}}
\newtheorem{proposition}{Proposition}[section]
\def\dsp{\def\baselinestretch{2}\large\normalsize}
\dsp
\begin{document}

% Declarations for Front Matter

\title{Finite Volume Methods For Incompressible Flow}
\author{Darryl M. Whitlow}
\degreeyear{2001}
\degreesemester{November}
\degree{Doctor of Philosophy}
\chair{Professor Jean-Jacques Chattot}
\othermembers{Professor Mohamed Hafez\\
Professor John Hunter}
\numberofmembers{3}
\prevdegrees{}
\field{Applied Mathematics}
\campus{Davis}

\maketitle

\begin{frontmatter}

\approvalpage
\begin{dedication}
\null\vfil
{\large
\begin{center}
To my family
\end{center}}
\vfil\null
\end{dedication}

\begin{acknowledgements}
I would first like to thank my family and friends for their support throughout a period of tribulation,
sacrifice and yes, reward.  Let me also thank each of my advisors who assisted me in the completion of my
dissertation.  While their intuition and accurate criticism is to be admired, it is their compassion, patience,
humility and imagination that I will cherish forever.
\paragraph{}
Also, I would like to acknowledge the support rendered by a fellowship from the National Physical Science
Consortium and by stipends and grants from NASA-Glenn and NASA-Ames.
\end{acknowledgements}
\begin{abstract}
Two finite volume methods are derived and applied to the solution of problems of incompressible
flow.  In particular, external inviscid flows and boundary-layer flows are examined.  The first
method analyzed is a cell-centered finite volume scheme.  It is shown to be formally first order accurate on equilateral triangles
and used to calculate inviscid flow over an airfoil.  The second method is a vertex-centered least-squares
method and is second order accurate.  It's quality is investigated for several types of inviscid flow
problems and to solve Prandtl's boundary-layer equations over a flat plate.  Future improvements and extensions
of the method are discussed.
\end{abstract}

\tableofcontents
\end{frontmatter}
\chapter{Introduction}
Finite volume methods have been used extensively in the field of computational fluid dynamics.  They have several advantages.  Among these are their ability to be used on both structured and unstructured meshes.  Usually, the methods can be derived in a 
straightfoward manner yet yield robust schemes with favorable properties of conservation of the fluxes in the flow field.
\paragraph{}
The first part of this work concerns a finite volume method where a variable is piecewise constant on a triangle.  The method yields a weak formulation of all necessary derivatives to achieve an accurate solution.  The method is applied to the solution of Laplace's equation.  The resulting difference equation for Laplace's equation is shown to have a truncation error that is in general no more than first order accurate on equilateral triangles.  However, it will be demonstrated that the error of the actual solution can in general achieve an accuracy as high as second order on arbitrary triangles.
\paragraph{}
In the second part of this work, a least-squares finite volume method is studied.  There has been a resurgence of interest in least-squares type methods.  These methods deal with the construction of some functional whose minimization is a solution of the governing equations.  Typically, these methods are derived in the finite element context though it has been done in a finite difference setting \cite{chattot2}. Finite volume approaches to the construction of the functional are relatively young in the CFD community and started with some interest in 1996 \cite{hub:least}.  
The methodology involves first
integrating the governing equations over each element in the computational domain.  Next,
the functional is constructed by summing the
squares of the discrete governing equations over each element.  The functional
is then minimized with respect to the unknowns.  In general, this results in a nonlinear system which is
solved by Newton's method.  The minimization is a solution of the governing equations.  Careful consideration must be taken to enforce the proper boundary conditions for external flow problems so that the resulting minimization gives a physically plausible solution.
\chapter{A Piecewise Constant Approximation Method}
In this chapter, a finite volume method is derived for the solution of Laplace's equation 
in two dimensions.  The method is applied to potential flow around a Joukowski airfoil in terms of the 
stream function. The solution is reached iteratively, using the Gauss-Seidel method, on an
unstructured triangular grid.

\section{Derivation}
\begin{figure}[htbp]
        \centering
         \includegraphics[totalheight=3in]{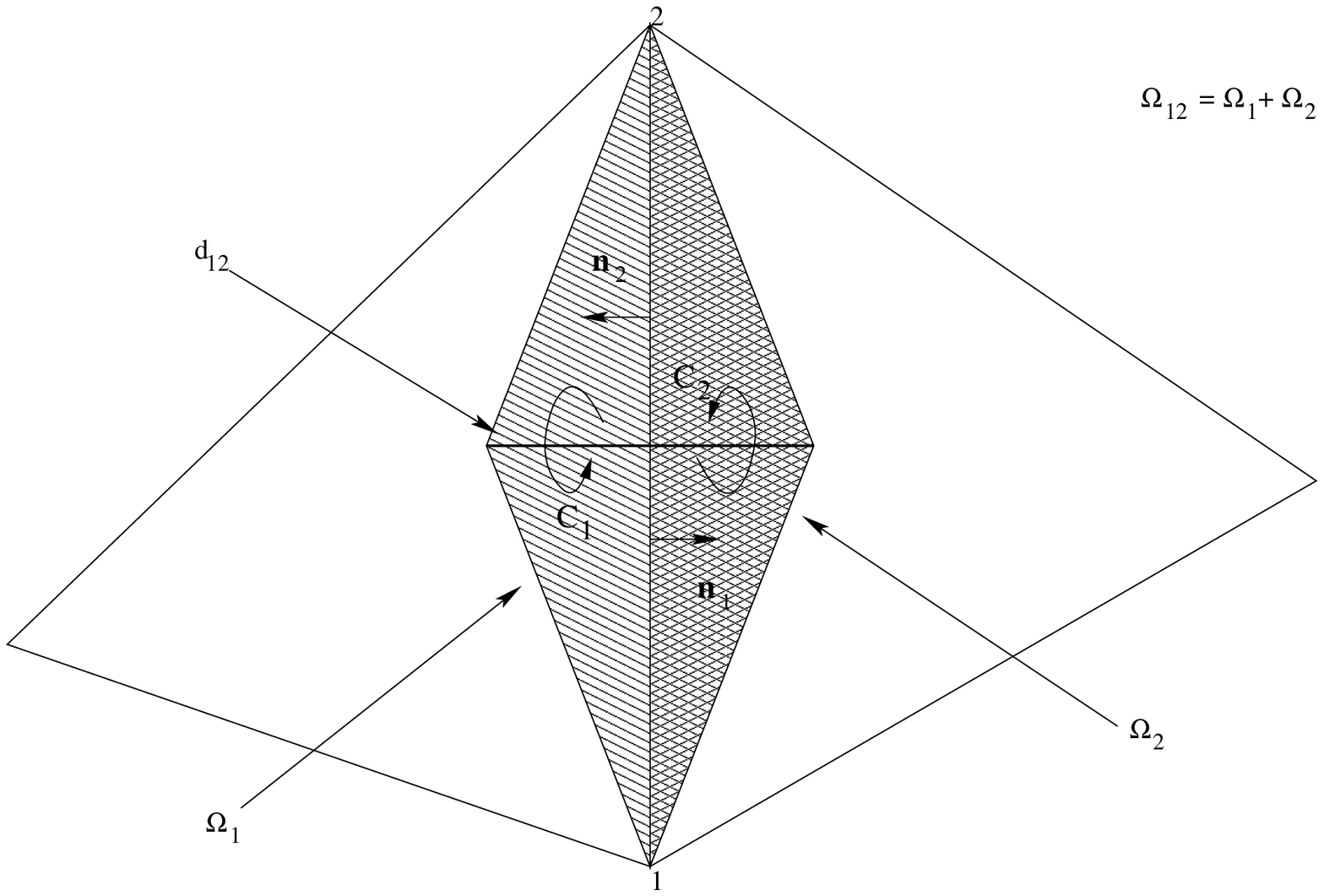}
                 \caption{}
         \label{2tri}
\end{figure}

To derive the scheme, let us first suppose that the unknown $\psi$ is constant on each triangular element.
The unknown is positioned at the circumcenter of the element.  We seek the gradient in the normal direction
across any two adjacent elements.  Let the shaded region in figure \ref{2tri} represent the control volume.  $C_1$ will be
the path along the shaded region of the element to the left and $C_2$ that of the shaded region of the element
to the right.
Integrating over this surface we have,
\begin{eqnarray}\int_\Omega\psi_x\,dx dy &=&\oint\limits_C \psi\,dy\nonumber\\
& = &\oint\limits_{C_2}\psi_2\,dy-\int_{y_2}^{y_1}\psi_2\,dy+\oint\limits_{C_1}\psi_1\,dy
-\int_{y_1}^{y_2}\psi_1\,dy\nonumber\\ 
&=&(\psi_2-\psi_1)(y_2-y_1)\nonumber\\
&\equiv & \psi_{x_{12}}\Omega\nonumber\\
\Longrightarrow \psi_{x_{12}} &=& \frac{(\psi_2-\psi_1)(y_2-y_1)}{\Omega}\label{1}\end{eqnarray} where the above line
integrals on $C_2$ and $C_1$ are zero.
Similarly we can show,
\begin{equation}\psi_{y_{12}} =-\frac{(\psi_2-\psi_1)(x_2-x_1)}{\Omega}\label{2}\end{equation}
From \eqref{1} and \eqref{2} we get,
\begin{equation}\vec{\nabla} \psi_{12}=\frac{(\psi_2-\psi_1)\mathbf n_1l_{12}}{\Omega}\end{equation} 
where $l_{12}$ is the distance between points 1 and 2.  The consistency condition is 
$\Omega$ must be chosen such that $(\vec{\nabla} \psi_{12})\cdot \mathbf d_{12}=\psi_2-\psi_1$ where $\mathbf d_{12}=d_{12}\mathbf n_1$ and
$d_{12}$ is the distance between the two circumcenters.  This
implies that,
\begin{equation}\frac{(\psi_2-\psi_1)l_{12}d_{12}}{\Omega}=\psi_2-\psi_1 \Longrightarrow \Omega=
l_{12}d_{12}=2\Omega_{12}\end{equation}  Therefore,
\begin{equation}\vec{\nabla} \psi_{12}=\frac{(\psi_2-\psi_1)\mathbf n_1}{d_{12}}\end{equation}
\begin{figure}[h]
        \centering
         \includegraphics[totalheight=3in]{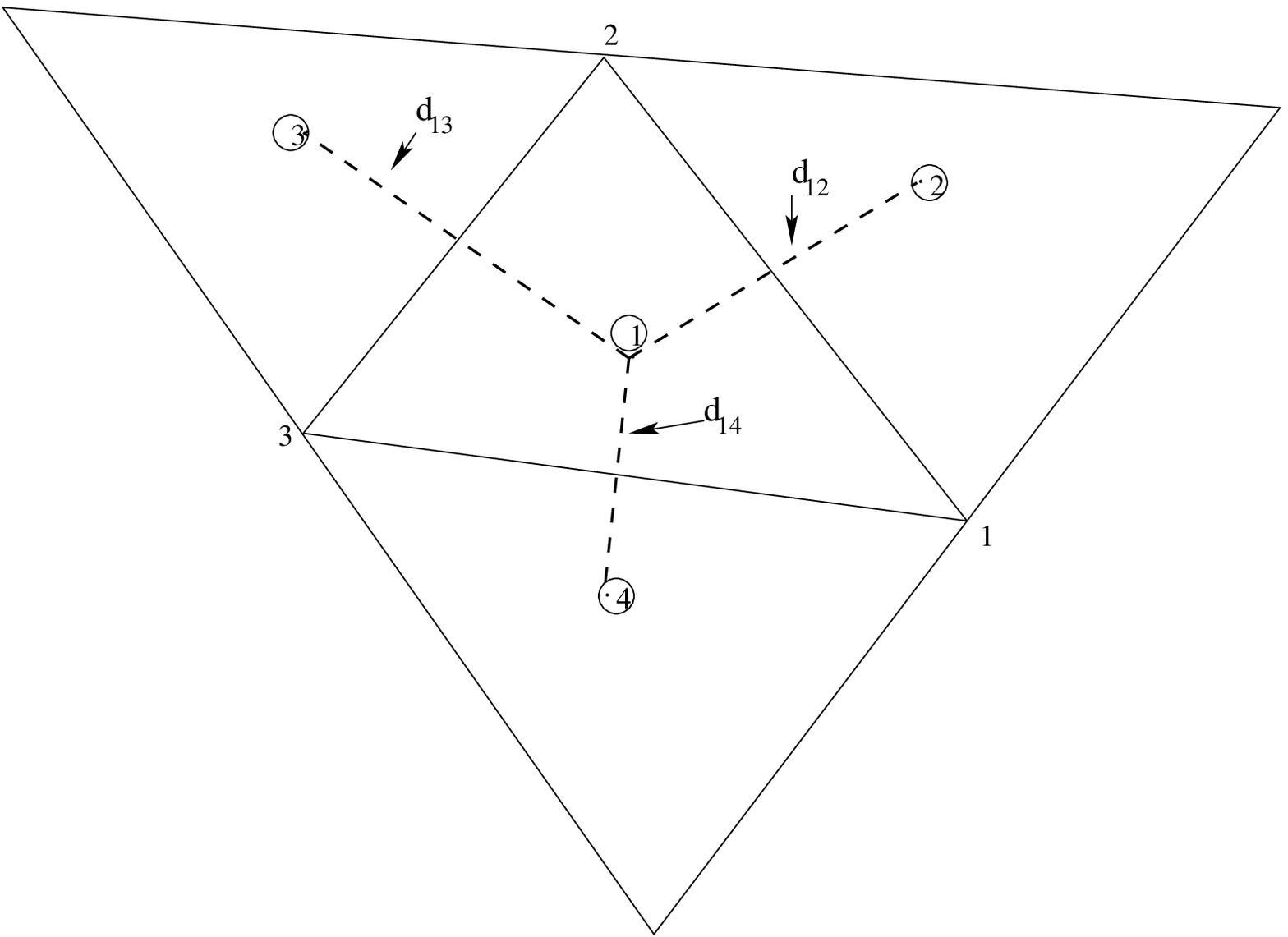}
                 \caption{}
         \label{stencil}
\end{figure}
Now, it is straightforward to obtain the stencil for the Laplacian.  Consider a set of elements
similar to those in figure \ref{stencil} where the unknowns are again at the circumcenter of each element. 
If we integrate the Laplacian over the center element numbered 1 we have,
\begin{eqnarray}0 &=&\int_{\Omega_1}\Delta\psi\,dx dy\nonumber\\ &=&\oint\limits_{C_{123}} \vec{\nabla}\psi\cdot\mathbf n\,dl\nonumber\\
&=&\frac{(\psi_2-\psi_1)\mathbf n_{12}}{d_{12}}\cdot\mathbf n_{12}l_{12}+\frac{(\psi_3-\psi_1)\mathbf n_{13}}{d_{13}}\cdot\mathbf n_{13}l_{23}+\frac{(\psi_4-\psi_1)\mathbf n_{14}}{d_{14}}\cdot\mathbf n_{14}l_{31}\nonumber\\
&=&\frac{\psi_4\,l_{31}}{d_{14}}+\frac{\psi_3\,l_{23}}{d_{13}}+\frac{\psi_2\,l_{12}}{d_{12}}-
\psi_1(\frac{l_{31}}{d_{14}}+\frac{l_{23}}{d_{13}}+\frac{l_{12}}{d_{12}})\label{6}\end{eqnarray}
The local truncation error at the circumcenter of element 1 is defined as
\begin{equation}\label{trunc}
e_{_\Delta}=\overline{\Delta\phi_1}=\frac{\int_{\Omega_1}\Delta\phi\,dx dy}{\Omega_1}
\end{equation}
\section{Numerical Verification of the Formal Accuracy}
If the elements in figure \ref{stencil} are equilateral triangles, a Taylor expansion about the circumcenter of the
first element shows that equation \eqref{6} is formally first order accurate (see appendix A).  As long as all of the
elements contain a circumcenter (all angles $<\frac{\pi}{2}$) the accuracy holds as is shown in \cite{herbin}.  When
the circumcenter is not available (i.e. right triangles) the centroid is used instead.  

\paragraph{}
Care should be exercised when verifying the accuracy of cell-centered type methods.  This is due to the fact that an ordinary refinement of the grid does not preserve the orientation of the computational stencil and the centers from the previous grid may or may not
be available.  

\begin{figure}[htbp]
\begin{minipage}[t]{.5\linewidth}
\centering
\includegraphics[width=2.8in]{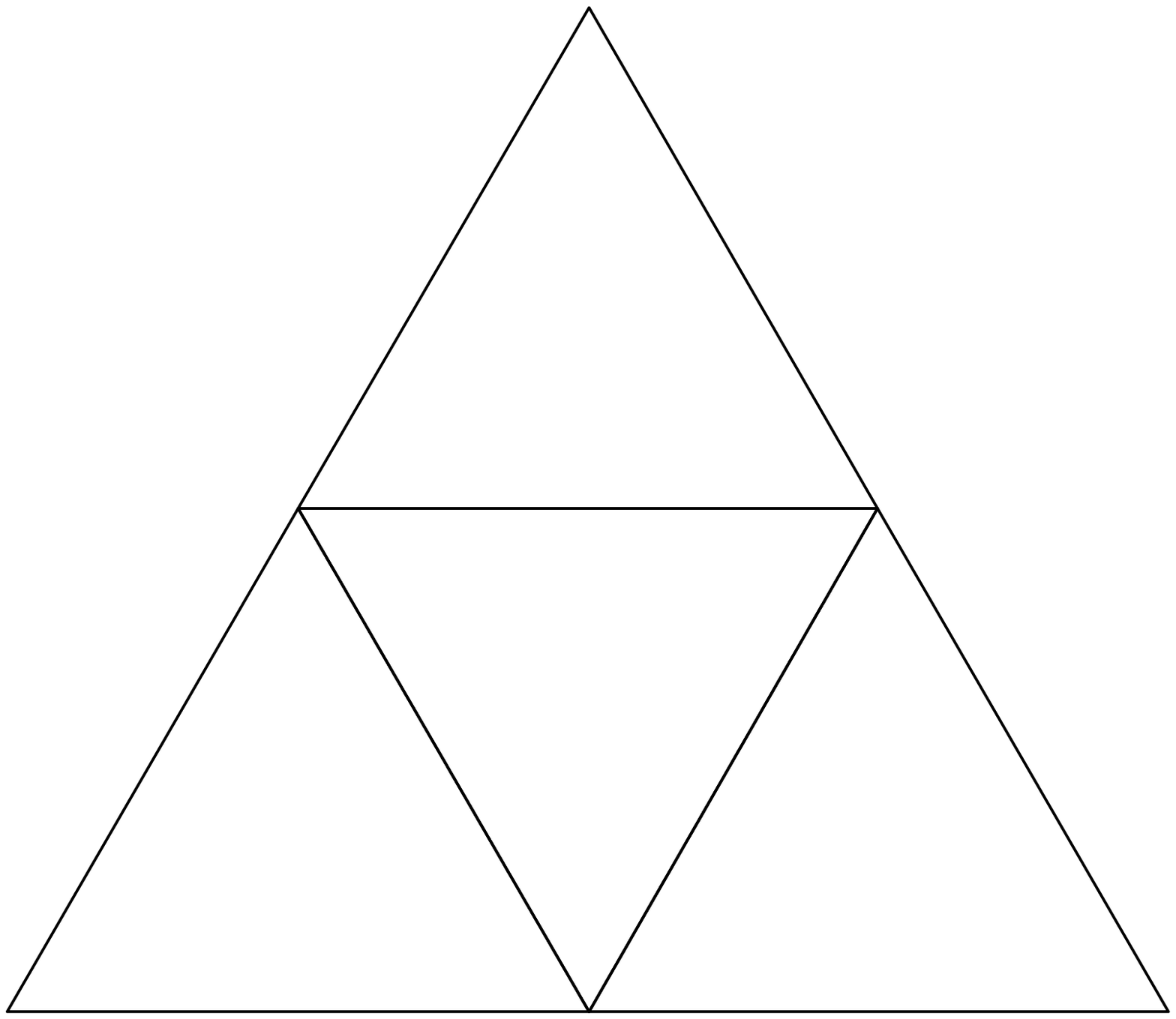}
\caption{First Grid\label{firstequil}} 
\end{minipage}%
\begin{minipage}[t]{.5\linewidth}
\centering
\includegraphics[width=2.8in]{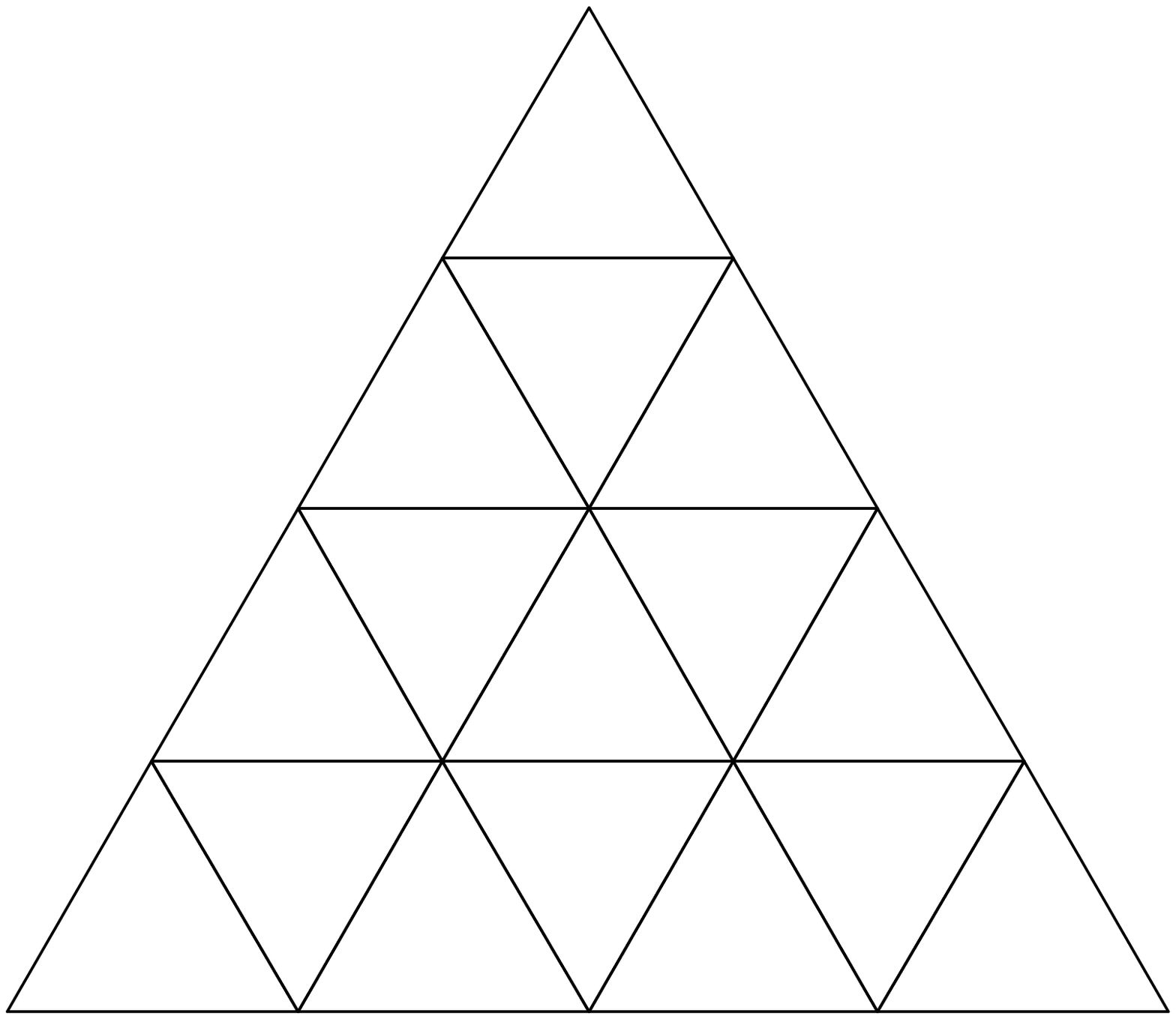}
\caption{Refined Grid\label{secondequil}} 
\end{minipage}
\end{figure}

\begin{figure}[htbp]
\begin{minipage}[t]{.5\linewidth}
\centering
\includegraphics[width=2.8in]{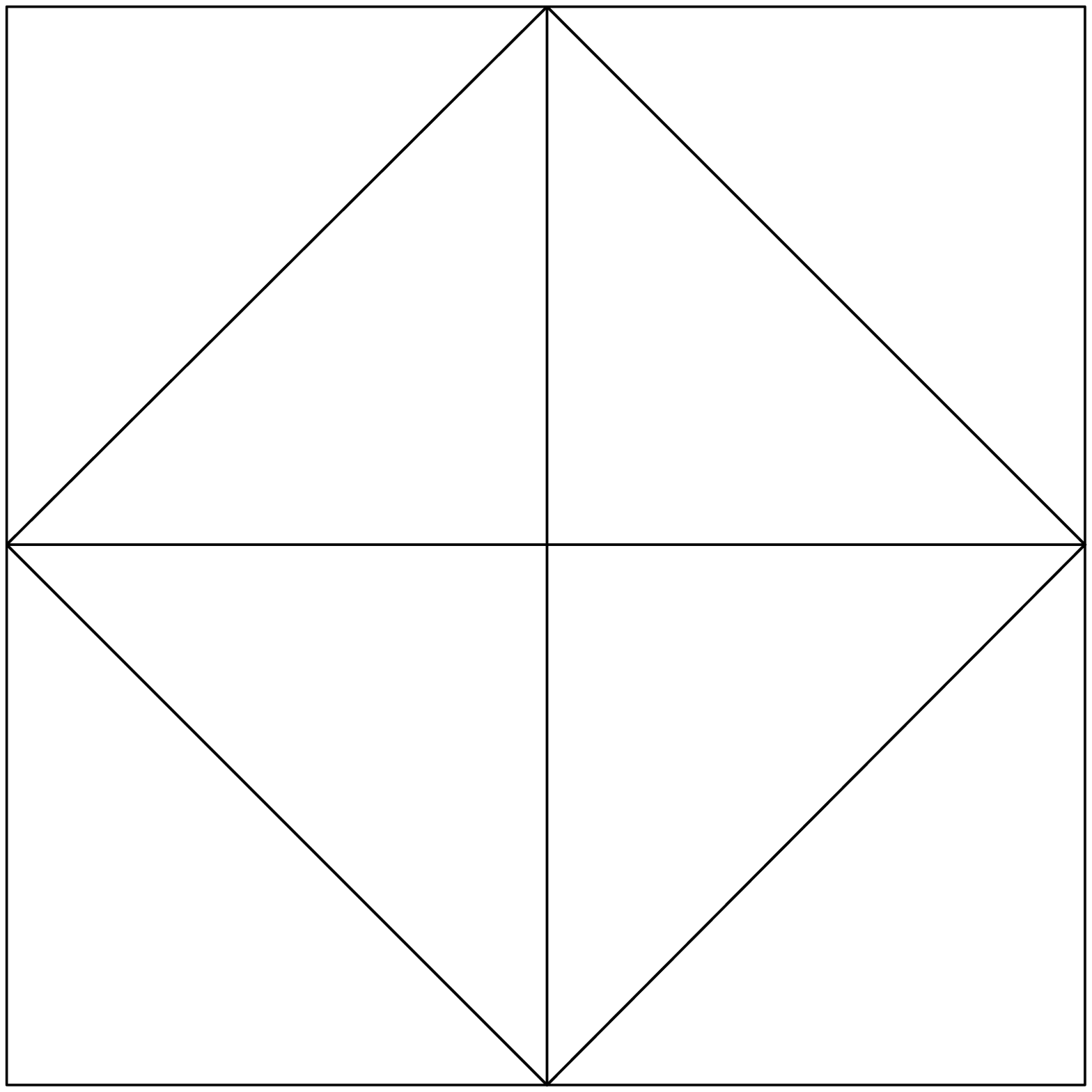}
\caption{First Grid}\label{firstsq} 
\end{minipage}%
\begin{minipage}[t]{.5\linewidth}
\centering
\includegraphics[width=2.8in]{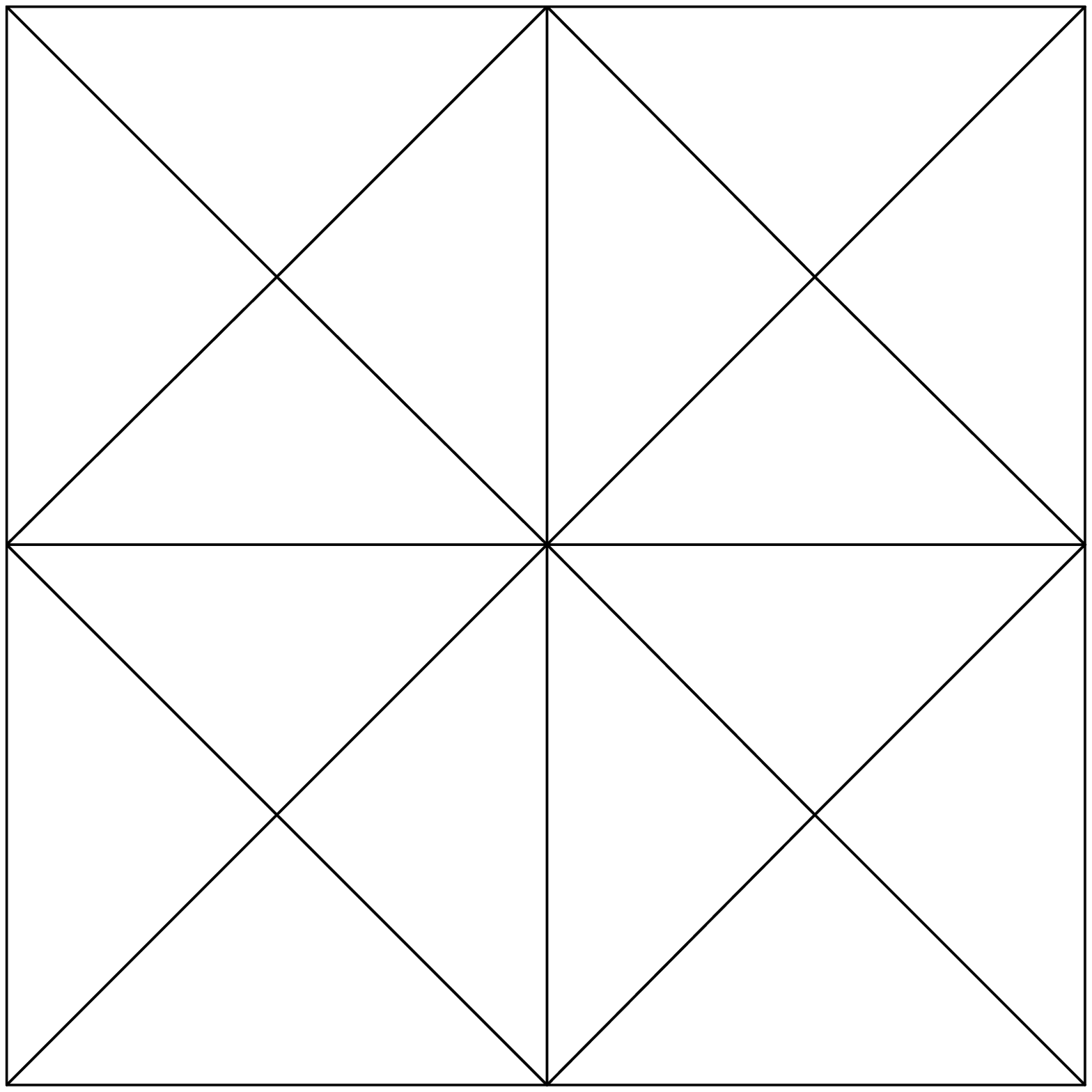}
\caption{Refined Grid\label{secondsq}} 
\end{minipage}
\end{figure}
\begin{figure}[htbp]
\centering
\includegraphics[totalheight=3in]{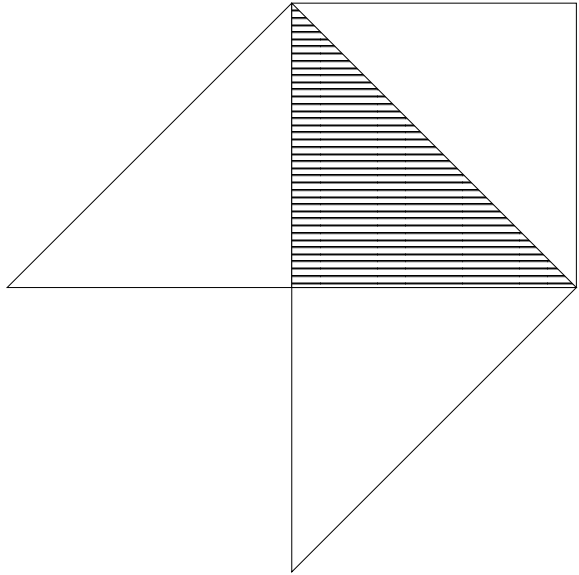}
\caption{Grid For Error Analysis On Right Triangles\label{compd}}
\end{figure}
\begin{figure}[htbp]
\centering
\includegraphics[width=4in]{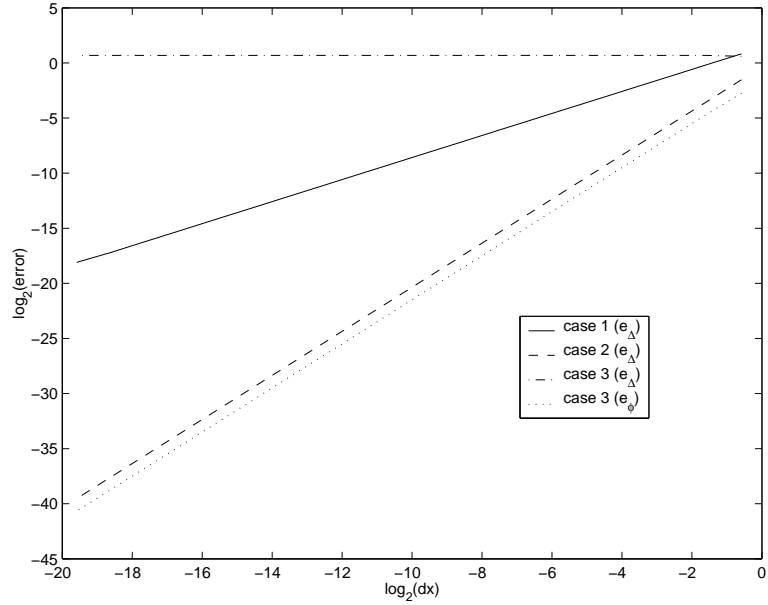}
\caption{Errors corresponding to a test patch of right triangles\label{er2}} 
\end{figure}

To clarify, observe the circumcenter of the center element in figure \ref{firstequil}.  With respect to equation (\ref{6}), the stencil for the unknown at this point
has undergone a ${180^\circ}$ rotation in the refined grid of figure \ref{secondequil}. Even worse, consider the case of right triangles where the unknowns are at the centroids.  The centroids in figure \ref{firstsq} aren't centroids of any elements in figure  \ref{secondsq} and in fact are along the edges of the triangles in the refined grid.  Because of these possible problems the observed local truncation error of a cell-centered method in general will not vary monotonically with mesh refinements.  Moreover, it can be expected to vary erratically yet, with a trend that will depend on the order of the method.
\paragraph{}
A way to avoid this problem is to use an approach that is computationally efficient and can be used for cell-centered and cell-vertex methods alike.  As outlined in \cite{chattot}, rather than refining a mesh, create one computational patch.  That is, create a triangulation with the least amount of elements that will allow the method to be fully implemented with respect to one unknown within the domain while enforcing a Dirichlet boundary condition.  The truncation error is then computed 
on this grid and the grid is rescaled by a factor of $\frac{1}{2}$ about the coordinates of the unknown.  The truncation error is then recomputed and compared to that of the previous
grid.  This process is repeated until the formal accuracy of the method is apparent or machine zero is reached.  This technique is used on the grid illustrated in figure \ref{firstequil} and figure \ref{compd} where the shaded triangle's centroid is given the coordinates $(0,0)$.
\paragraph{}
Here, the order of the truncation error $e_{_\Delta}$ is demonstrated for the following problem:
\begin{equation}
\Delta\phi=\phi_{xx} + \phi_{yy}= 0 
\end{equation}
where the following cases are examined:
\begin{enumerate}
\item $\phi=(x-1)^3 -3(x-1)y^2$
\item $\phi=x^4+y^4-6x^2y^2$
\item $\phi=\frac{1}{\sinh(\pi)}(\sinh(\pi x)\sin(\pi y)+\sinh(\pi y)\sin(\pi x))$ 
\end{enumerate}
The circumcenter of an equilateral triangle is equivalent to it's circumcenter.  Therefore, the variables are placed at the centroids for both equilateral and right triangles.  For this analysis the error of the solution is defined as 
\begin{equation}
e_{_\phi}=|\phi - \phi_{exact}|_{centroid}
\end{equation}
and $dx$ is
taken as the distance from the centroid of the unknown to the closest bordering centriod in the grid. 
\paragraph{}
In figure (\ref{er1}) for a grid of equilateral triangles, the method is exact for case 1
(see Appendix A) and therefore is not included.  The slope of the line for case 2 is exactly 1 and for case 3 it has a limiting slope of 1.  Numerically,
this confirms the formal accuracy of the method on a grid consisting entirely of equilateral triangles. 
\begin{figure}[htbp]
\centering
\includegraphics[width=4in]{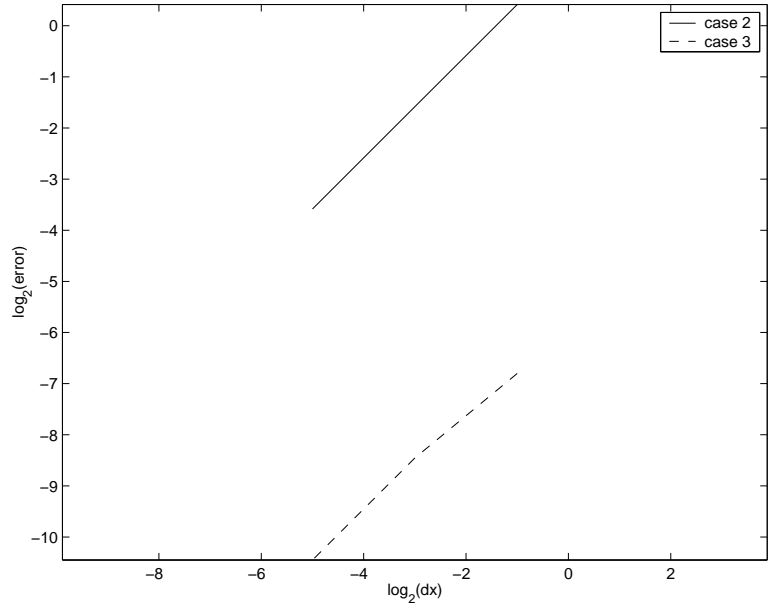}
\caption{Errors corresponding to a test patch of equilateral triangles\label{er1}}
\end{figure} 
\paragraph{}
In figure (\ref{er2}) for the grid consisting of right triangles, the method shows
first order behavior for case 1 but shows second order accuracy for case 2.  This is surprising since in general the method is not pointwise
consistent on a grid consisting of right triangles.  This means that the Laplacian from equation (\ref{6}) will not be consistent when the variables are placed at the centroids.  But, for these particular test cases a Taylor expansion of equation (\ref{6}) proves consistent with what is observed numerically. The third case is a true test because the Taylor expansion contains derivatives to all orders.  Note that the truncation error $e_{_\Delta}$  is zeroth order.  However,
the graph of the error $e_{_\phi}$ for that case is second order accurate.  This illustrates an important point
as discussed in \cite{chattot}.  Namely, a lack of consistency does not preclude convergence.  Second order accuracy of $e_{_\phi}$
in general cannot be guaranteed on arbitrary triangulations.  In \cite{lazarov} second order accuracy is discussed for triangulations with certain regularities and symmetries.

\section{Inviscid Incompressible Flow Over An Airfoil}
Let us apply this method to flow over a symmetric airfoil at an angle of attack.  Though the methodology applies to an
arbitrary airfoil, here it is restricted to the class of airfoils that can be generated using the K\'{a}rm\'{a}n-Trefftz conformal transformation (see \cite{vaz}).  For a stream function $\psi=\psi(x,y)$, the
governing equation is \begin{equation}\psi_{xx} + \psi_{yy} = 0 \end{equation}

\subsection{Numerical Boundary Conditions}
\begin{figure}[htbp]
        \centering
                \includegraphics[totalheight=3in]{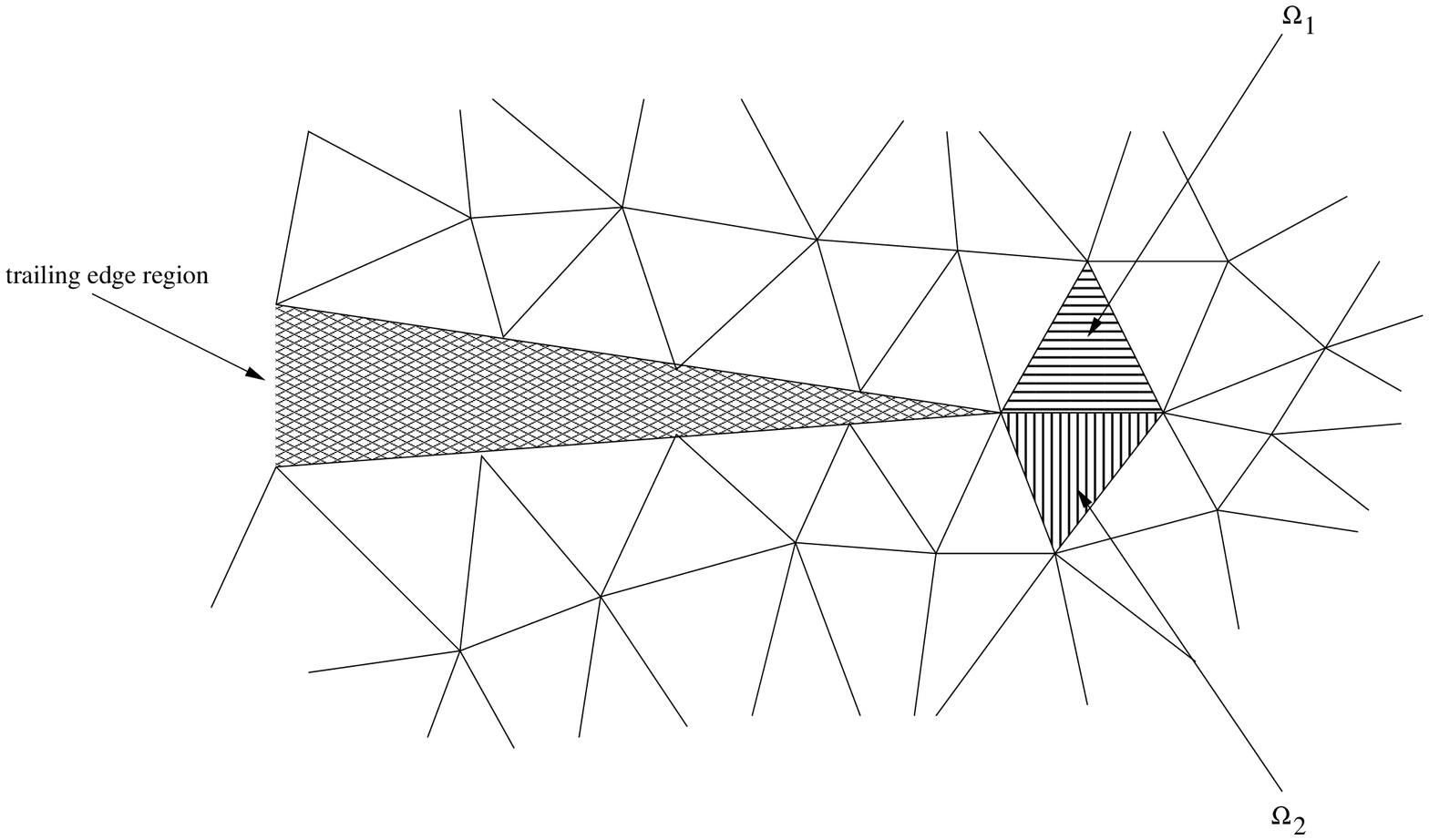}
                 \caption{Trailing Edge\label{kutta}}
\end{figure}

\begin{figure}[htbp]
        \centering
                \includegraphics[totalheight=3in]{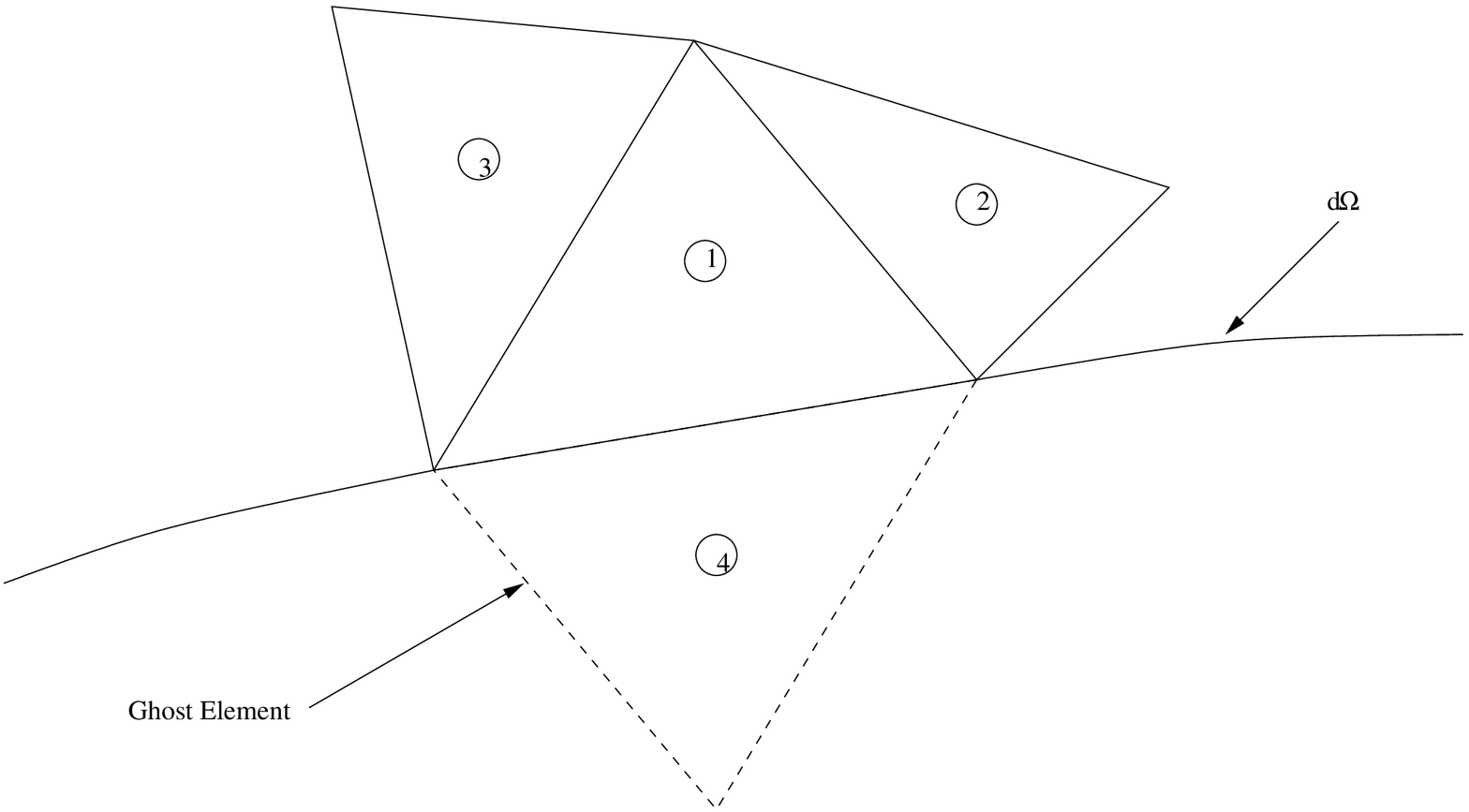}
                 \caption{\label{ghost}}
\end{figure}

In the far field we enforce uniform flow on all elements with a node on the farfield boundary.  Namely, 
\begin{equation}
\psi=y \cos{\alpha}- x \sin{\alpha}+\frac{\Gamma \ln(r)}{2 \pi}\label{streamfunction}
\end{equation} where
\begin{equation}
\Gamma=4\pi a \sin{\alpha}
\end{equation} 
and $a$ is the radius of the cylinder used to create the airfoil via the K\'{a}rm\'{a}n-Trefftz transformation.
  The angle of attack is denoted by $\alpha$ and is taken to be $15^{\circ}$.  
\paragraph{}
On the surface, the value of the streamfunction is a constant but is not known beforehand.  It is part of the solution and must be 
calculated iteratively, effectively giving a Kutta condition.  This is done by first taking the average
value of the two elements immediately downstream of the trailing edge that are parallel to the camber line (see figure \ref{kutta}).  
This
average is then placed along the surface of the airfoil and updated after each iteration through the system
of unknowns.  For each element with an edge on the boundary, a ghost element is created by reflection about
this edge.  Now, the value at the ghost element is given a value such that the average of the element and the
ghost element equals the value on the surface. 
To clarify, imagine element 4 is a ghost element of element 1 in
figure \ref{ghost}.  Then, 
\[\psi^*=\frac{\psi_1+\psi_4}{2}\] where $\psi^*$ is the value on the surface.  
Hence,
\[\psi_4=2\psi^*-\psi_1\]
Now, the stencil we get from equation \eqref{6} can be used for all elements on the surface.
\section{Results}  
\begin{figure}[htbp]
        \centering
         \includegraphics[totalheight=3in]{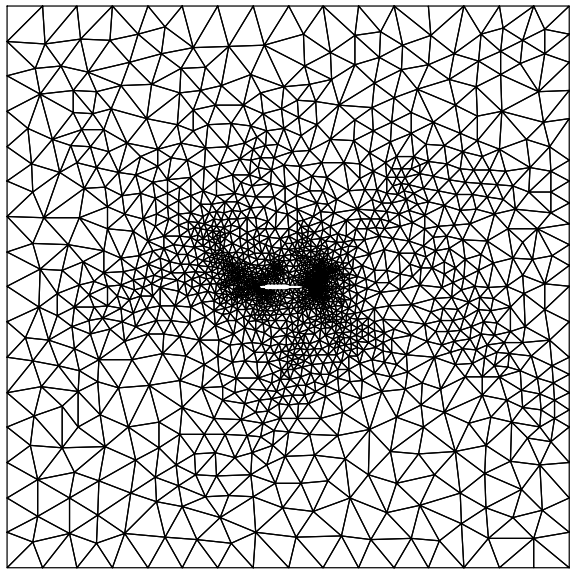}
                 \caption{Unstructured Grid About A Joukowski Airfoil \label{airgrid}}
\end{figure}

\begin{figure}[htbp]
        \centering
         \includegraphics[totalheight=2.5in,width=2.5in]{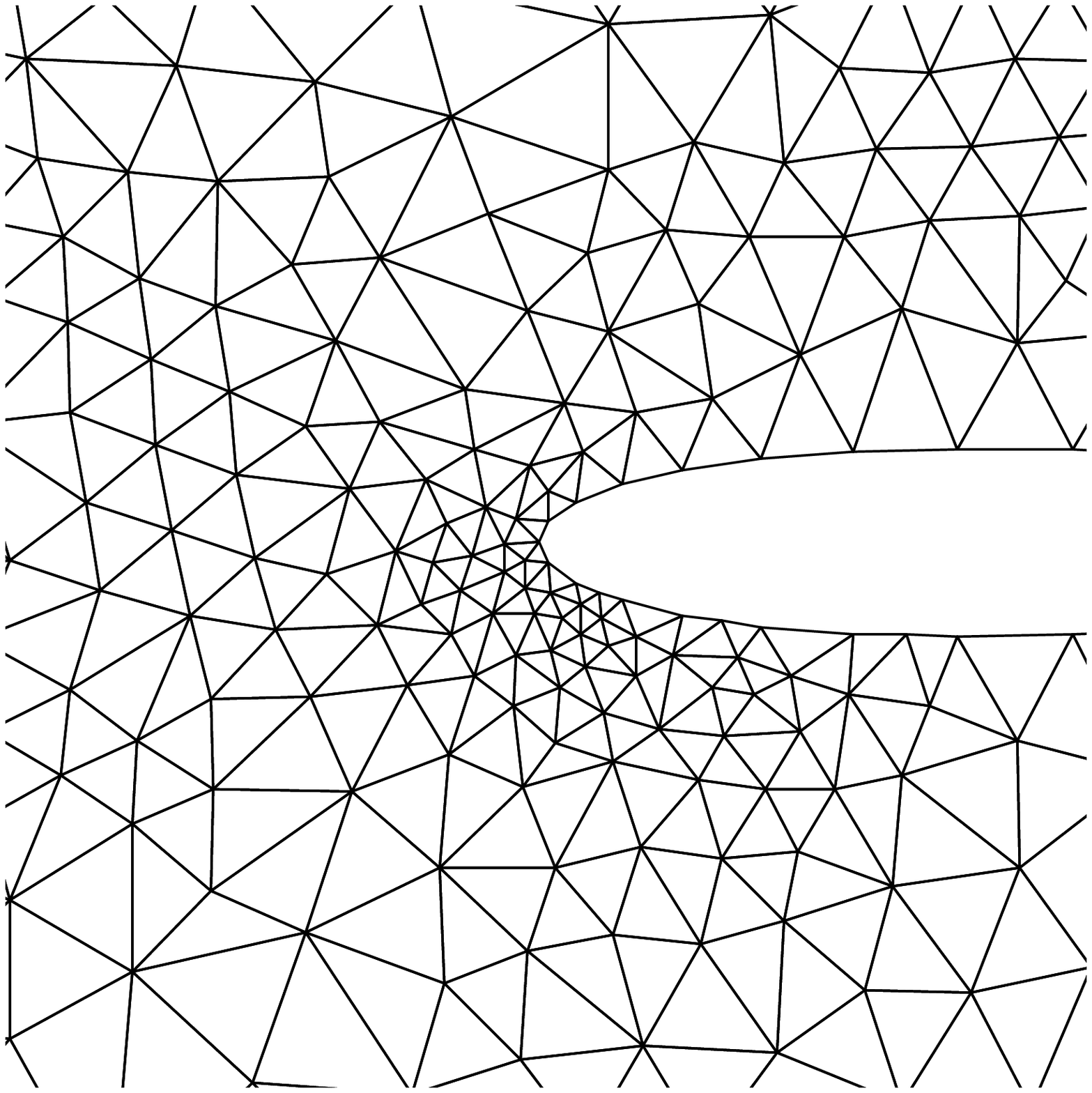}%
                \hspace{.25in}
                 \includegraphics[totalheight=2.5in,width=2.5in]{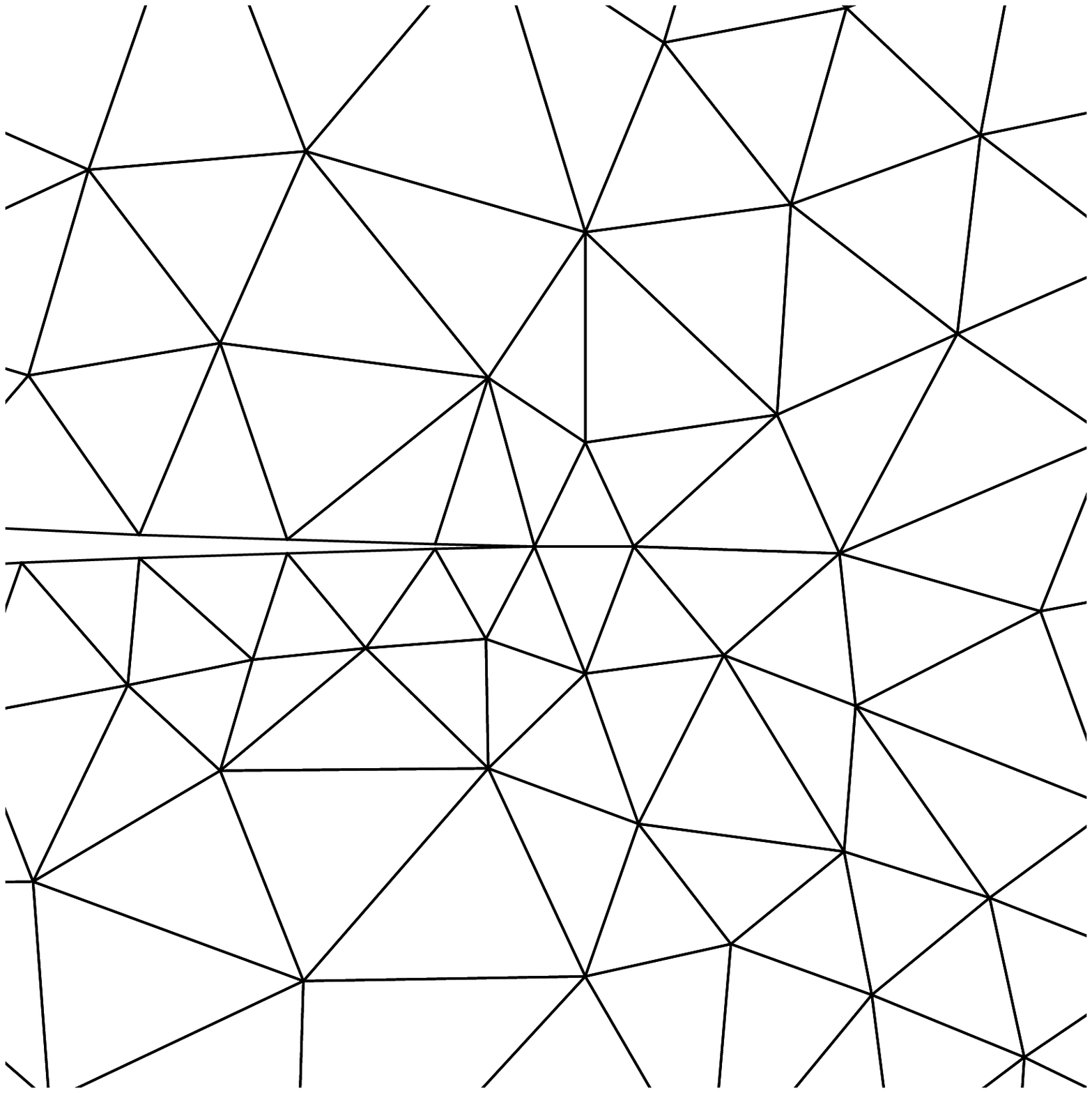}
                 \caption{Close-up of Leading and Trailing Edge\label{close}}
\end{figure}

\begin{figure}[htbp]
        \centering
                \includegraphics[totalheight=3.7in]{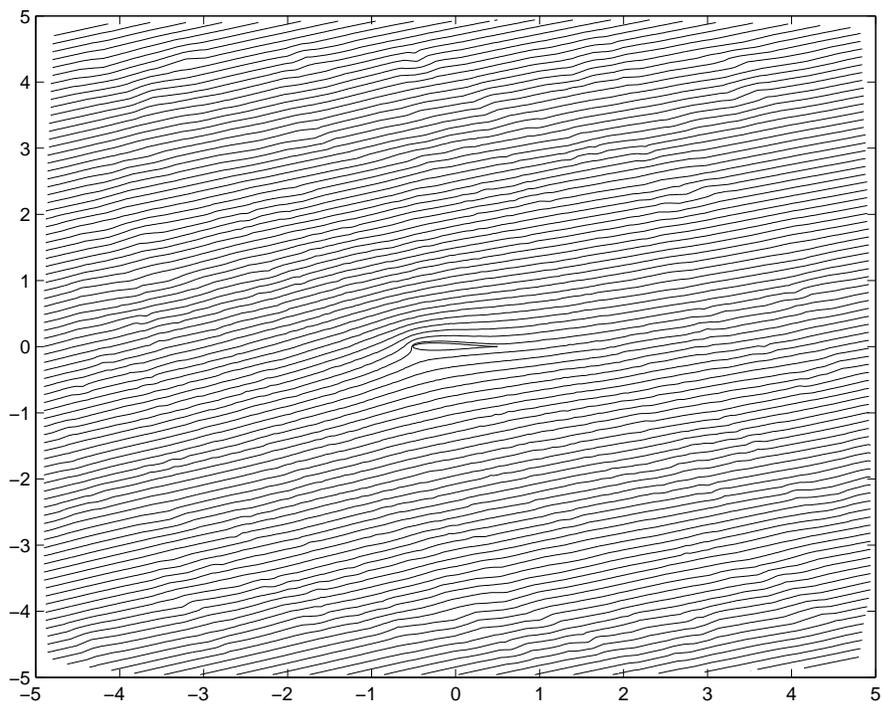}
                 \caption{Streamlines About A Joukowski Airfoil (173 surface elements)\label{jouk40}}
\end{figure}

The stream function was calculated for the flow over a symmetric Joukowski airfoil with a chord length of one and a twelve percent thickness ratio. Typical grids about the airfoil are shown in figures \ref{airgrid}  and \ref{close} were the computational domain is a
$[-5,5]\times[-5,5]$ box\footnotemark[1].  The solutions were obtained using the Gauss-Seidel iterative
method. The iterations were terminated once the value of $\psi^*$ reached a plateau. As a test of convergence, the
value of the circulation around the airfoil is compared to the exact value $\Gamma= 0.888215341$ for 3 grids with varying numbers of elements on the surface of the airfoil.
The results are shown in table \ref{table0}. The method captures the correct solution throughout the
domain. A contour plot of the solution to the stream function is shown in figure \ref{jouk40} for
a Joukowski airfoil with 173 surface elements.\\ 

\footnotetext[1]{All grids produced with the aid of the ``triangle'' mesh generator \cite{shewchuk}.}
\begin{table}[H]\centering \label{table0}
\begin{tabular}{|c|c|} \hline
surface points & $\Gamma$\\ 
\hline \hline
52& 0.945274938\\ 
98& 0.905162059\\
173&0.880039787\\ 
 \hline 
\end{tabular}
\caption{Error table}
\end{table}
\begin{figure}[htbp]
        \centering
                \includegraphics[totalheight=3.7in]{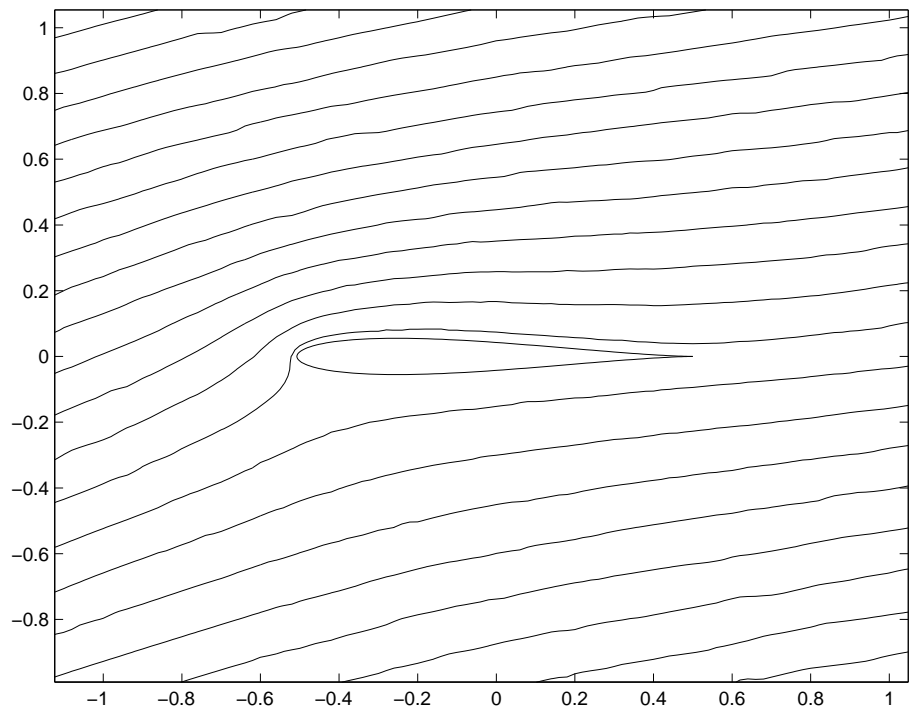}
                 \caption{Close-up About A Joukowski Airfoil (173 surface elements)\label{jouk41}}
\end{figure}

\chapter{A Least-Squares Finite Volume Method}
In light of the restrictions of the method derived in the previous chapter, we seek a higher order method that
also allows greater flexibility in the type of triangulation used.  A way to achieve this is by placing the unknowns
of the governing equations at the vertices of the triangle and assume a linear variation over each element.  
\section{Preliminary Derivations For A Scalar Equation}
\begin{figure}[htbp]
        \centering
         \includegraphics[totalheight=3in]{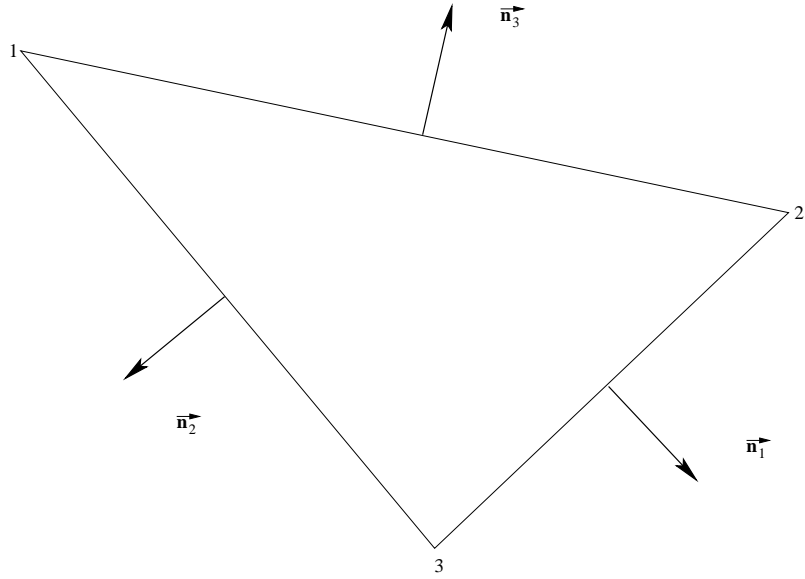}
         \caption{linear element \label{element}}
\end{figure}

As mentioned earlier, the governing equations are approximated over each element using a finite
volume approximation.  Figure \ref{element} represents a typical triangular element for the method developed here.
The variables are placed at the vertices which gives a linear approximation over each element.  Consider
an arbitrary variable $\phi=\phi(x,y)$ given at the three vertices.  The gradient is given by:

\begin{eqnarray*}\int\limits_{\Omega}\vec{\nabla}\phi\,d\Omega=\oint\limits_{\partial\Omega}\phi\vec{n}\,dl
\end{eqnarray*}
Since $\phi$ varies linearly over the element,
\begin{eqnarray*}\vec{\nabla}\phi\Omega&=&\oint\limits_{\partial\Omega}\phi\vec{n}\,dl
\end{eqnarray*}

Using the trapeziodal rule on each side,
\begin{eqnarray*}\vec{\nabla}\phi\Omega&\approx&\frac{\phi_2+\phi_1}{2}\vec{n}_3
+\frac{\phi_3+\phi_2}{2}\vec{n}_1+\frac{\phi_1+\phi_3}{2}\vec{n}_2
\end{eqnarray*}
where the normal vectors have been scaled to the length of the corresponding edge.  Since for any triangle
$\vec{n}_1+\vec{n}_2+\vec{n}_3=\vec{0}$, we have,
\begin{eqnarray*}
\vec{\nabla}\phi\Omega&=&\frac{-(\phi_1\vec{n}_1+\phi_2\vec{n}_2+\phi_3\vec{n}_3)}{2}
\end{eqnarray*}
It follows that the approximation for the gradient of $\phi$ 
in the triangle is,
\begin{eqnarray} 
\vec{\nabla}\phi&=&\frac{-1}{2\Omega}\sum_{i=1}^{3}\phi_i\vec{n}_i\label{gradient}
\end{eqnarray} 
Consequently,
\begin{eqnarray} 
\phi_x&=&\frac{-1}{2\Omega}\sum_{i=1}^{3}\phi_in_{x_i}\label{dx}\\
\phi_y&=&\frac{-1}{2\Omega}\sum_{i=1}^{3}\phi_in_{y_i}\label{dy}
\end{eqnarray} 
Note that since $\phi$ is linear on an element, $\vec{\nabla} \phi$ is just a constant vector.
Now, consider the equation
\begin{equation} 
\vec{c}\cdot\vec{\nabla} \phi=0\label{gov1}
\end{equation}
where $\vec{c} = (a,b) = (f(x,y,\phi),g(x,y,\phi))$.  Integrating (\ref{gov1}) over an element $\Omega_{_T}$ 
we have,
\begin{eqnarray}
\int\limits_{\Omega_{_T}}\vec{c}\cdot\vec{\nabla}\phi\,d\Omega&=&\left(\int\limits_{\Omega_{_T}}\vec{c}\,d\Omega\right)\cdot\vec{\nabla}\phi\nonumber\\
                                                              &=&\overline{\vec{c}}\cdot\vec{\nabla}\phi\,\Omega_{_T}\label{dum}
\end{eqnarray}
Where, 
\begin{equation}
\overline{\vec{c}}=\frac{1}{\Omega_{_T}}\int\limits_{\Omega_{_T}}\vec{c}\,d\Omega\label{ave}
\end{equation}
In general, it can be complicated to evaluate equation (\ref{ave}) exactly.  However, for the equations approximated in this work, $\vec{c}$ varies at most linearly with respect to $x,y$ and $\phi$.
Henceforth, assuming $\vec{c}$ varies linearly
on an element we get,
\begin{equation}
\overline{\vec{c}}=\frac{1}{\Omega_{_T}}\int\limits_{\Omega_{_T}}\vec{c}\,d\Omega=\frac{1}{3}\sum_{j=1}^{3}\vec{c}_j\label{ave1}
\end{equation}
Combining the results of (\ref{gradient}), (\ref{dum}) and (\ref{ave1}) we have,
\begin{equation}
\int\limits_{\Omega_{_T}}\vec{c}\cdot\vec{\nabla}\phi\,d\Omega=\frac{-1}{2}\sum_{i=1}^{3}\phi_i\overline{\vec{c}}\cdot\vec{n}_i=R_{_T}\Omega_{_T}\label{eq1}
\end{equation}
where $R_T$ is commonly called the cell residual.
\section{The Discrete Functional}
The goal is to apply (\ref{eq1}) on each triangle in the domain.  For any given triangulation there are always more nodes than
triangles and thus, we have more unknowns than available equations.  To circumvent this inadequacy, one might try supplementing the
resulting underdetermined system of equations with additional equations that conform to the physics of the problem and yield a unique
solution \cite{alcrudo}.  In general, the process of determining the correct additional equations is complicated.  In this work a least-squares 
method is chosen instead.
\paragraph{} So, let the functional $\mathcal{I}$ be represented by,
\begin{equation}
\mathcal{I}=\frac{1}{2}\sum_TR_{_T}^2\Omega_{_T}\label{func}
\end{equation}
where the sum is over all of the triangles in the domain.  This form of $\mathcal{I}$ is by no means obligatory but as mentioned in 
\cite{baines:mov}, has the nice property that if $\vec{c}$ is a constant vector, then $\mathcal{I}$ is in fact twice the square of 
the $L_2$ norm of $R_T=\vec{c}\cdot\vec{\nabla}\phi$:
\begin{equation}
\begin{split}
\int\limits_{\Omega}(\vec{c}\cdot\vec{\nabla}\phi)^2\,d\Omega=\sum_T\int\limits_{\Omega_{_T}}(\vec{c}\cdot\vec{\nabla}\phi)^2\,d\Omega
&=\sum_T(\vec{c}\cdot\vec{\nabla}\phi)^2_{_T}\int\limits_{\Omega_{_T}}d\Omega\\
&=\sum_T(\vec{c}\cdot\vec{\nabla}\phi)^2_{_T}\Omega_{_T}\\
&=\sum_TR_{_T}^2\Omega_{_T}\\
&=2\mathcal{I}\end{split}
\end{equation}

\begin{figure}[htbp]
        \centering
         \includegraphics[totalheight=3in]{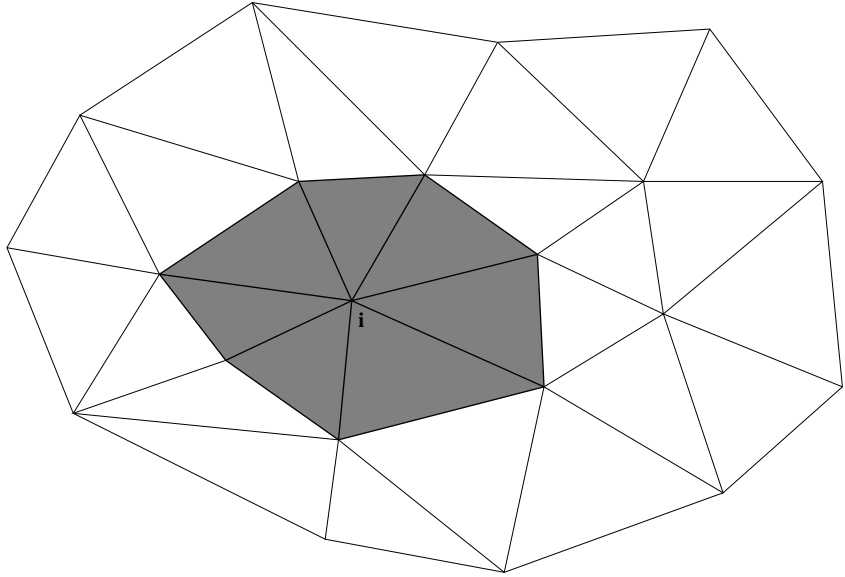}
         \caption{$\Omega{_i}$ of node $i$ (shaded region) \label{dual}}
\end{figure}

\paragraph{}
In \cite{darcha} the authors constructed a functional by calculating the residuals of all elements ajoining a particular node 
(see figure (\ref{dual})).  For such a functional $\mathcal{I}^*$ we have,
\begin{equation}
\mathcal{I}^*=\frac{1}{2}\sum_{i=1}^N R_{_i}^2\Omega_{_i}=\frac{1}{2}\sum_{i=1}^N \sum_{T_i}R_{_{T_i}}^2\Omega_{_{T_i}}
\end{equation}
where $N$ is the total number of nodes in the domain.  Each element gives three contributions to the
double summation of the right hand side.  Therefore we get,
\begin{equation}
\mathcal{I}^*=\frac{1}{2}\sum_{i=1}^N \sum_{T_i}R_{_{T_i}}^2\Omega_{_{T_i}}=\frac{3}{2}\sum_TR_T^2\Omega_{_T}=3\mathcal{I}
\end{equation}
So we find that either approach to building up the functional will give the same minimization though one can use his discretion
in deciding which approach favors the intended method of minimization and numerical implementation.

\subsection{Minimization}
We seek a minimization of equation (\ref{func}) with respect to the variable $\phi$.  Let
\begin{equation}
\mathcal{I}_{_T}=\frac{1}{2}R_{_T}^2\Omega_{_T}
\end{equation}  Then,
\begin{equation}
\mathcal{I}=\sum_T\mathcal{I}_{_T} 
\end{equation} 
and
\begin{equation}
\frac{\partial\mathcal{I}}{\partial\phi_i}=\sum_T\frac{\partial\mathcal{I}_{_T}}{\partial\phi_i} 
\end{equation} 
On each element $\Omega_{_T}$ containing the node indexed $i$ we have,
\begin{equation}
\frac{\partial\mathcal{I}_{_T}}{\partial\phi_i}=R_{_T}\Omega_{_T}\frac{\partial R_{_T}}{\partial\phi_i}
=R_{_T}\Omega_{_T}\frac{\partial}{\partial\phi_i}(\overline{\vec{c}}\cdot\vec{\nabla}\phi)
=R_{_T}\Omega_{_T}\left(\frac{\partial \overline{\vec{c}}}{\partial\phi_i}\cdot\vec{\nabla}\phi-\frac{1}{2\Omega_{_T}}
\overline{\vec{c}}\cdot\vec{n}_i\right)\label{mini}
\end{equation} 
Collecting all of these contributions at each node permits a steepest descent method of minimization for the discrete functional.
However, in this work Newton's method was chosen as the method of minimization.  Consequently, the Jacobian of the resulting nonlinear
system of equations must be calculated.  Let $j$ index the nodes of an element $\Omega_{_T}$.  Taking the partial derivative of equation
\eqref{mini} with respect to $\phi_j$ gives
\begin{equation}
\begin{split}
\frac{\partial ^2 \mathcal{I}_{_T}}{\partial\phi_i\partial\phi_j}
&=\frac{\partial}{\partial\phi_j}\left(R_{_T}\Omega_{_T}\left(\frac{\partial \overline{\vec{c}}}{\partial\phi_i}\cdot\vec{\nabla}\phi-\frac{1}{2\Omega_{_T}}\overline{\vec{c}}\cdot\vec{n}_i\right)\right)\\
&=\left(\frac{\partial \overline{\vec{c}}}{\partial\phi_j}\cdot\vec{\nabla}\phi-\frac{1}{2\Omega_{_T}}\overline{\vec{c}}\cdot\vec{n}_j\right) 
\left(\frac{\partial \overline{\vec{c}}}{\partial\phi_i}\cdot\vec{\nabla}\phi-\frac{1}{2\Omega_{_T}}\overline{\vec{c}}\cdot\vec{n}_i\right)\Omega_{_T}\\
& \quad +\left(\frac{\partial ^2 \overline{\vec{c}}}{\partial\phi_j\partial\phi_i}\cdot\vec{\nabla}\phi
-\frac{1}{2\Omega_{_T}} \left(\frac{\partial \overline{\vec{c}}}{\partial\phi_i}\cdot\vec{n}_j+\frac{\partial \overline{\vec{c}}}{\partial\phi_j}\cdot\vec{n}_i\right)\right)R_{_T}\Omega_{_T}\label{jac}
\end{split}
\end{equation}
Note that we have the symmetry property:
\begin{equation}\label{symm}
\frac{\partial ^2 \mathcal{I}_{_T}}{\partial\phi_i\partial\phi_j}=\frac{\partial ^2 \mathcal{I}_{_T}}{\partial\phi_j\partial\phi_i}
\end{equation}
\subsubsection{Newton's Method}
Let us call $J$ the Jacobian whose entries are formed via equation \eqref{jac}.  Let
\begin{equation}
\vec{\mathcal{F}}=(\frac{\partial \mathcal{I}}{\partial \phi_1},\frac{\partial \mathcal{I}}{\partial \phi_2},\ldots,\frac{\partial \mathcal{I}}{\partial \phi_N})^T
\end{equation} Now if 
\begin{equation} 
\vec{\phi}=(\phi_1,\phi_2,\ldots,\phi_N)^T
\end{equation} 
using Newton's method, we solve $\vec{\phi}$ from
\begin{equation}
\vec{\phi}^{new}=\vec{\phi}^{old}+\vec{q}\label{correct}
\end{equation}
where $\vec{q}$ is the solution to the system:
\begin{equation}
J\vec{q}=-\vec{\mathcal{F}}
\end{equation}
The convergence criterion is either $\|\vec{q}\|$ (Euclidean norm) has to be within a desired tolerance or that the calculation of
${\mathcal{I}}$ has reached a steady value to within a desired tolerance.  Occasionally, it might be beneficial to under-relax
on the correction of \eqref{correct} in which case we have
\begin{equation}
\vec{\phi}^{new}=\vec{\phi}^{old}+\sigma\vec{q}
\end{equation}
where $0<\sigma\leq 1$.
\section{Two-Dimensional Systems Of Equations}
Extending the proposed least-squares method to a two-dimensional system of equations is straightforward.  
Let $\mathbf u =(u_1,u_2,\ldots,u_n)^t \in \Omega\subset\mathbb R ^n$ and 
$\mathbf A=\mathbf A(\mathbf u), \mathbf B=\mathbf B(\mathbf u)$ be $n\times n$
matrices allowing a two-dimensional $n\times n$ system to be written as
\begin{equation}
\mathbf A \mathbf u_x+\mathbf B \mathbf u_y=\mathbf 0
\end{equation}
Each element of $\mathbf u$ is piecewise linear over the elements in the domain.
Integrating over an element $\Omega_{_T}$,
\begin{equation}
\int\limits_{\Omega_{_T}}(\mathbf A \mathbf u_x+\mathbf B \mathbf u_y)\,d\Omega
=\vec{\nabla}\mathbf u\cdot\left(\int\limits_{\Omega_{_T}}(\mathbf A,\mathbf B)\,d\Omega\right)
\end{equation}
We assume a linear variation of $\mathbf A$, $\mathbf B$ with respect to $\mathbf u$ as in the scalar case.
Consequently we get
\begin{align}
\bar{\mathbf A}=\mathbf A(\bar{\mathbf u})\\
\bar{\mathbf B}=\mathbf B(\bar{\mathbf u})
\end{align}
where $\bar{\mathbf u}=\frac{1}{3}(\mathbf u_1 + \mathbf u_2 + \mathbf u_3)$. For the residual we have,
\begin{equation}
\int\limits_{\Omega_{_T}}(\mathbf A \mathbf u_x+\mathbf B \mathbf u_y)\,d\Omega
=(\bar{\mathbf A}\mathbf u_x+\bar{\mathbf B}\mathbf u_y)_{_T}\Omega_{_T}=\mathbf R_{_T}\Omega_{_T}
\end{equation}
Now the discrete functional is
\begin{equation}
\mathcal{I}=\frac{1}{2}\sum_T\mathbf R_T^{^t}\mathbf R_{_T}\Omega_{_T}\label{funcsys}
\end{equation}
Let 
\begin{equation}
(u,v)=\left\{(u_l,u_m)\in\mathbf u: l,m=1,2,3,\ldots,n\right\} 
\end{equation}
Consider the variation of the functional on an element $\Omega_{_T}$ with respect to $u_i$.  We have,

\begin{equation}
\begin{split}
\frac{\partial\mathcal{I}_{_T}}{\partial u_i}
=\mathbf R_{_T}^t\Omega_{_T}\frac{\partial \mathbf R_{_T}}{\partial u_i}
&=\mathbf R_{_T}^t\Omega_{_T}\frac{\partial(\bar{\mathbf A}\mathbf u_x+\bar{\mathbf B}\mathbf u_y)_{_T}}{\partial u_i}\\
&=\mathbf R_{_T}^t\Omega_{_T}\left(\frac{\partial\bar{\mathbf A}}{\partial u_i}\mathbf u_x+\bar{\mathbf A}\frac{\partial \mathbf u_x}{\partial u_i}+
\frac{\partial\bar{\mathbf B}}{\partial u_i}\mathbf u_y+\bar{\mathbf B}\frac{\partial \mathbf u_y}{\partial u_i}\right)_{_T}
\end{split}
\end{equation}
which indeed reduces to \eqref{mini} when $n=1$.  As for the entries of the Jacobian we get
\begin{equation}
\begin{split}
\frac{\partial ^2 \mathcal{I}_{_T}}{\partial v_j\partial u_i}
&=\left(\left(\frac{\partial\bar{\mathbf A}}{\partial v_j}\mathbf u_x+\bar{\mathbf A}\frac{\partial \mathbf u_x}{\partial v_j}+
\frac{\partial\bar{\mathbf B}}{\partial v_j}\mathbf u_y+\bar{\mathbf B}\frac{\partial \mathbf u_y}{\partial v_j}\right)^t_{_T}\right.\\
& \quad \cdot\left(\frac{\partial\bar{\mathbf A}}{\partial u_i}\mathbf u_x+\bar{\mathbf A}\frac{\partial \mathbf u_x}{\partial u_i}+
\frac{\partial\bar{\mathbf B}}{\partial u_i}\mathbf u_y+\bar{\mathbf B}\frac{\partial \mathbf u_y}{\partial u_i}\right)_{_T}\\
& \quad +\mathbf R_{_T}^t\left(\frac{\partial^2\bar{\mathbf A}}{\partial v_j \partial u_i}
\mathbf u_x+\frac{\partial\bar{\mathbf A}}{\partial v_j}\frac{\partial \mathbf u_x}{\partial u_i}+\frac{\partial\bar{\mathbf A}}{\partial u_i}\frac{\partial \mathbf u_x}{\partial v_j}\right.\\
& \left.\quad +\left.\frac{\partial^2\bar{\mathbf B}}{\partial v_j \partial u_i}
\mathbf u_y+\frac{\partial\bar{\mathbf B}}{\partial v_j}\frac{\partial \mathbf u_y}{\partial u_i}+\frac{\partial\bar{\mathbf B}}{\partial u_i}\frac{\partial \mathbf u_y}{\partial v_j}\right)_{_T}\right)\Omega_{_T}
\end{split}
\end{equation}

\chapter{Inviscid Incompressible Flow}
This chapter concerns the numerical solution of Laplace's equation.  Namely,
\[\Delta \phi = 0 \] for some potential $\phi$.  This is the governing equation
of inviscid incompressible flow Using $\phi_x=u$ and $\phi_y=v$ it can be written as the following system:
\begin{equation}
\begin{split}
0&=u_x + v_y\\
0&=v_x - u_y
\end{split}
\end{equation}
which are also known as the Cauchy-Riemann equations.  The contribution to the functional on
each element is 
\begin{equation}
\mathcal{I}_{_T}=\mathbf R_T^{^t}\mathbf R_{_T}\Omega_{_T}=\frac{1}{2}\left((u_x + v_y)^2+(v_x - u_y)^2\right)\Omega_{_T}
\end{equation}
At a node $i$ of the element we have,
\begin{align}
\frac{\partial\mathcal{I}_{_T}}{\partial u_i}&=-\frac{1}{2}((u_x + v_y)n_{x_i}-(v_x - u_y) n_{y_i})\label{didu}\\
\frac{\partial ^2\mathcal{I}_{_T}}{\partial u_j \partial u_i}
&=\frac{1}{4\Omega_{_T}}(n_{x_i}n_{x_j}+n_{y_i}n_{y_j})\\
\frac{\partial ^2\mathcal{I}_{_T}}{\partial v_j \partial u_i}
&=\frac{1}{4\Omega_{_T}}(n_{x_i}n_{y_j}-n_{y_i}n_{x_j})
\end{align}
\hrulefill\\ 
\begin{align}
\frac{\partial\mathcal{I}_{_T}}{\partial v_i}&=-\frac{1}{2}((u_x + v_y)n_{y_i}+(v_x - u_y) n_{x_i})\\
\frac{\partial ^2\mathcal{I}_{_T}}{\partial u_j \partial v_i}
&=\frac{1}{4\Omega_{_T}}(n_{y_i}n_{x_j}-n_{x_i}n_{y_j})\\
\frac{\partial ^2\mathcal{I}_{_T}}{\partial v_j \partial v_i}
&=\frac{1}{4\Omega_{_T}}(n_{y_i}n_{y_j}+n_{x_i}n_{x_j})\label{didvv}
\end{align}
where $j$ is indexes the nodes of the element.

Equations (\ref{didu})- (\ref{didvv}) provide us with a complete set of equations for all of the inner nodes in the domain.  The
boundary nodes may not use these equations and will be detailed for each problem.  As a numerical verification of the order of accuracy of the method a test is done on a unit square with Dirichlet boundary conditions.

\section{Test Case On A Unit Square}
\begin{figure}[htbp]
\centering
\includegraphics[totalheight=2in]{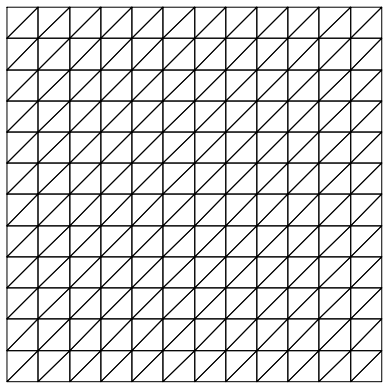}
         \caption{Initial grid \label{reggrid}}
\end{figure}

\begin{figure}[htbp]
\centering
\includegraphics[totalheight=3in]{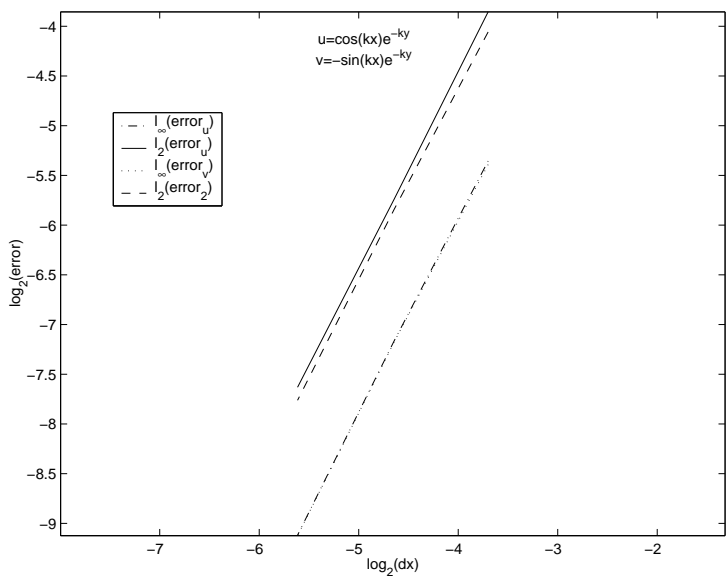}
         \caption{Errors for the sinusoidal wall\label{error1}}
\end{figure}

\begin{figure}[htbp]
\begin{minipage}[t]{.5\linewidth}
\centering
\includegraphics[totalheight=3in]{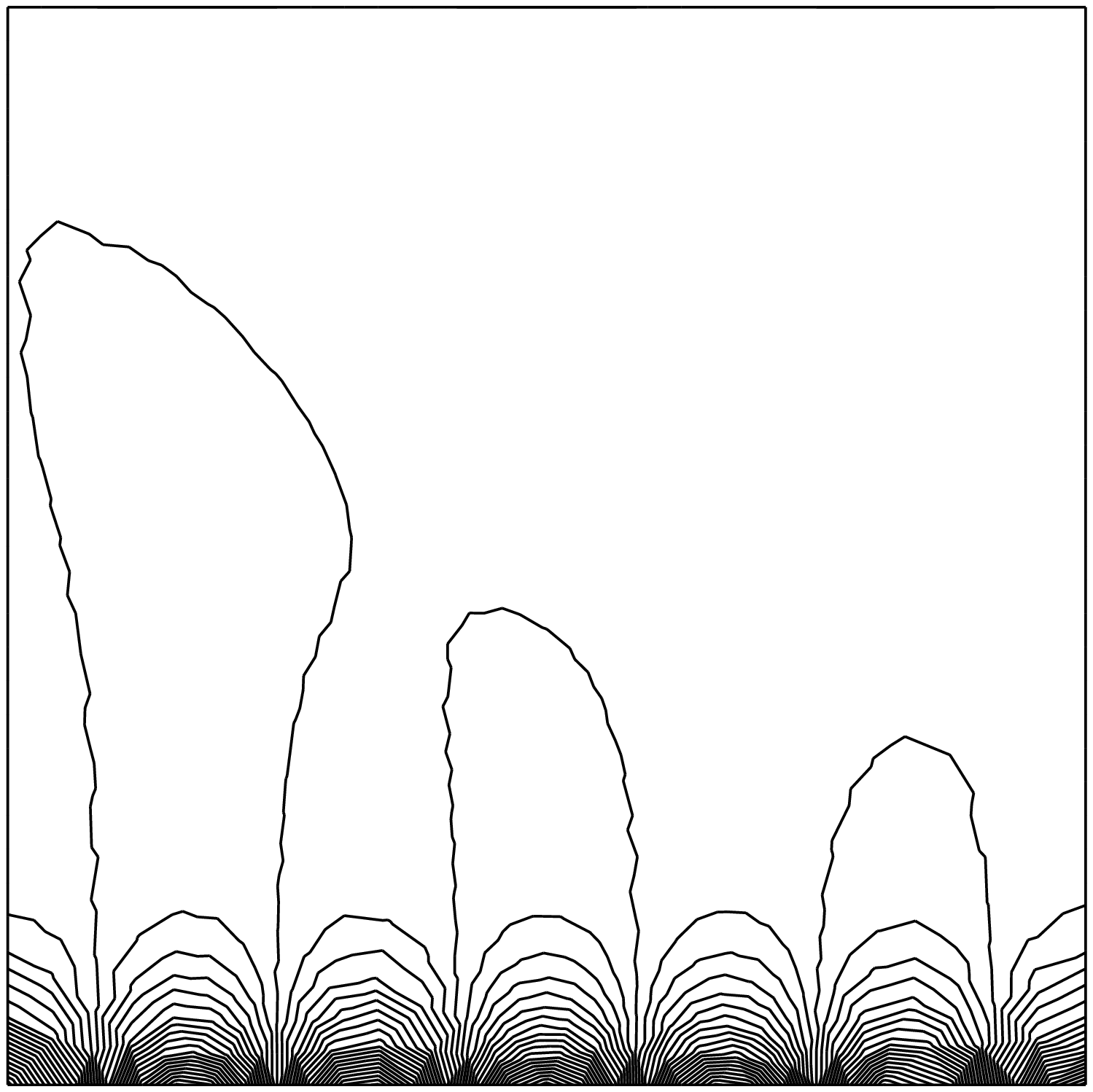}
         \caption{u (numerical)\label{unum1}}
\end{minipage}%
\begin{minipage}[t]{.5\linewidth}
\centering
\includegraphics[totalheight=3in]{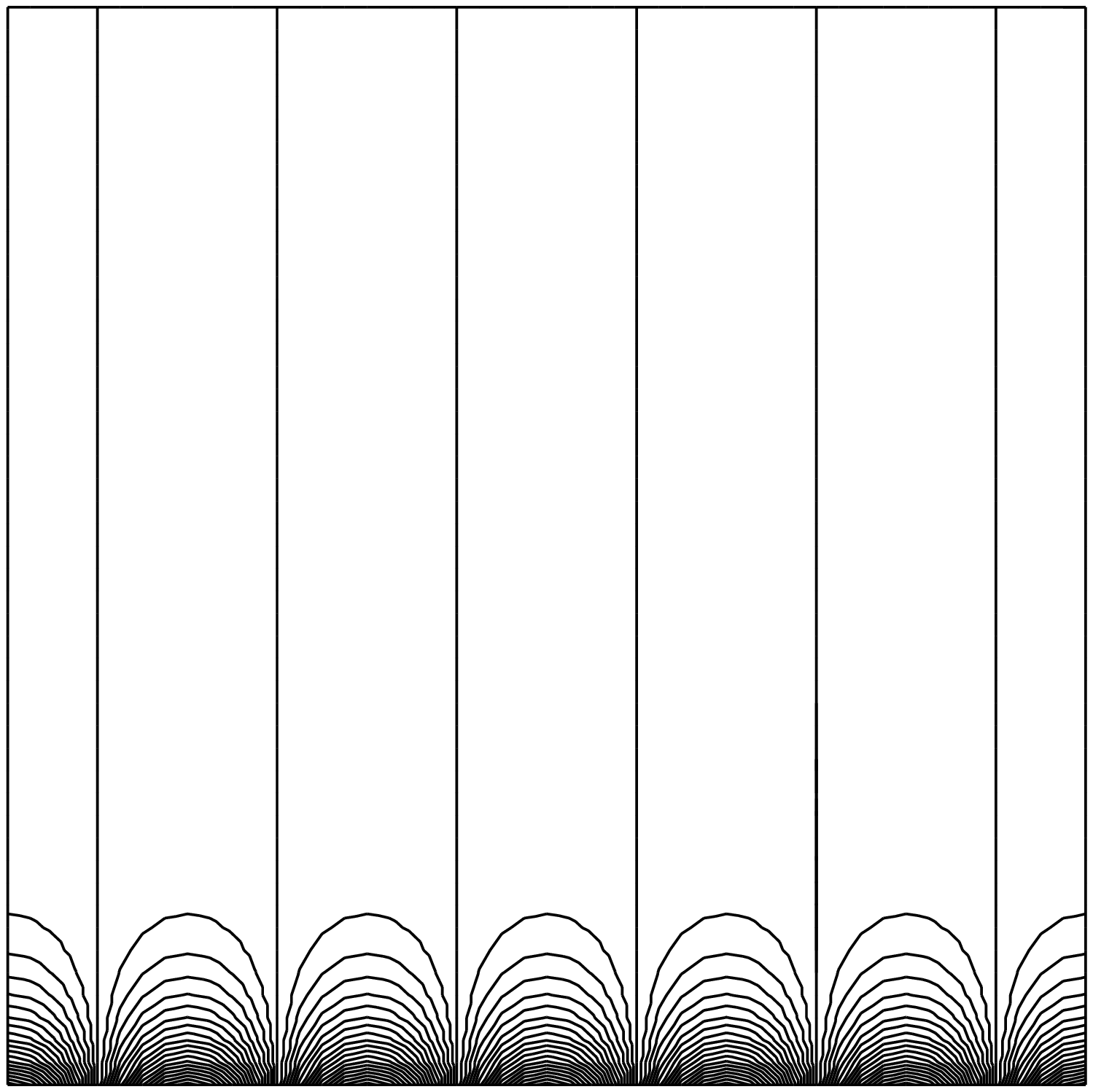}
         \caption{u (analytical)\label{uex1}}
\end{minipage}
\end{figure}

\begin{figure}[htbp]
\begin{minipage}[t]{.5\linewidth}
\centering
\includegraphics[totalheight=3in]{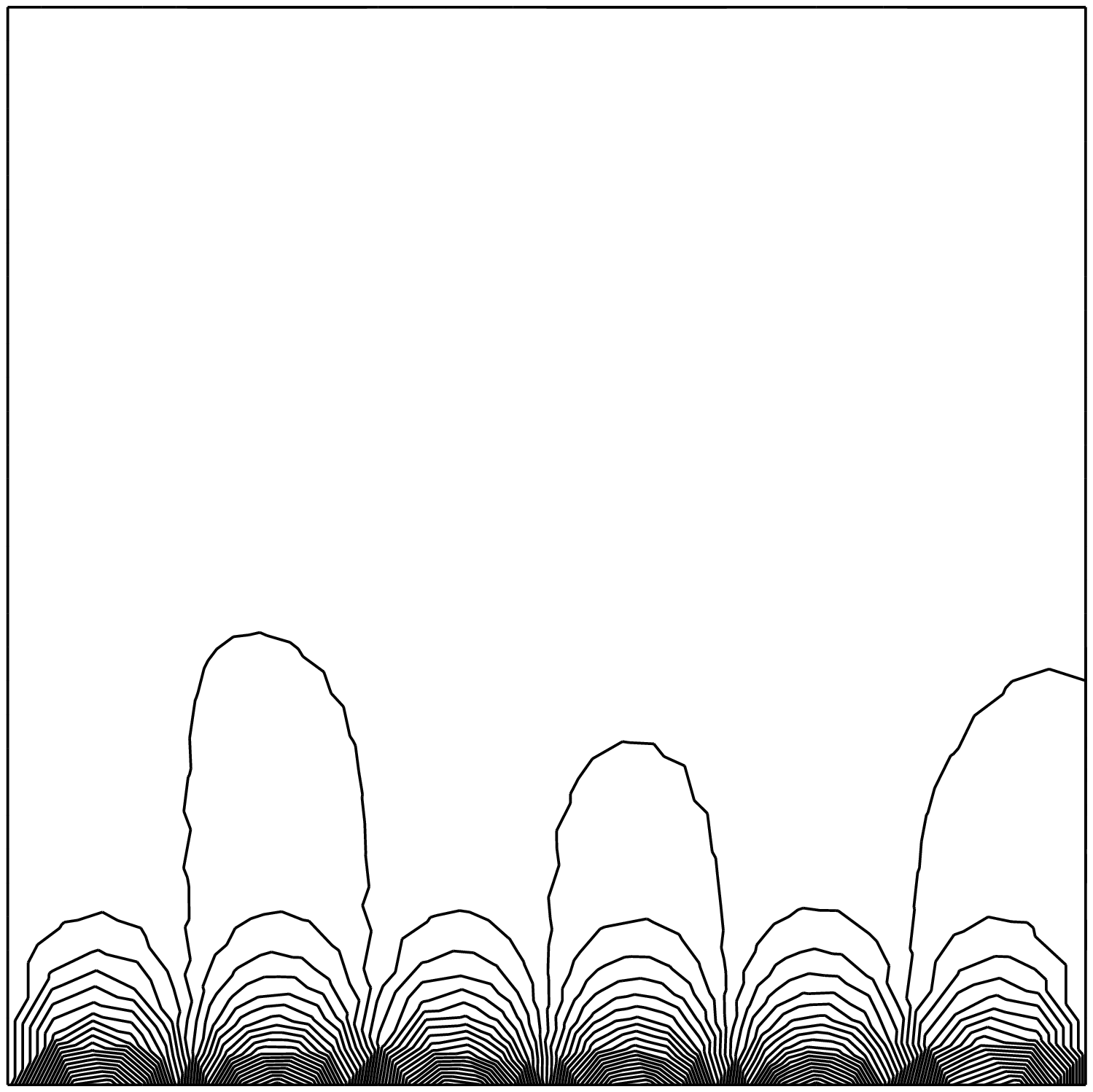}
         \caption{v (numerical)\label{vnum1}}
\end{minipage}%
\begin{minipage}[t]{.5\linewidth}
\centering
\includegraphics[totalheight=3in]{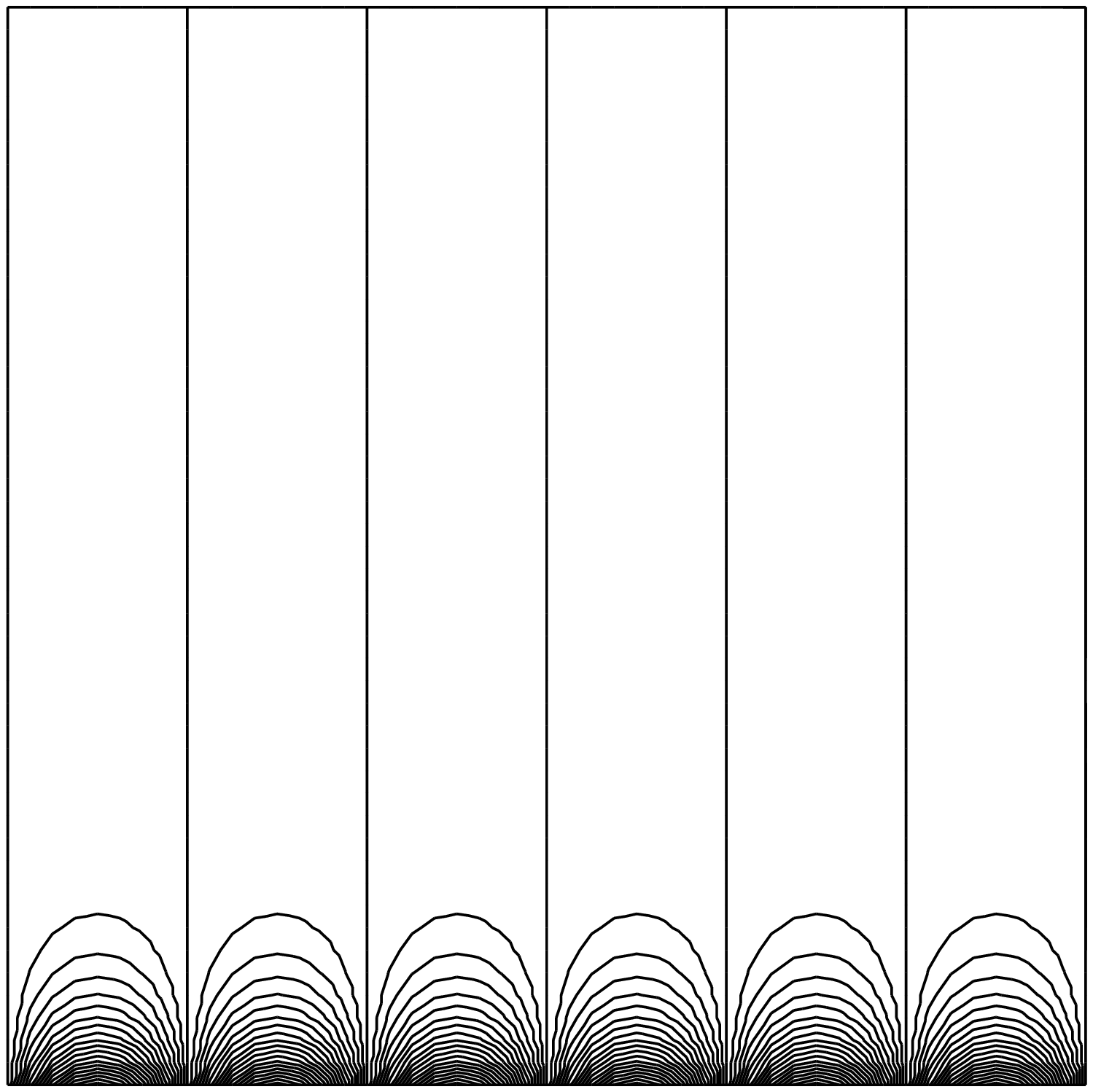}
         \caption{v (analytical)\label{vex1}}
\end{minipage}
\end{figure}

Given a potential $\phi=\frac{1}{k}e^{-ky} \sin(kx)$ we have
\begin{align}
u=\phi_x &=e^{-ky} \cos(kx)\label{equa1}\\
v=\phi_y &=-e^{-ky} \sin(kx)\label{equa2}
\end{align}
where $u$ and $v$ satisfy the Cauchy-Riemann equations.
Actually, the potential corresponds to the flow past an infinite sinusoidal wall.  For this test, $k=6\pi$.
If we are given Dirichlet boundary conditions for the variables they are implemented in a least-squares fashion by summing
the squares of the conditions.  So for this test problem, if $i$ corresponds to a boundary node we have,
\begin{equation}
\mathcal{I}_i=\frac{(u_i-e^{-ky_i} \cos(kx_i))^2+(v_i+e^{-ky_i} \sin(kx_i))^2}{2}
\end{equation}
giving
\begin{align}
\frac{\partial \mathcal{I}_i}{\partial u_i}&=u_i-e^{-ky_i} \cos(kx_i)\label{dir1}\\
\frac{\partial^2 \mathcal{I}_i}{\partial u_i\partial u_i}&=1\\
\frac{\partial \mathcal{I}_i}{\partial v_i}&=v_i+e^{-ky_i} \sin(kx_i)\\
\frac{\partial^2 \mathcal{I}_i}{\partial v_i\partial v_i}&=1\label{dir2}
\end{align}
\paragraph{}
Now, the system of equations can be formulated using the following process:
\begin{description}
\item loop $T=1,T_{max}$
  \begin{description}
  \item loop $i\in \Omega_{_T}$
    \begin{description}
      \item if($i\ni d\Omega$) 
         \begin{description}
          \item loop $j\in \Omega_{_T}$
          \item   update functional variation $\mathcal{F}$ and Jacobian $J$ using (\ref{didu})- (\ref{didvv})
           \item end loop
           \item else update functional variation and Jacobian using (\ref{dir1})- (\ref{dir2})
           \end{description}
    \end{description}
  \item end loop
  \end{description}
\item end loop
\end{description}
Here, $T_{max}$ is the total number of triangles in the domain.
\paragraph{}
Since we have Dirichlet boundary conditions, the resulting Jacobian is symmetric and positive definite.  Therefore,
$J\vec{q}=-\vec{\mathcal{F}}$ is solved using the preconditioned conjugate gradient method.  A Jacobi preconditioner is used.
This preconditioner is defined as the matrix $P$ where
\begin{eqnarray}
P_{l,m}&=&\left\{ \begin{array}{r@{\quad,\quad}l} \frac{1}{J_{l,m}}& l=m\\ 0 & otherwise\end{array}\right.\label{jacobi}
\end{eqnarray} 
The initial grid is a structured $13\times 13$ right
triangulation of the unit square (see figure \ref{reggrid}).  In this test, upon each grid refinement the solution 
converged in no more than 3 corrections with $\|\vec{q}\|<tol=10^{-8}$.    Due to the boundary conditions and grid configuration, the method yielded a 
Jacobian in which the velocity components decoupled entirely (see Appendix B).   Whether or not the grid was structured, the preconditioned conjugate gradient method consistently converged in at most, half as many iterations as there are unknowns
in the system.  In general this property cannot be guaranteed.  The errors are computed using
\begin{equation}
{\|f\|}_\infty=\max|f(x,y)|
\end{equation}
and
\begin{equation}
{\|f\|}_p=\left[\sum_i\sum_{T_j}\Omega_{_{T_j}}|f_i|^p\right]^{\frac{1}{p}}
\end{equation}
where $\Omega_{_{T_j}}$ is an element abutting the node $i$.  For this problem and those following, Newton's method is always
initialized with a guess that depicts uniform flow.
\paragraph{}
From figure \ref{error1} it is clear that second order accuracy is observed for this method.  This test was repeated for larger
values of $k$ and unsurprisingly, the method was still observed to be second order accurate though the errors did increase for corresponding increases of $k$.  For
completeness a case is shown with $k=6\pi$ on a fully unstructured grid.  Table \ref{table1} shows the errors in the velocity components.  Figures (\ref{unum1})-(\ref{vex1}) look reasonable and can
be improved with a finer grid.  These results show we can expect the least-squares method to yield nice results for smooth flows if
Dirichlet boundary conditions are enforced.  

\begin{figure}[htbp]
\centering
\includegraphics[totalheight=3.5in]{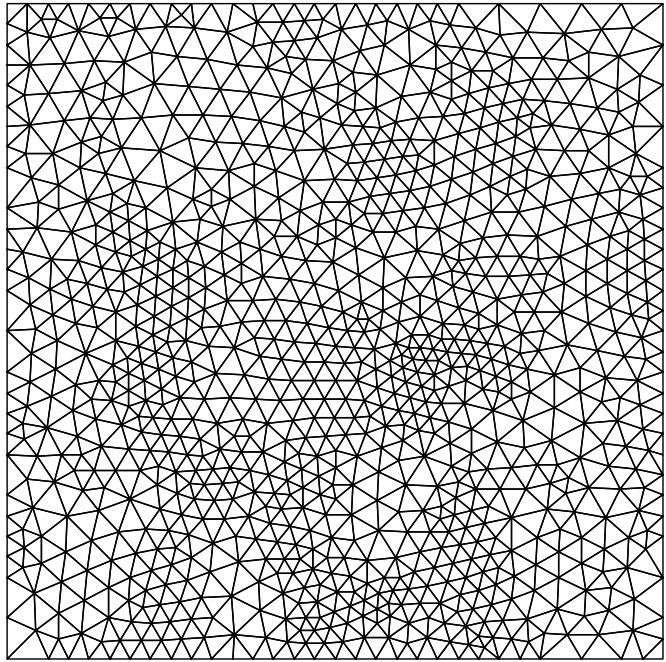}
         \caption{Unstructured grid (976 points) \label{unstruc1}}
\end{figure} 
\begin{table}[htbp]\centering
\begin{tabular}{|c|c|} \hline
Error & Value\\ \hline \hline
${\|e_u\|}_2$ & 0.005356007786650060 \\ 
${\|e_u\|}_\infty$ & 0.036229598337966351\\ 
${\|e_v\|}_2$ & 0.004727775128293196\\ 
${\|e_v\|}_\infty$ & 0.025036120477625323 \\ \hline
\end{tabular}
\caption{Error table (unstructured grid)}\label{table1}
\end{table}

\section{Inviscid Flow Over A Parabolic Profile}
In the last section the flow was reasonably smooth.  Of course we are also interested in how the least-squares method handles
flow with singular behavior.  For a parabolic profile there is a singularity at the leading and trailing edge of the profile.
The analytical solutions are (see \cite{chattot3}):
\begin{align}
u_{exact}&= \frac{k}{a^2}\left[ x\ln\left(\frac{r_2}{r_1}\right)-y(\theta_2-\theta_1)\right]-\frac{2k}{a}\\
v_{exact}&=-\frac{k}{a^2}\left[ y\ln\left(\frac{r_2}{r_1}\right)+x(\theta_2-\theta_1)\right]
\end{align}
where
\begin{align}
r_{2,1}=\sqrt{(x\pm a)^2 +y^2}\\
\intertext{and}
\theta_{2,1}=\arctan\left(\frac{y}{x\pm a}\right)
\end{align}
For solutions computed here $a=\frac{1}{2}$ and $k=-\frac{1}{4}$.

\subsection{Boundary Conditions}
The computational domain is a rectangle with dimension $[-2,2]\times[0,2]$. For $v$,  Dirichlet conditions are specified everywhere 
and for $u$, everywhere except the bottom side.  There, the functional for the governing equations is minimized with respect to $u$.  In a least-squares formulation we have:\\\\
\underline{Left, Right, Top}\\
\begin{align}
\mathcal{I}_i&=\frac{(u-u_{exact})^2_i+(v-v_{exact})^2_i}{2}\\
\frac{\partial \mathcal{I}_i}{\partial u_i}&=(u-u_{exact})_i\\
\frac{\partial^2 \mathcal{I}_i}{\partial u_i\partial u_i}&=1\\
\frac{\partial \mathcal{I}_i}{\partial v_i}&=(v-v_{exact})_i\\
\frac{\partial^2 \mathcal{I}_i}{\partial v_i\partial v_i}&=1
\end{align}
\hrulefill\\ 
\underline{Bottom}\\
$v$:
\begin{align}
\mathcal{I}_i&=\frac{(v-v_{exact})^2_i}{2}\\
\frac{\partial \mathcal{I}_i}{\partial v_i}&=(v-v_{exact})_i\\
\frac{\partial^2 \mathcal{I}_i}{\partial v_i\partial v_i}&=1
\end{align}
$u$:
\begin{align}
\mathcal{I}_{_T}&=\frac{1}{2}\left((u_x + v_y)^2+(v_x - u_y)^2\right)\Omega_{_T}\\
\frac{\partial\mathcal{I}_{_T}}{\partial u_i}&=-\frac{1}{2}((u_x + v_y)n_{x_i}-(v_x - u_y) n_{y_i})\\
\frac{\partial ^2\mathcal{I}_{_T}}{\partial u_j \partial u_i}
&=\frac{1}{4\Omega_{_T}}(n_{x_i}n_{x_j}+n_{y_i}n_{y_j})\\
\frac{\partial ^2\mathcal{I}_{_T}}{\partial v_j \partial u_i}
&=\frac{1}{4\Omega_{_T}}(n_{x_i}n_{y_j}-n_{y_i}n_{x_j})
\end{align}
\subsection{Results}
\begin{figure}[htbp]
\centering
\includegraphics[width=6in]{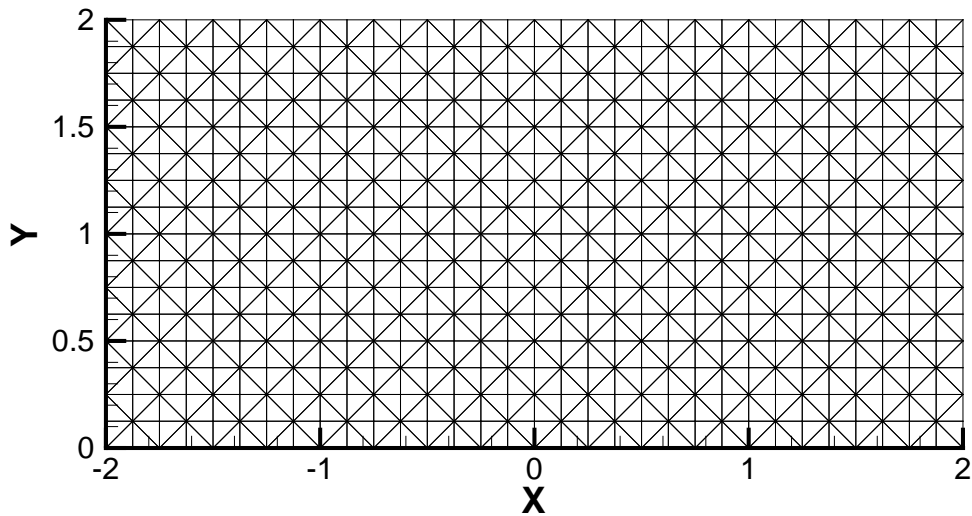}
\caption{Initial Grid Used For The Computations \label{grid2}}
\end{figure}

\begin{figure}[htbp]
        \centering
         \includegraphics[totalheight=3in]{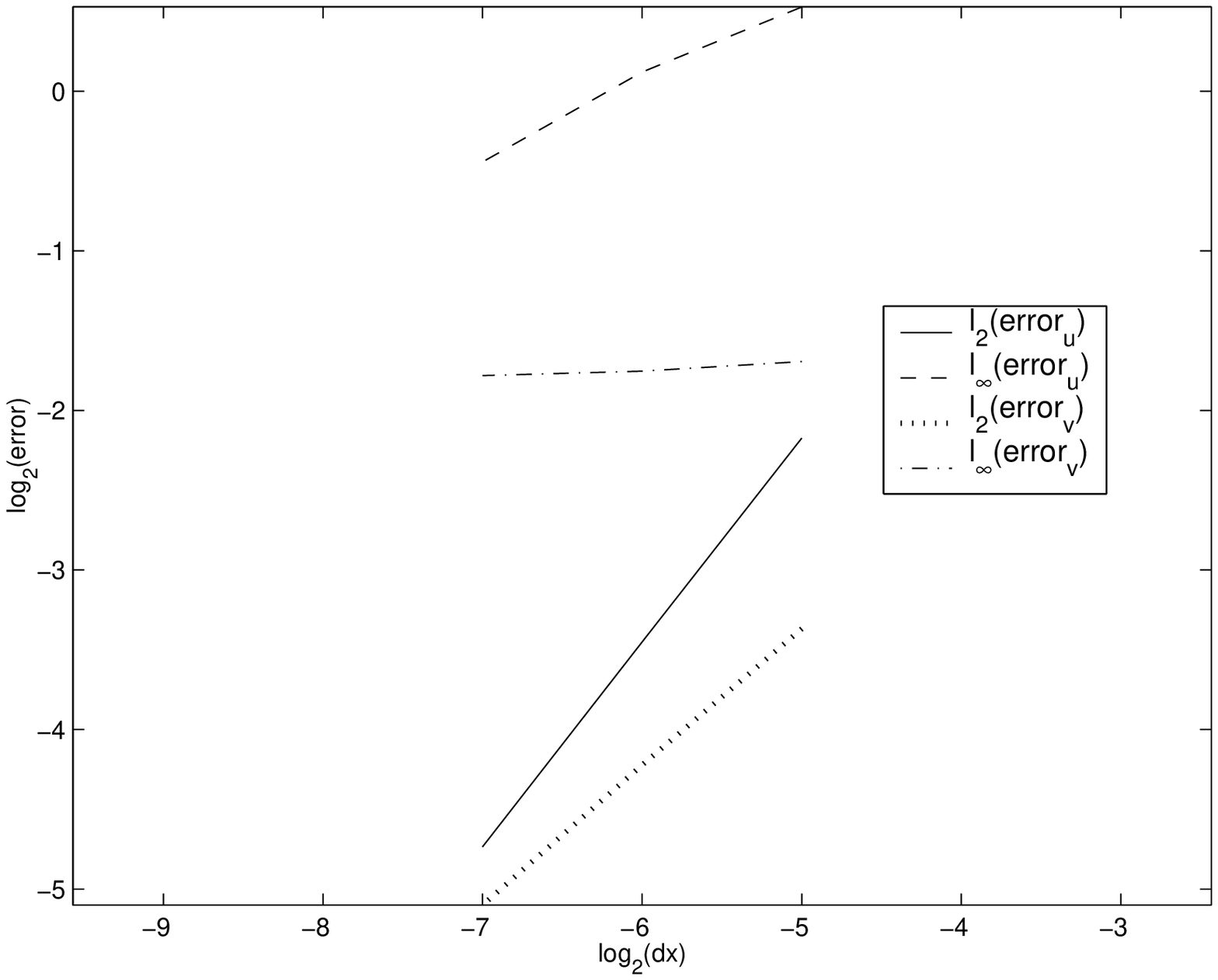}
\caption{Velocity Distribution On The Lower Boundary \label{error2}}
\end{figure}

\begin{figure}[htbp]
\centering
\includegraphics[totalheight=3in]{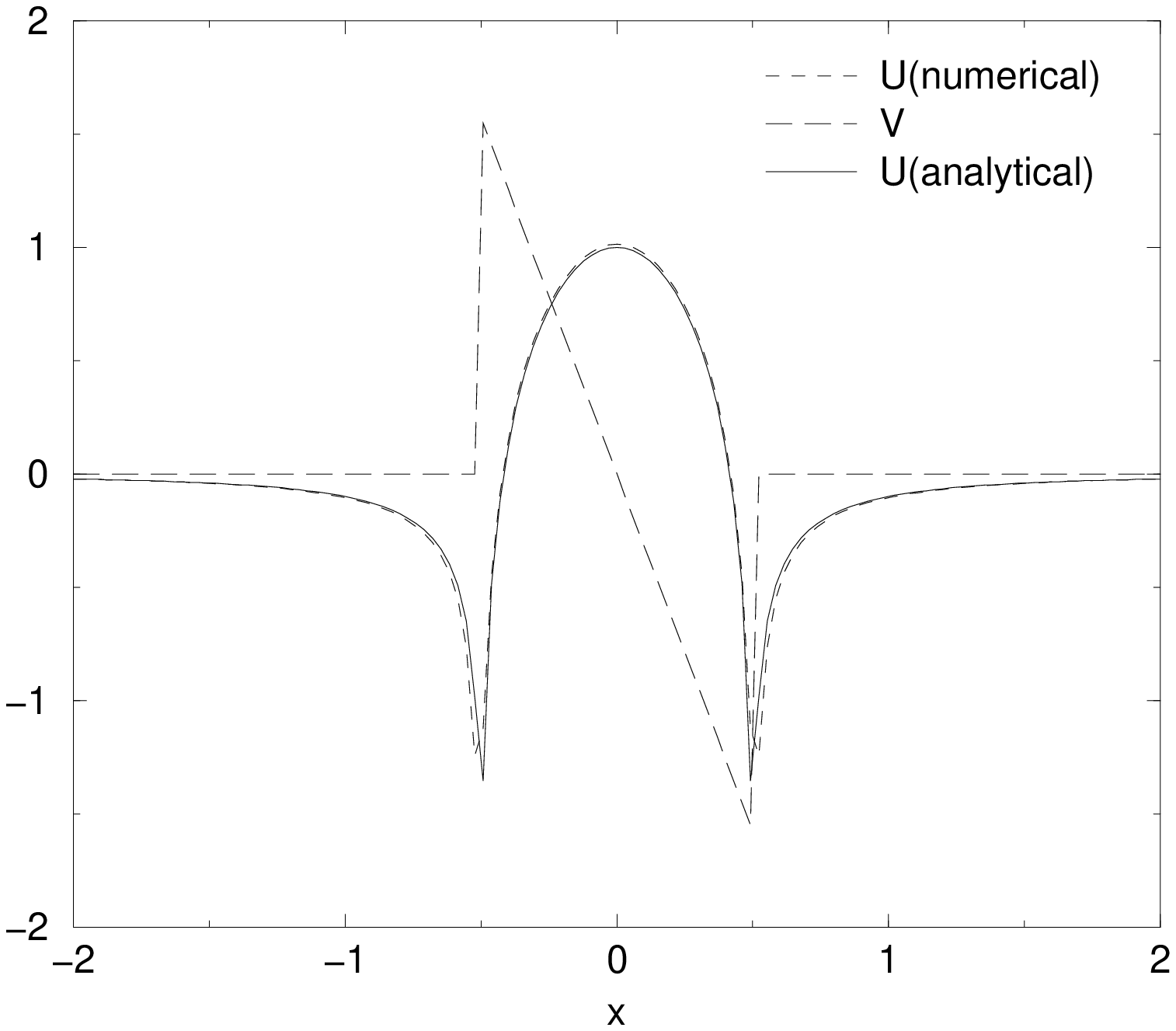}
\caption{Velocity Distribution On The Lower Boundary \label{profile}}
\end{figure}

Again, the preconditioned conjugate gradient method is used to solve the system of equations.  The initial grid is shown in figure 
\ref{grid2}.  From figure \ref{error2} observe that the error does not tend to zero.  The method appears to be zeroth order in $l_{\infty}$
for $v$ and $u$.  This is due to the logarithmic behavior of $u$ and the jump in $v$ on the boundary $y=0$. 
The $l_2$ norm of $u$ is almost 2.  In the $l_2$ norm the order of the method for $u$ seems to be 1 greater than that of $v$.  This is surprising since $v$ is known everywhere on the boundary.  Note the excellent agreement in $u$
and $u_{exact}$ on the profile as illustrated in figure \ref{profile} for the finest grid.  The remaining figures show contour plots
corresponding to the finest grid.

\begin{figure}[htbp]
        \centering
         \includegraphics[totalheight=3.7in]{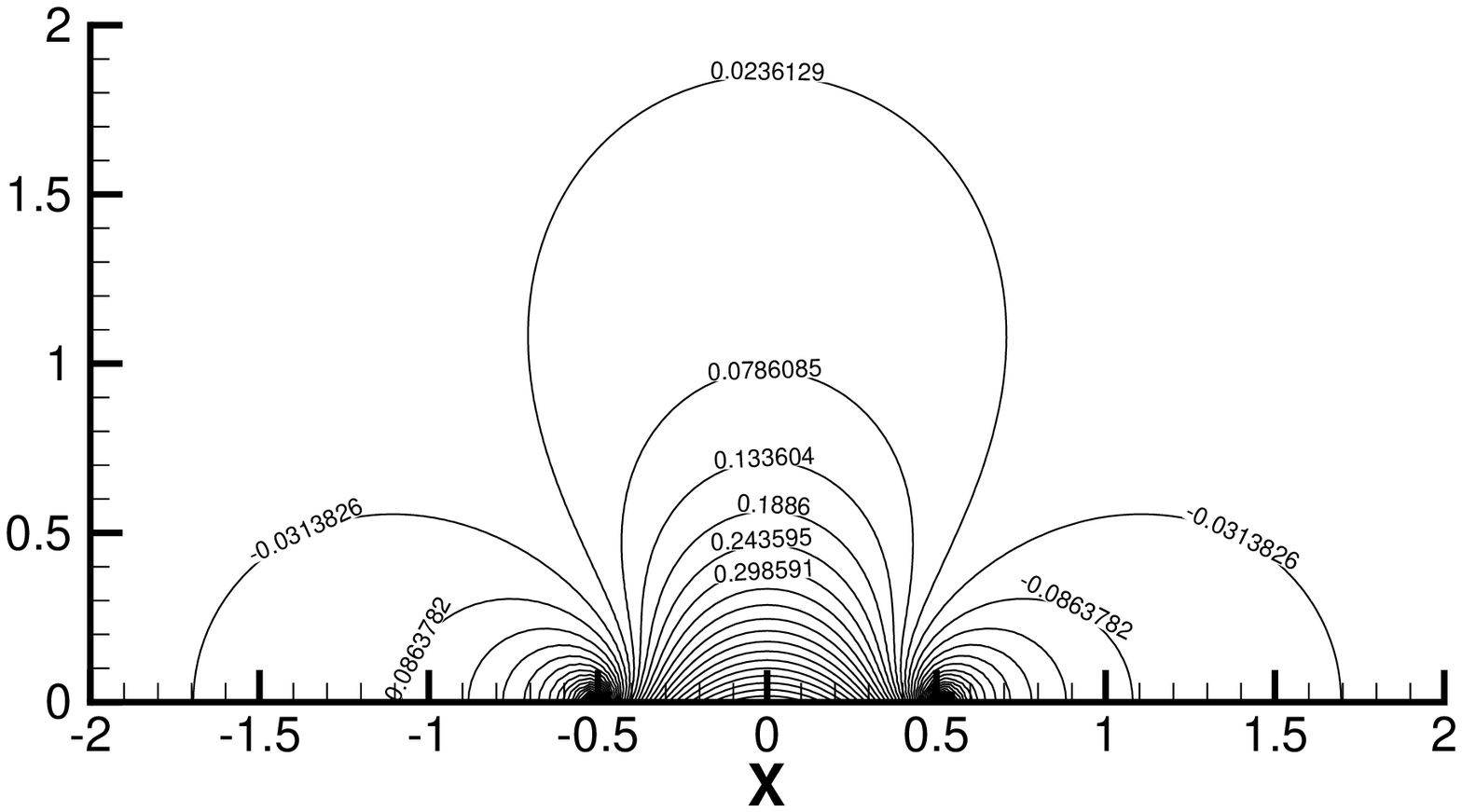}\caption{u (numerical)} 
\end{figure}
\begin{figure}[htbp]
        \centering
         \includegraphics[totalheight=3.7in]{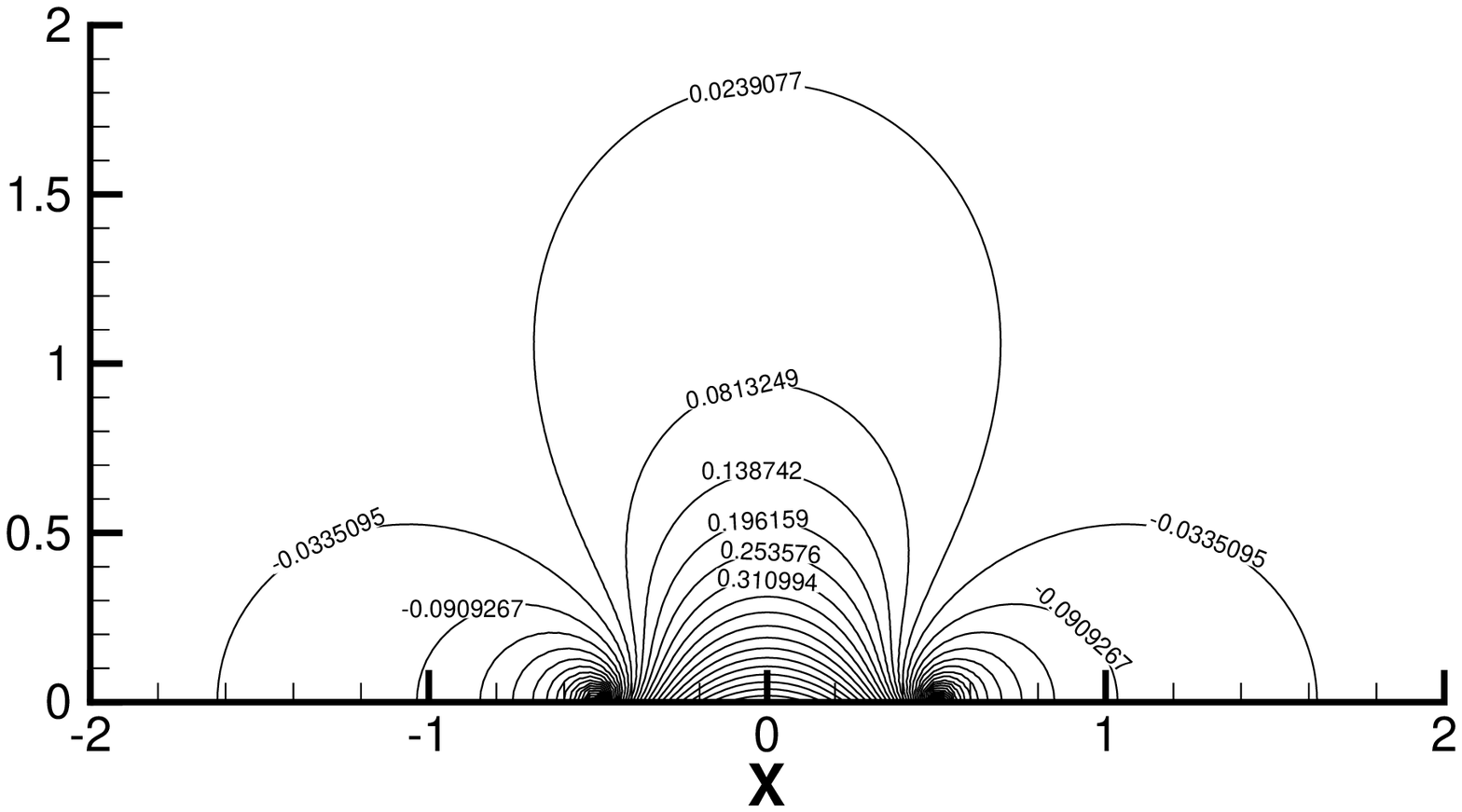}\caption{u (analytical)} 
\end{figure}

\begin{figure}[htbp]
        \centering
         \includegraphics[totalheight=3.7in]{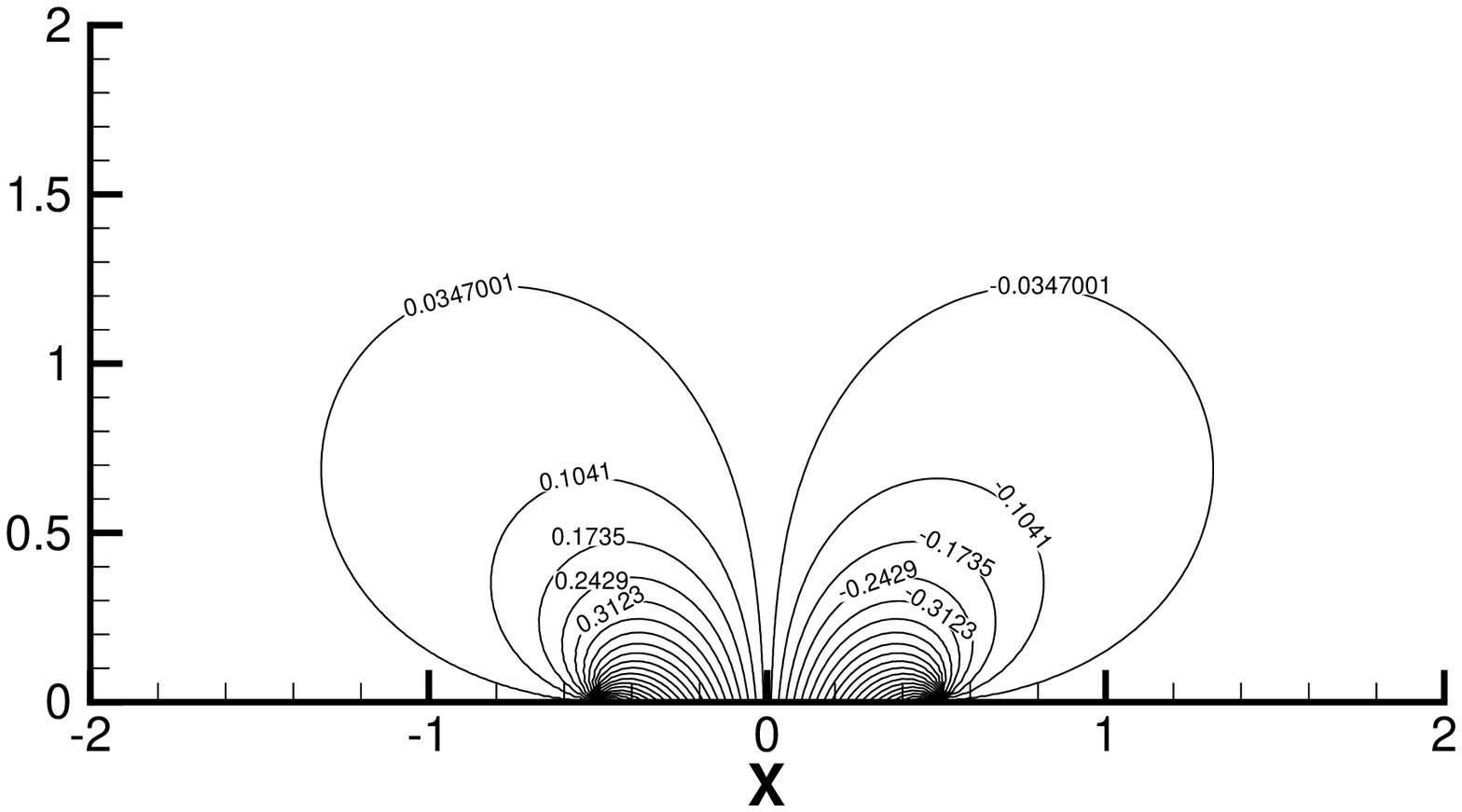}\caption{v (numerical)} 
\end{figure}
\begin{figure}[htbp]
        \centering
         \includegraphics[totalheight=3.7in]{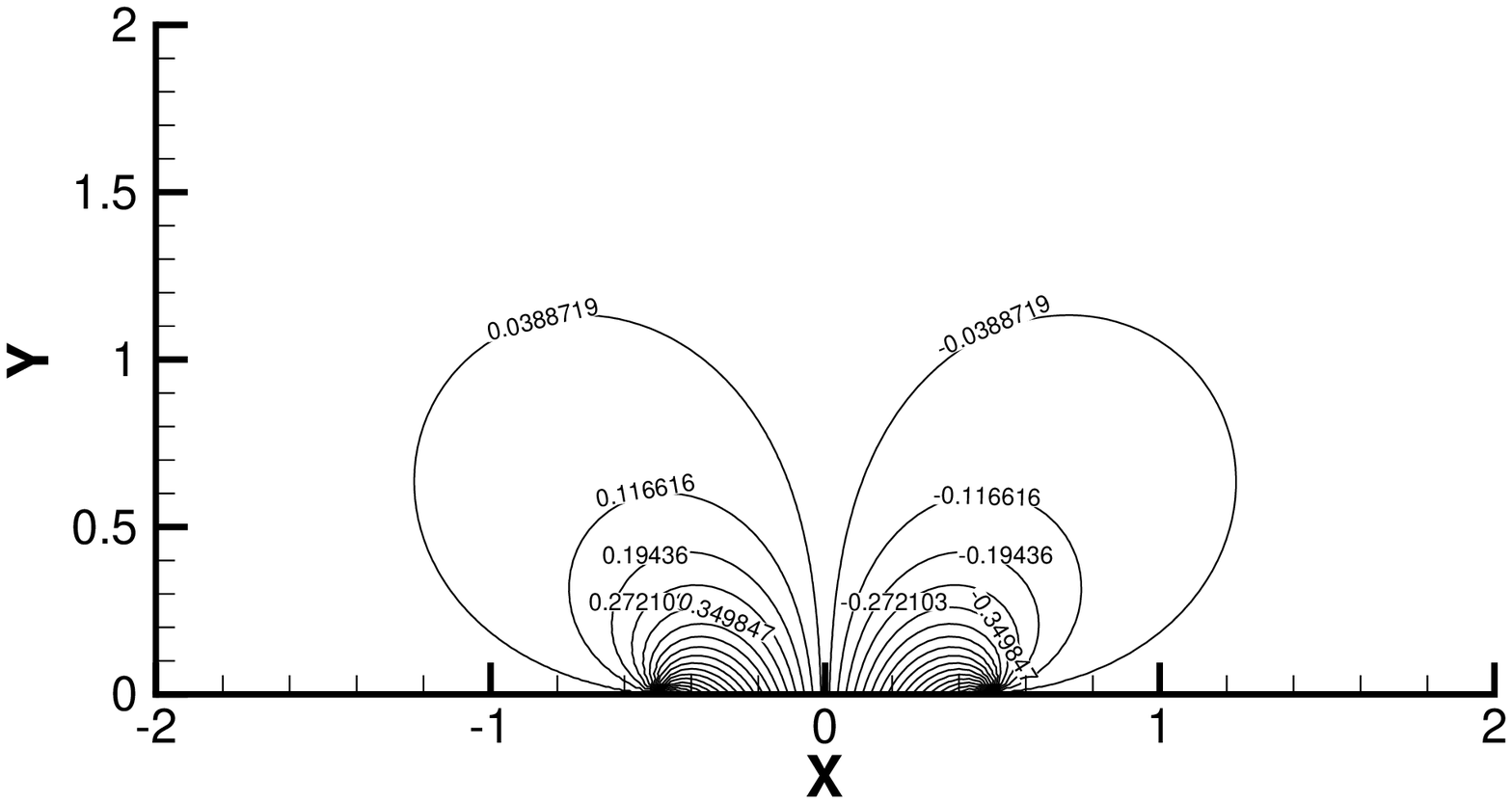}\caption{v (analytical)} 
\end{figure}
\section{Inviscid Flow Over A Cylinder}
In this section the quality of the least-squares method is studied for flow over a cylinder
with circulation.  The flow is again governed by the Cauchy-Riemann equations.
The analytical solutions for the velocity components are:

\begin{align}
u_{exact}&=1-\frac{y\Gamma}{2\pi r^2}-a^2(\frac{1}{r^2}-2\frac{y^2}{r^4})\\
v_{exact}&=\frac{x\Gamma}{2\pi r^2}-\frac{2a^2xy}{r^4}
\end{align}
where $\Gamma$ denotes the circulation.  As usual, for all inner nodes in the domain equations (\ref{didu})- (\ref{didvv}) apply.

\subsection{Boundary Conditions}
\underline{Farfield}\\
In the farfield the exact solutions are enforced:
\begin{align}
\mathcal{I}_i&=\frac{(u-u_{exact})^2+(v-v_{exact})^2_i}{2}\label{cyl1}\\
\frac{\partial \mathcal{I}_i}{\partial u_i}&=(u-u_{exact})_i\\
\frac{\partial^2 \mathcal{I}_i}{\partial u_i\partial u_i}&=1\\
\frac{\partial \mathcal{I}_i}{\partial v_i}&=(v-v_{exact})_i\\
\frac{\partial^2 \mathcal{I}_i}{\partial v_i\partial v_i}&=1\label{cyl2}
\end{align}
\hrulefill\\ 
\underline{Cylinder's Surface}\\
For nodes on the surface, in addition to enforcing equations (\ref{didu})- (\ref{didvv}) used in general for the inner nodes, contributions due to the boundary conditions are added.  This is consistent with the least-squares approach and provides the correct minimization.  For inviscid flow the necessary boundary condition on the surface is $\mathbf u \cdot \mathbf n =0$.  This is
enforced at a node $i$ on the surface in a least-squares formulation:
\begin{align}
\mathcal{I}_i&=\frac{(u n_{x_i}+v n_{y_i})^2}{2}\label{norm1}\\
\frac{\partial \mathcal{I}_i}{\partial u_i}&=(u n_{x_i}+v n_{y_i})n_{x_i}\\
\frac{\partial^2 \mathcal{I}_i}{\partial u_i\partial u_i}&=(n_{x_i})^2\\
\frac{\partial^2 \mathcal{I}_i}{\partial v_i\partial u_i}&=n_{x_i} n_{y_i}\\
\frac{\partial \mathcal{I}_i}{\partial v_i}&=(u n_{x_i}+v n_{y_i})n_{y_i}\\
\frac{\partial^2 \mathcal{I}_i}{\partial v_i\partial v_i}&=(n_{y_i})^2\\
\frac{\partial^2 \mathcal{I}_i}{\partial u_i\partial v_i}&=\frac{\partial^2 \mathcal{I}_i}{\partial v_i\partial u_i}\label{norm11}
\end{align}
Here, the normal vector is actually of unit length in contrast to the convention established earlier.  For a disc centered at the origin it is simply
\begin{equation}
\mathbf n_i = \frac{\langle x_i,y_i \rangle}{(x_i)^2+(y_i)^2}
\end{equation}
Enforcement of the tangency condition compromises the symmetry of the Jacobian.  Therefore, GMRES along with the aforementioned
Jacobi preconditioner (equation (\ref{jacobi})) is used within Newton's method as opposed to the preconditioned conjugate gradient method.
\subsection{Results}
\begin{figure}[H]
\centering
\includegraphics[totalheight=3.6in]{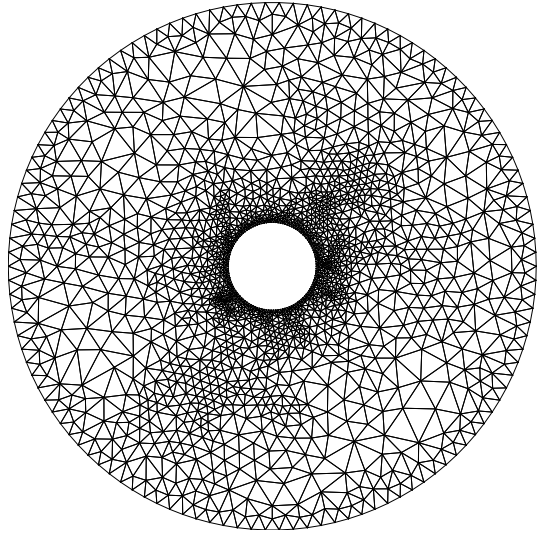}
\caption{Grid 1 (1902 nodes 3568 elements)}\label{cylgrid1}
\end{figure}

\begin{figure}[H]
        \centering
\includegraphics[totalheight=3.7in]{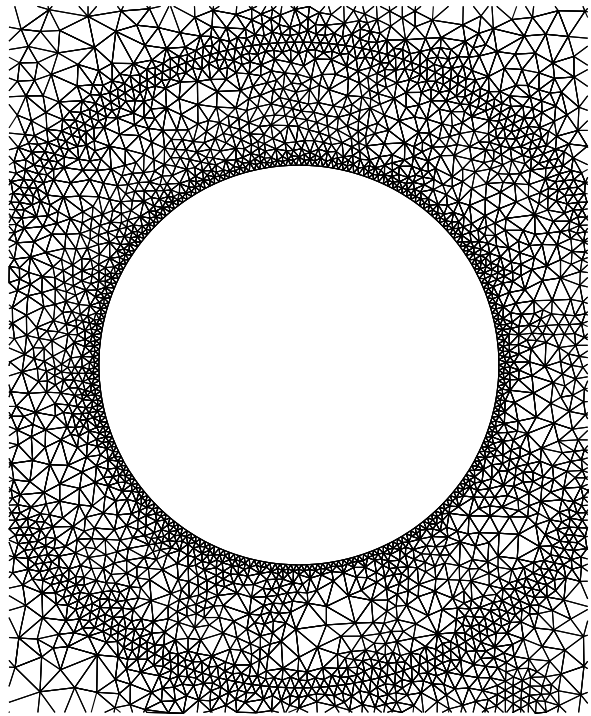}
\caption{Close Up of Grid 2 ($a=.5$, 238 equi-spaced surface points)}\label{cylgrid2}
\end{figure}

\begin{figure}[htbp]
\centering
\includegraphics[totalheight=3.7in]{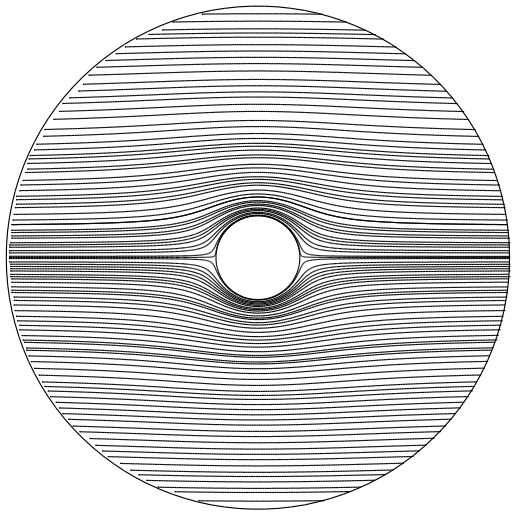}
\caption{Streamlines ($\Gamma=0$) \label{gamma0}}
\end{figure}

\begin{figure}[htbp]
        \centering
                \includegraphics[totalheight=3.7in]{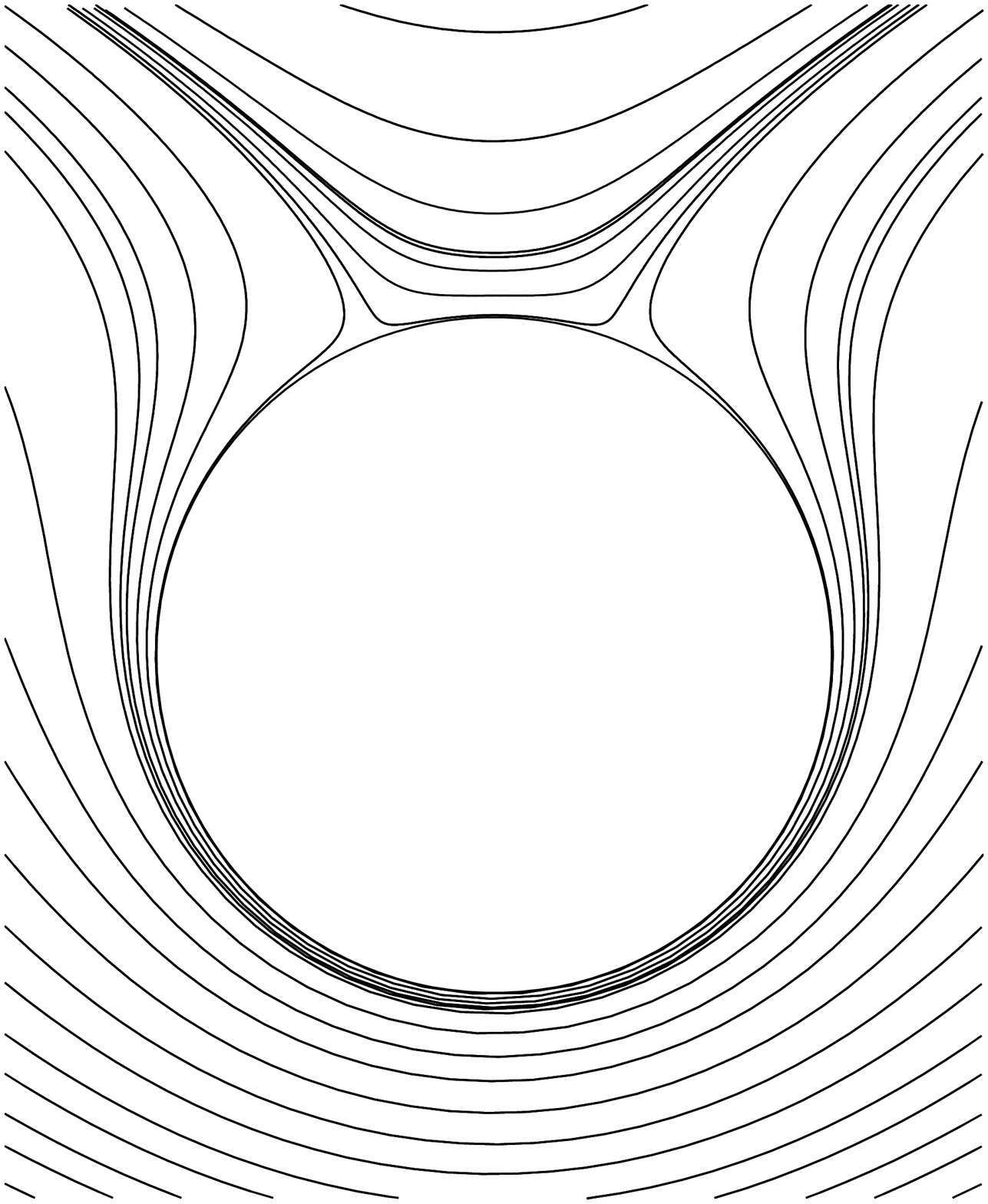}
\caption{Streamlines ($\Gamma=4\pi a$)}\label{1801}
\end{figure}

\begin{figure}[htbp]
        \centering
                \includegraphics[totalheight=3.7in]{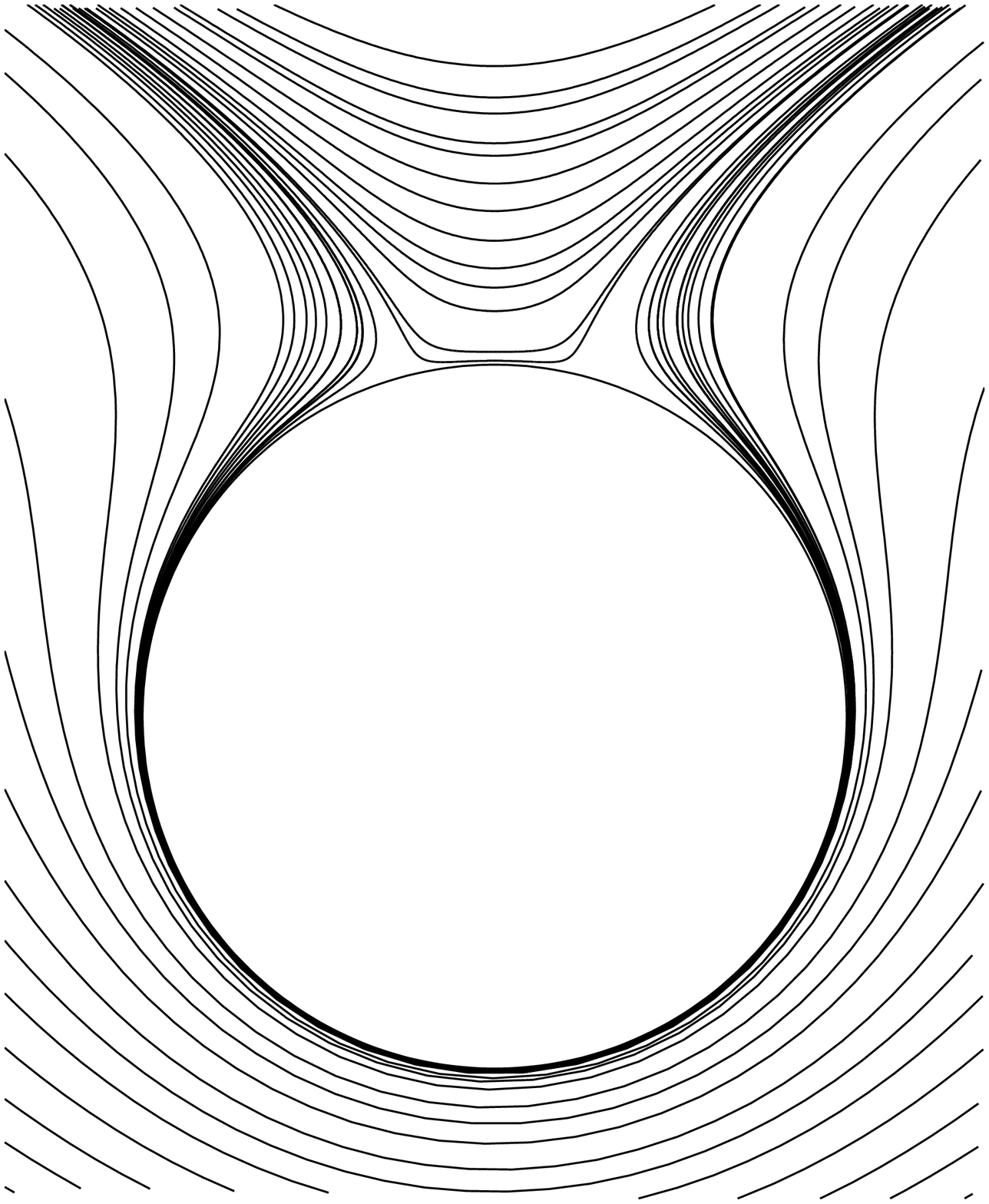}
\caption{Streamlines ($\Gamma=4\pi a$)}\label{1802}
\end{figure}

\begin{figure}[htbp]
        \centering
                \includegraphics[totalheight=3.7in]{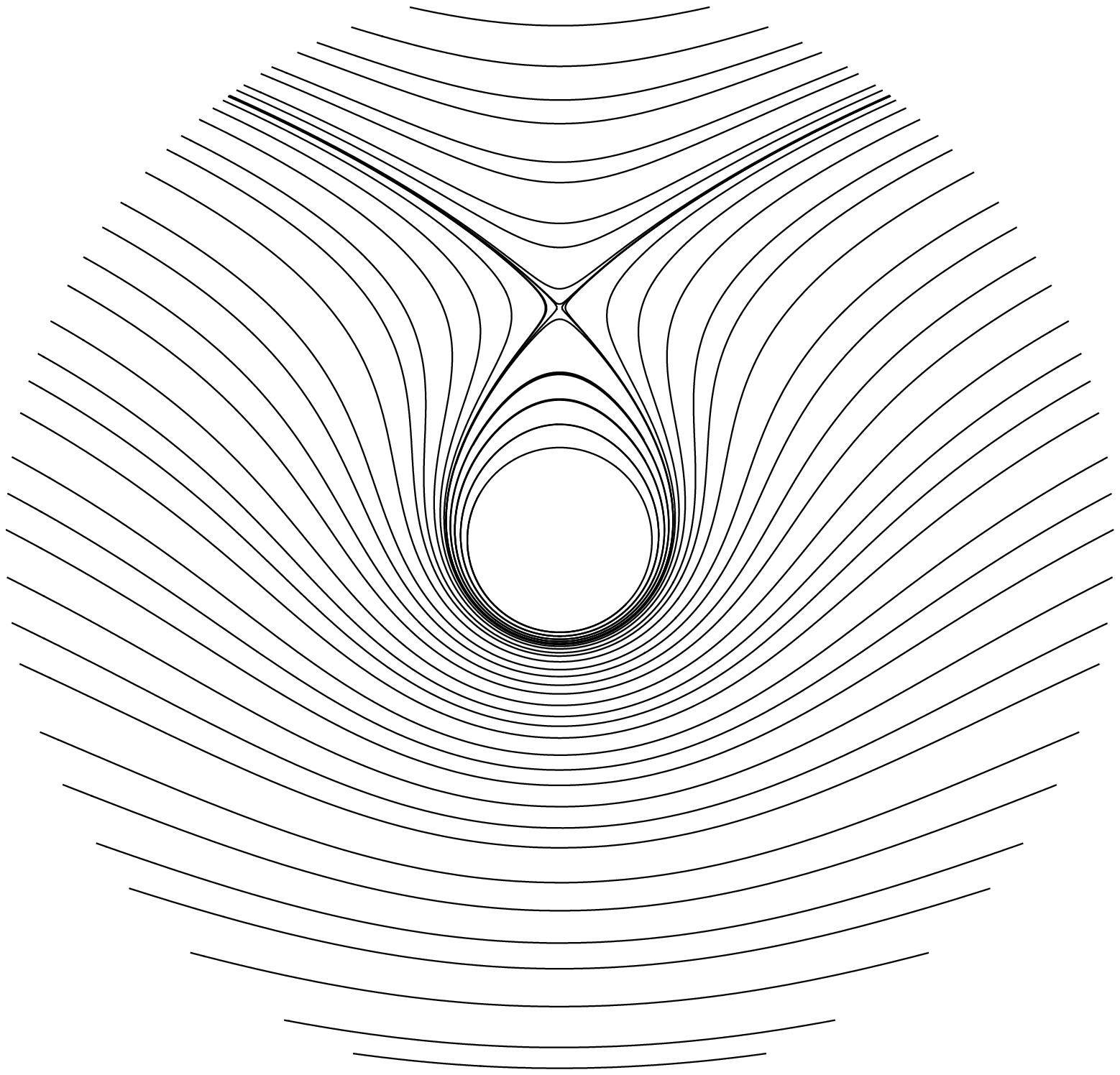}
\caption{Streamlines ($\Gamma=6\pi a$) \label{stag1}}
\end{figure}

\begin{figure}[htbp]
        \centering
                \includegraphics[totalheight=3.7in]{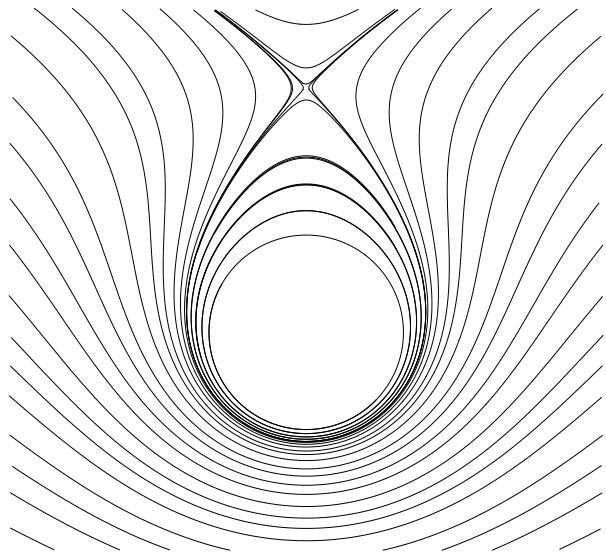}
\caption{Streamlines ($\Gamma=6 \pi a$), Stagnation point at $(x,y)=(0.0004,1.262)$ \label{stag2}}
\end{figure}

Several cases were done for different values of the circulation.  The quality with respect to the errors and the ability of the
method to catch the correct stagnation angles on the cylinder were analyzed.  Analytically, the stagnation angles are at $\theta$
and $\pi-\theta$ where $\theta$ is given by:
\begin{equation}
\theta=\sin ^{-1}\left(\frac{\Gamma}{4 \pi au_{\infty}}\right)
\end{equation}
For values $\Gamma > 4\pi au_{\infty}$ there is one stagnation point above the cylinder at the point (see \cite{white})
\begin{equation}
y=\frac{a}{2}(\beta +(\beta^2-4)^{\frac{1}{2}})
\end{equation}
where 
\begin{equation}\beta=\frac{\Gamma}{2 \pi au_{\infty}}\end{equation}
For these calculations $a=.5$ and $u_{\infty}=1$.
Two grids were used in the calculations and are illustrated in figures \ref{cylgrid1}-\ref{cylgrid2}.  The farfield has a radius
of 3 in both grids.
\paragraph{}
From tables \ref{table3} and \ref{table4} it is clear that the errors drop approximately by a factor slightly smaller than $\frac{1}{2}$ in each velocity component on the finer mesh.  On both meshes, the errors increase as $\Gamma$ was increased.  Also, note
that the errors in $u$ are always larger than the errors in $v$ for all values of $\Gamma$ tested and on either grid.  In Table
\ref{table5} the data suggests the error of the location of the stagnation points on the surface of the cylinder increases as $\Gamma$ increases on both grids and is
illustrated in the streamline plots of figures \ref{gamma0}-\ref{1802}.  A stagnation point was not found on either grid for $\Gamma=4 \pi a$ which is visible in figures \ref{1801} and \ref{1802}.  For $\Gamma= 6 \pi a$ the stagnation point found on the finest grid 
 (see figures \ref{stag1} and \ref{stag2}) is located at $(x,y)=(0.0004,1.262)$ which compares favorably to the analytical values of $(x,y)=(0,1.3)$.  This suggests that the difficulty of
resolving the correct stagnation points is not simply dependent on the magnitude of $\Gamma$ but also, whether or not the stagnation 
points lie on the surface of the cylinder.  
\begin{figure}[H]
\centering
\includegraphics[totalheight=3.7in]{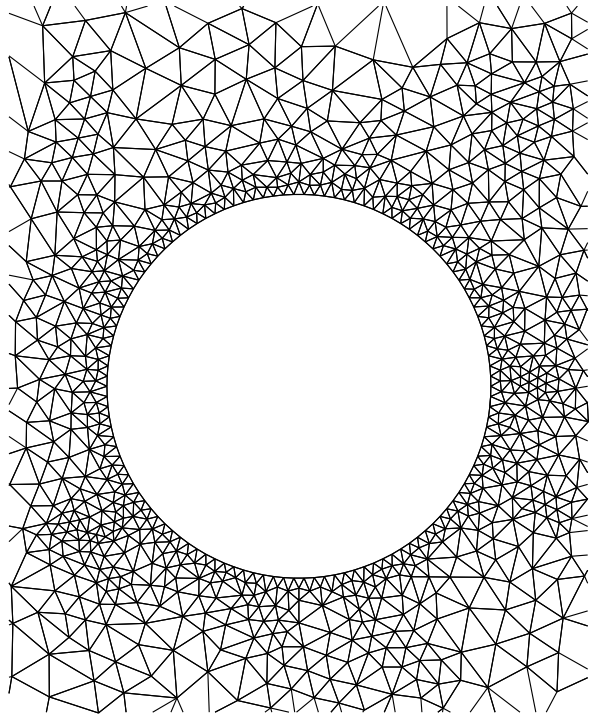}
\caption{Close Up of Grid 1 ($a=.5$, 118 equi-spaced surface points)}
\end{figure}

\begin{figure}[H]
        \centering
                \includegraphics[totalheight=3.5in]{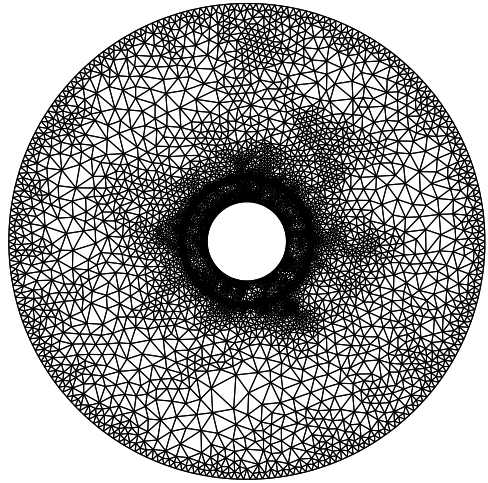}
\caption{Grid 2 (5512 nodes 10548 elements)}
\end{figure}

\begin{minipage}[H]{.35\linewidth}
\centering
\begin{table}[H]

\begin{tabular}{|c|c|c|} \hline
$\Gamma$ & Error & Value\\ \hline \hline
$0$ & ${\|e_u\|}_2$ &  0.071990\\ 
& ${\|e_u\|}_\infty$ & 0.082418\\ 
& ${\|e_v\|}_2$ & 0.023744\\ 
& ${\|e_v\|}_\infty$ &  0.042451\\ \hline
$2 \pi a$ & ${\|e_u\|}_2$ &  0.116207\\ 
& ${\|e_u\|}_\infty$ & 0.177046\\ 
& ${\|e_v\|}_2$ & 0.093906\\ 
& ${\|e_v\|}_\infty$ &  0.121053\\ \hline 
 & ${\|e_u\|}_2$ &  0.175243\\ 
$2 \pi a \sqrt{3}$& ${\|e_u\|}_\infty$ & 0.249885\\ 
& ${\|e_v\|}_2$ & 0.160710\\ 
& ${\|e_v\|}_\infty$ &  0.189222\\ \hline 
$4 \pi a$ & ${\|e_u\|}_2$ &  0.198391\\ 
& ${\|e_u\|}_\infty$ & 0.276553\\ 
& ${\|e_v\|}_2$ & 0.185400\\ 
& ${\|e_v\|}_\infty$ &  0.214179\\ \hline 
$6 \pi a$ & ${\|e_u\|}_2$ &  0.287657\\ 
& ${\|e_u\|}_\infty$ & 0.376061\\ 
& ${\|e_v\|}_2$ & 0.277899\\ 
& ${\|e_v\|}_\infty$ &  0.310920\\ \hline 
\end{tabular}
\caption{Error table (Grid 1)}\label{table3}
\end{table}
\end{minipage}
\begin{minipage}[H]{.35\linewidth}
\centering
\begin{table}[H]
\begin{tabular}{|c|c|c|} \hline
$\Gamma$ & Error & Value\\ \hline \hline
$0$ & ${\|e_u\|}_2$ &  0.028821\\ 
& ${\|e_u\|}_\infty$ & 0.039837\\ 
& ${\|e_v\|}_2$ & 0.010564\\ 
& ${\|e_v\|}_\infty$ & 0.021029\\ \hline
$2 \pi a$ & ${\|e_u\|}_2$ &  0.053179\\ 
& ${\|e_u\|}_\infty$ & 0.082475\\ 
& ${\|e_v\|}_2$ & 0.044475\\ 
& ${\|e_v\|}_\infty$ &  0.054971\\ \hline 
 & ${\|e_u\|}_2$ &  0.078290\\ 
$2 \pi a \sqrt{3}$& ${\|e_u\|}_\infty$ & 0.113686\\ 
& ${\|e_v\|}_2$ & 0.072028\\ 
& ${\|e_v\|}_\infty$ &  0.084053\\ \hline 
$4 \pi a$ & ${\|e_u\|}_2$ &  0.087903\\ 
& ${\|e_u\|}_\infty$ & 0.125113\\ 
& ${\|e_v\|}_2$ & 0.082176\\ 
& ${\|e_v\|}_\infty$ &  0.094758\\ \hline 
$6 \pi a$ & ${\|e_u\|}_2$ &  0.124608\\ 
& ${\|e_u\|}_\infty$ & 0.167751\\ 
& ${\|e_v\|}_2$ & 0.120150\\ 
& ${\|e_v\|}_\infty$ &  0.134961\\ \hline 
\end{tabular}
\caption{Error table (Grid 2)}\label{table4}
\end{table}
\end{minipage}

\begin{table}[H]\centering
\begin{tabular}{|c|c|c|} \hline
$\Gamma$ &$\theta_{analytical}$ & $\theta_{numerical}$\\ \hline 
0 & $0^{\circ},180^{\circ}$ & $0^{\circ},180^{\circ}$\\ \hline
$2 \pi a$ & $30^{\circ}, 150^{\circ}$ & $27.4576^{\circ}, 152.5424^{\circ}$(grid 1)\\ 
& & $28.7395^{\circ}, 149.7479^{\circ}$(grid 2) \\ \hline
$2 \pi a \sqrt{3}$ & $60^{\circ}, 120^{\circ}$ & $54.9152^{\circ}, 125.0847^{\circ}$(grid 1)\\ 
& & $57.4790^{\circ}, 122.5210^{\circ}$(grid 2)\\ \hline
\end{tabular}
\caption{Error table (Grid 1 and 2)}\label{table5}
\end{table}

\begin{figure}
\begin{minipage}[htbp]{.5\linewidth}
\centering
\includegraphics[width=2.8in]{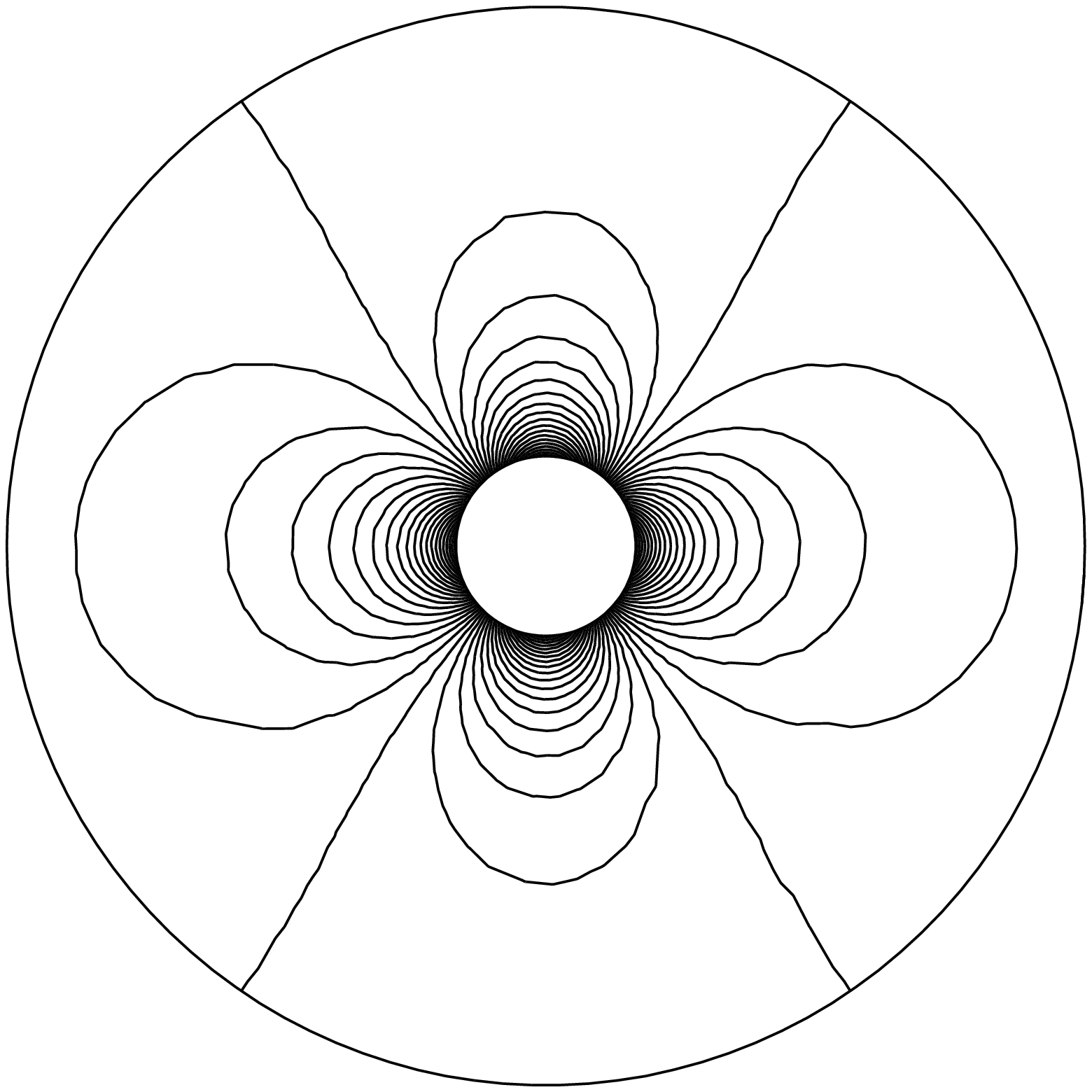}
\caption{Plot of u ($\Gamma=0$)} 
\end{minipage}%
\begin{minipage}[htbp]{.5\linewidth}
\centering
\includegraphics[width=2.8in]{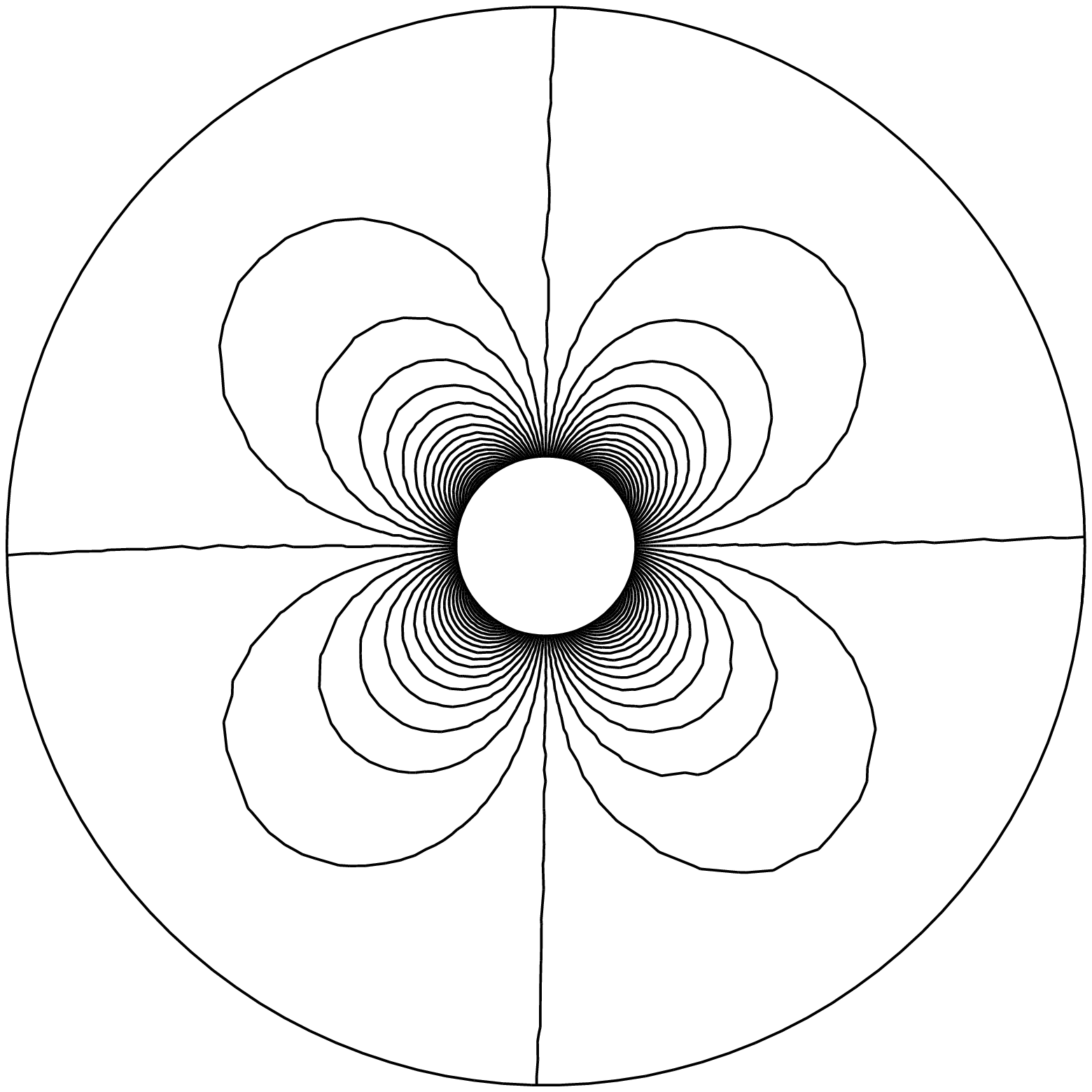}
\caption{Plot of v ($\Gamma=0$)} 
\end{minipage}
\end{figure}

\begin{figure}[htbp]
\centering
\includegraphics[totalheight=3.7in]{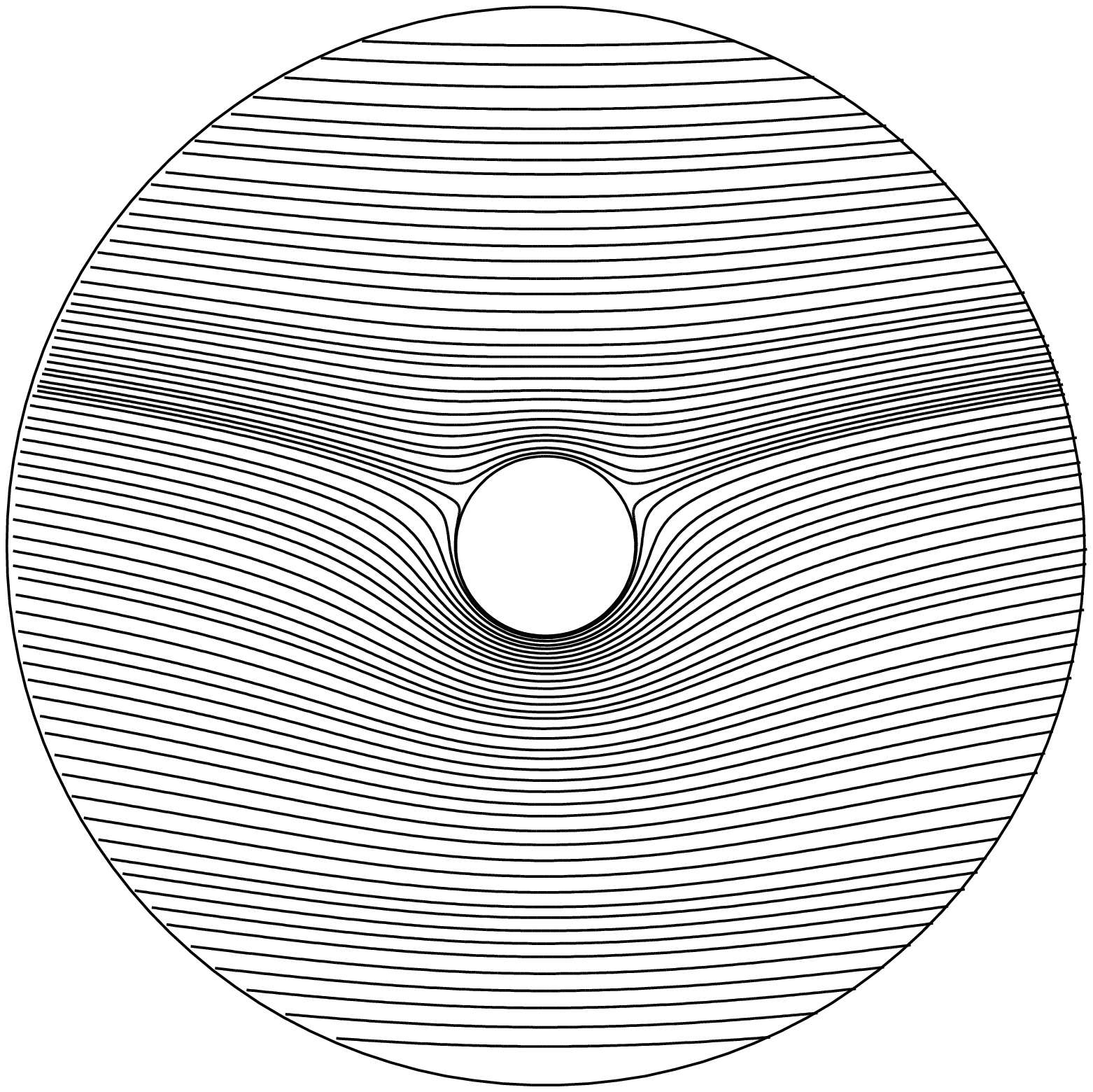}
\caption{Streamlines ($\Gamma=2\pi a$)}
\end{figure}

\begin{figure}[htbp]
\centering
\includegraphics[totalheight=3.7in]{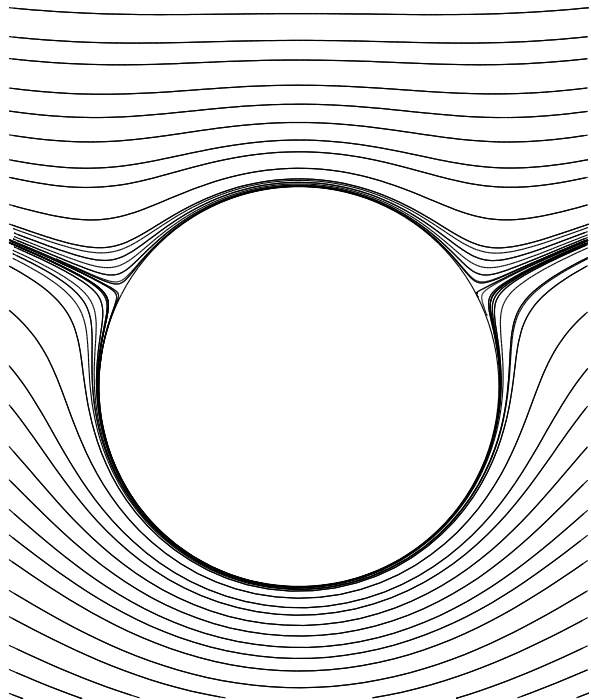}
\caption{Streamlines ($\Gamma=2\pi a$), Stagnation points at $27.4576^{\circ}$ and $152.5424^{\circ}$}
\end{figure}

\begin{figure}[htbp]
        \centering
                \includegraphics[totalheight=3.7in]{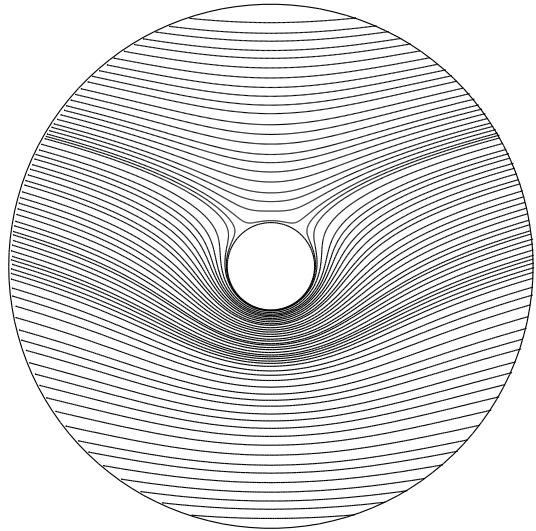}
\caption{Streamlines ($\Gamma=2 \pi a \sqrt{3}$)}
\end{figure}

\begin{figure}[htbp]
        \centering
                \includegraphics[totalheight=3.7in]{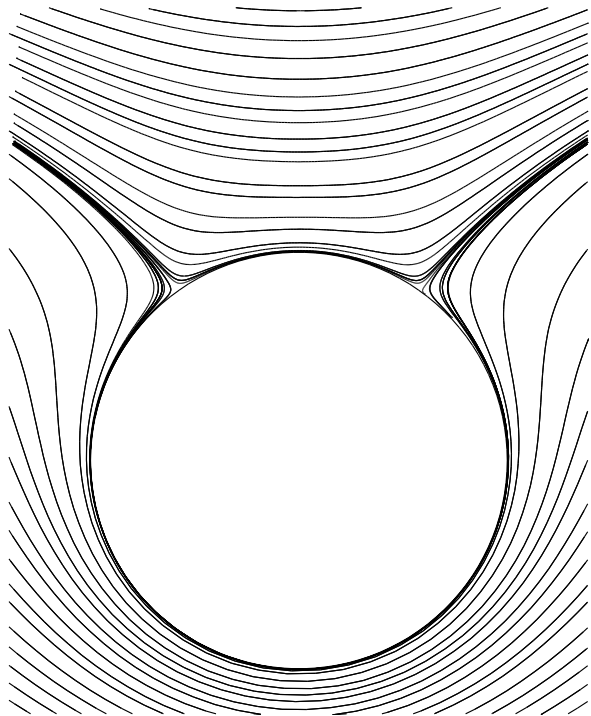}
\caption{Streamlines ($\Gamma=2 \pi a \sqrt{3}$), Stagnation points at $54.9152^{\circ}$ and $125.0847^{\circ}$}
\end{figure}

\begin{figure}[htbp]
        \centering
                \includegraphics[totalheight=3.7in]{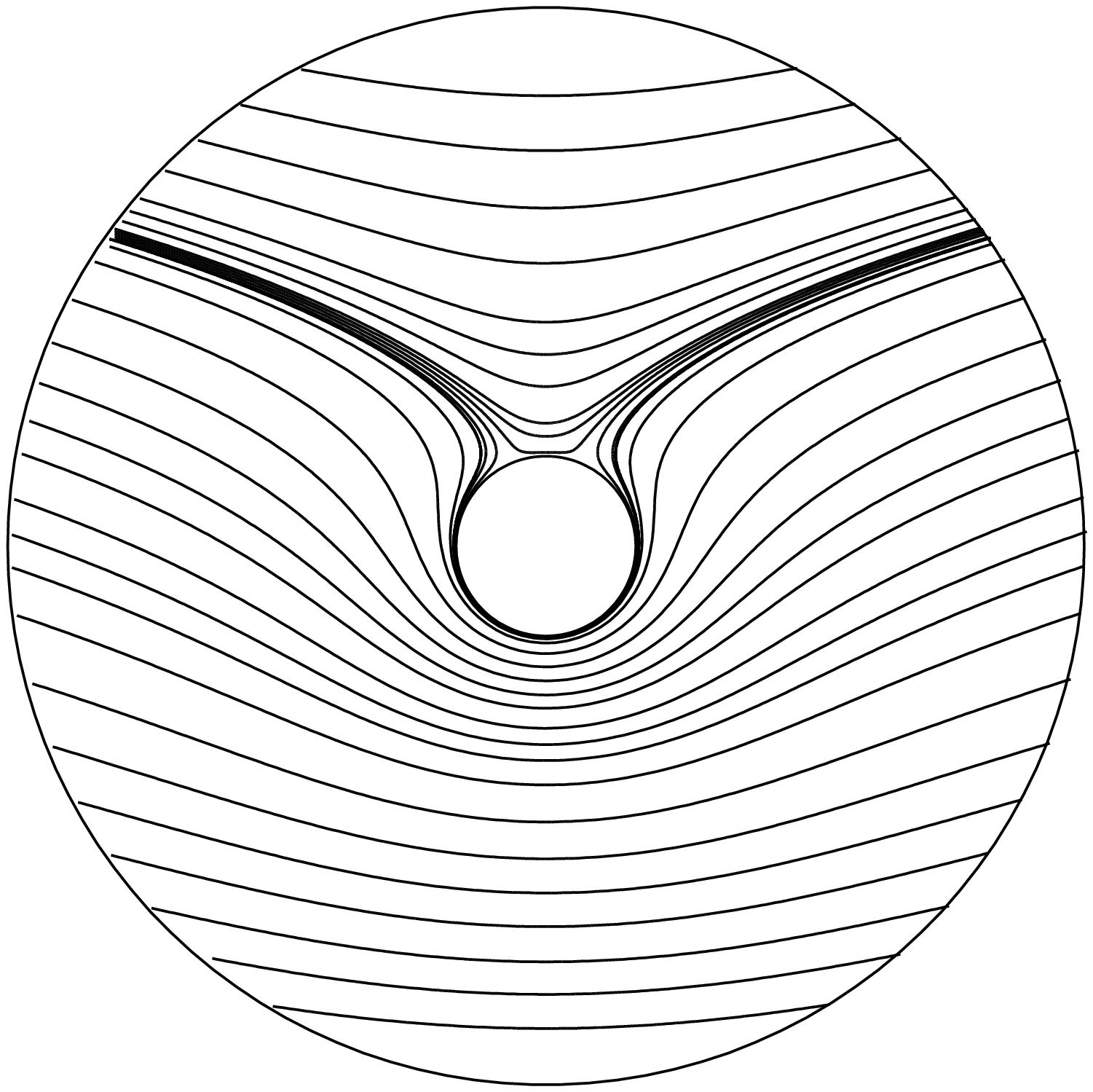}
\caption{Streamlines ($\Gamma=4\pi a$)}
\end{figure}

\begin{figure}[htbp]
        \centering
                \includegraphics[totalheight=3.7in]{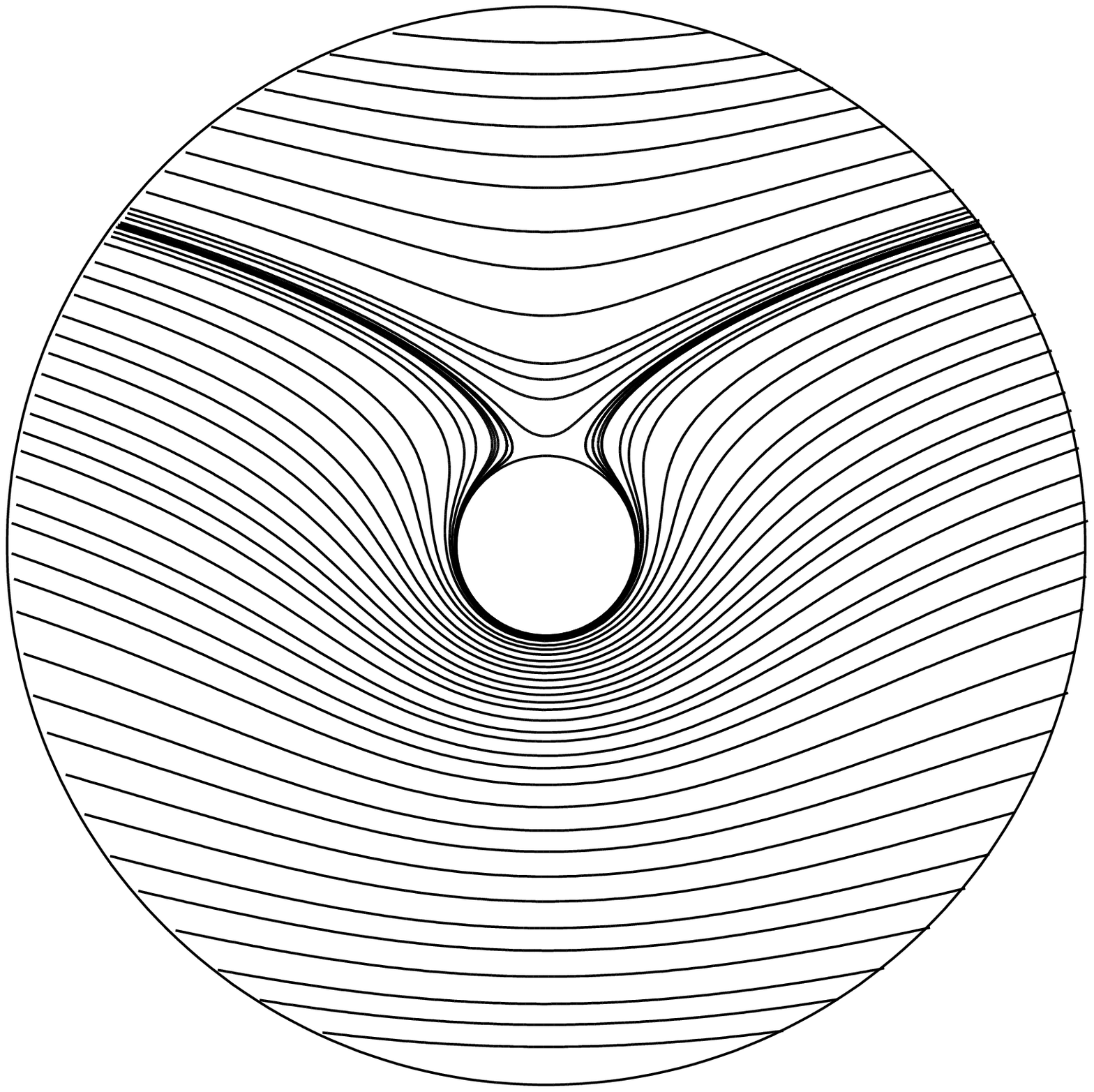}
\caption{Streamlines ($\Gamma=4\pi a$)}
\end{figure}


\section{Inviscid Flow Over An Airfoil}
\begin{figure}[htbp]
        \centering
                \includegraphics[width=6in]{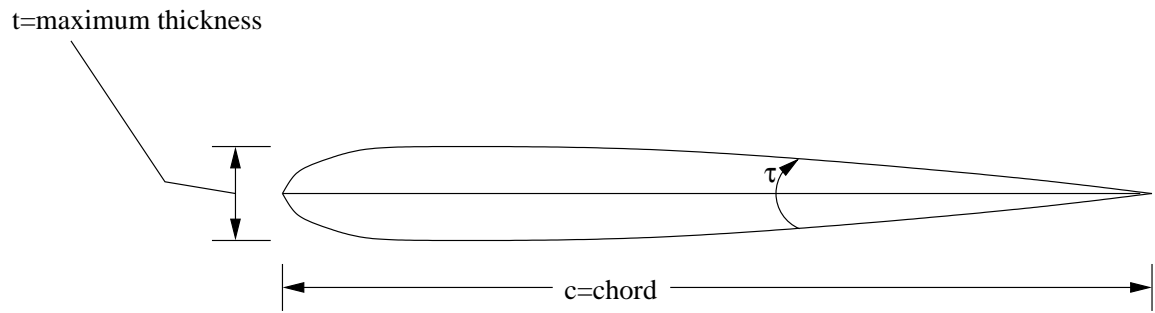}
\caption{airfoil nomenclature}\label{karman}
\end{figure}

In this section the Cauchy-Riemann equations are again used to approximate incompressible inviscid flow over an airfoil.
The treatment of the problem is much the same as in the previous section for flow over the cylinder.  All that is needed
are some modifications to the boundary conditions on the surface and the addition of a Kutta condition.  We again will restrict
the actual computations to airfoils that can be generated using the K\'{a}rm\'{a}n-Trefftz conformal transformation since we have analytical solutions for the flow field. The methodology is nonetheless applicable to an arbitrary airfoil.  However, since flow over an
airfoil is dependent on the flow at the trailing edge, the K\'{a}rm\'{a}n-Trefftz class of airfoils are good for testing the robustness of a numerical method because the velocity is always zero at the trailing edge.  Specifically (refer to figure \ref{karman}), K\'{a}rm\'{a}n-Trefftz airfoil is chosen with a chord length of 1, a maximum thickness to chord ratio of 0.12, a trailing edge angle $\tau=10^{\circ}$ and no camber.
\subsection{Boundary Conditions}
\begin{figure}[htbp]
        \centering
                \includegraphics[width=6in]{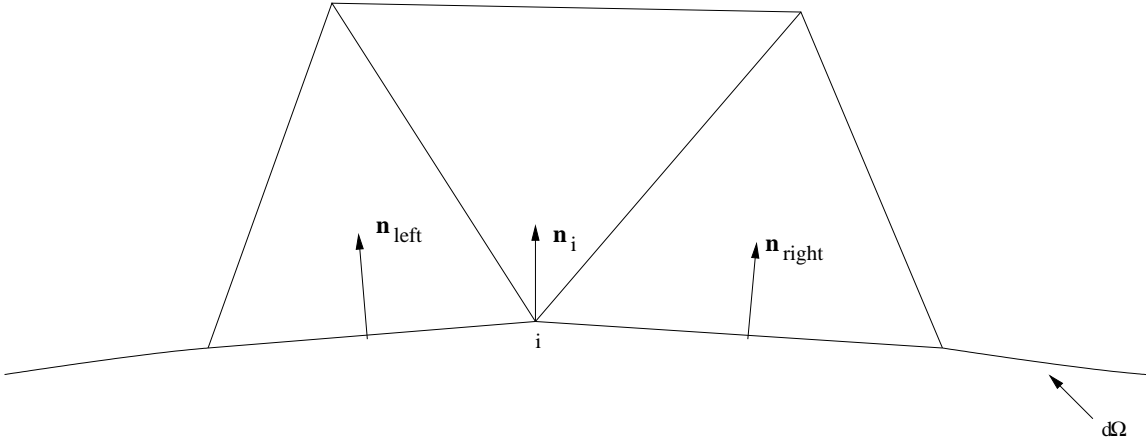}
\caption{Normal vectors on the surface of an airfoil}\label{airnorm}
\end{figure}

\begin{figure}[htbp]
        \centering
                \includegraphics[totalheight=3in]{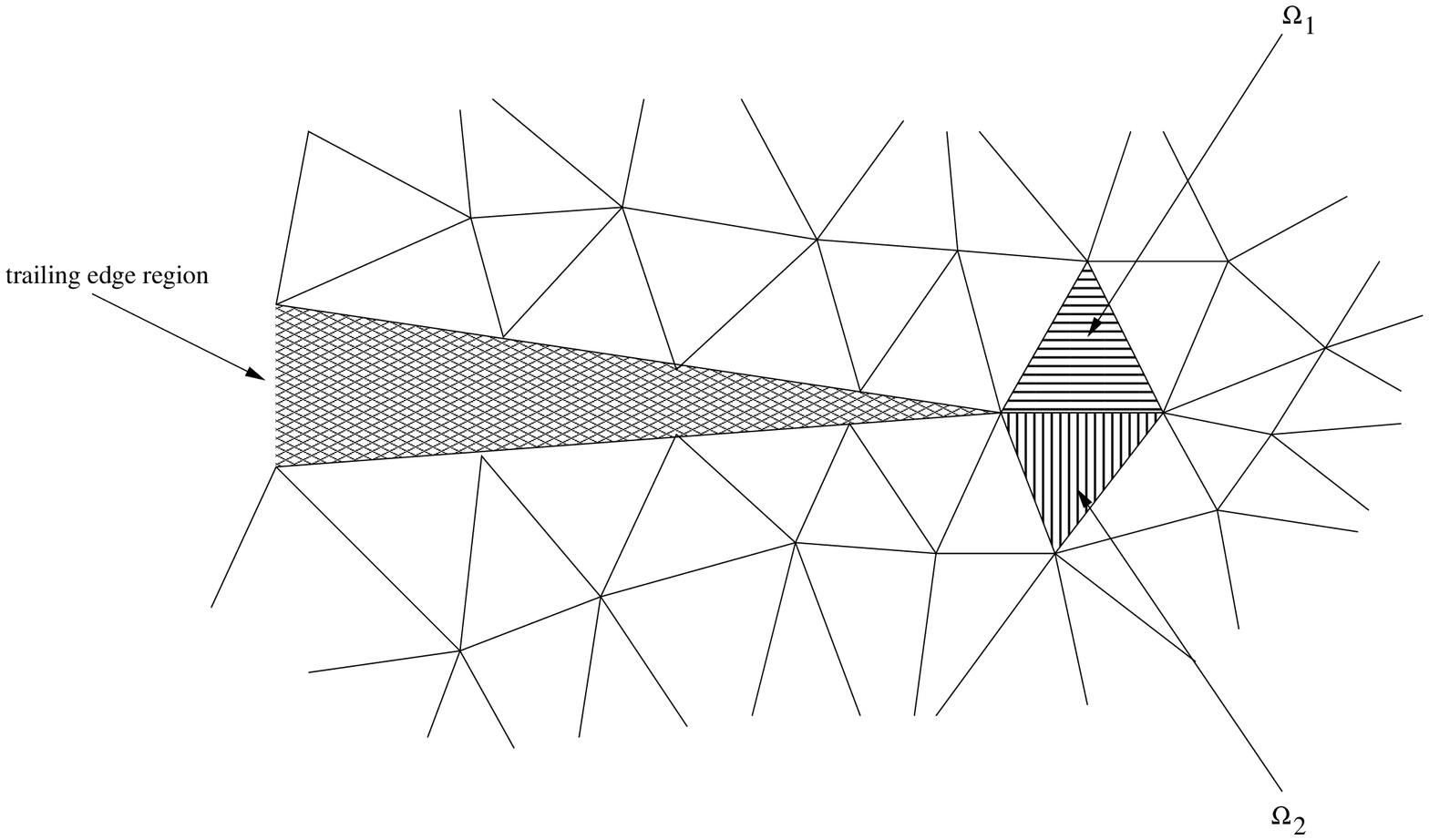}
                 \caption{Trailing Edge\label{kutta2}}
\end{figure}

\noindent \underline{Farfield}\\
In the farfield equations (\ref{cyl1})-(\ref{cyl2}) are enforced where
\begin{align}
u_{exact}&= \cos{\alpha}+\frac{\Gamma}{2 \pi r^2}\\
v_{exact}&= \sin{\alpha}-\frac{\Gamma}{2 \pi r^2}
\end{align}
which are readily produced from the stream function of equation (\ref{streamfunction}).  As before, $\alpha$ is the angle of attack and $\Gamma$ is the circulation given by,
\begin{equation}
\Gamma=4\pi a \sin{\alpha}
\end{equation}
where $a$ is the radius of the cylinder used to create the airfoil via the K\'{a}rm\'{a}n-Trefftz transformation.  For all airfoils in
this section $a=0.273094$. \\
\underline{Airfoil's Surface}\\
For nodes on the surface, in addition to enforcing equations (\ref{didu})- (\ref{didvv}) we again enforce the tangency condition $\mathbf u \cdot \mathbf n =0$ via equations (\ref{norm1})- (\ref{norm11}).  However, unlike the cylinder in general there is no analytic
expression for the outward normal vector on the surface of an airfoil.  Consider a node $i$ on the surface of an airfoil as depicted in figure \ref{airnorm}.  In order
to approximate the outward normal vector $\mathbf n_i$ at the node, first the normal vectors are calculated to the left and right of the node corresponding to the sides adjacent to the surface of triangles with the common vertex node $i$.  Then the normal at node $i$ is
calculated from the arithmetic average of these two vectors.\\
\underline{Kutta Condition}\\
With reference to figure \ref{kutta2}, the Kutta condition is simply
\begin{equation}
(\mathbf u_1- \mathbf u_2) \cdot \mathbf n=0
\end{equation}
where $\mathbf u_1$ and $\mathbf u_2$ are the average velocities in $\Omega_1$ and $\Omega_2$ respectively.  The normal vector
$\mathbf n$ is taken as the normal vector of the common edge of the two shaded elements in the figure.  So, in addition to enforcing
equations (\ref{didu})- (\ref{didvv}) on these two elements we need to minimize,
\begin{equation}
\mathcal{I}_{kutta}=\frac{((\mathbf u_1- \mathbf u_2)\cdot \mathbf n)^2}{2}
\end{equation}
We only need to minimize at the nodes that aren't on the interface of $\Omega_1$ and $\Omega_2$. Assume nodes $i$ and $j$ are such nodes of the respective elements $\Omega_1$ and $\Omega_2$.
Then, at node $i$ we have,
\begin{align}
\frac{\partial\mathcal{I}_{kutta}}{\partial u_i}&=\frac{n_x}{3}(\mathbf u_1- \mathbf u_2)\cdot \mathbf n\\
\frac{\partial ^2\mathcal{I}_{kutta}}{\partial u_i \partial u_i}
&=\frac{n_x^2}{9}\\
\frac{\partial ^2\mathcal{I}_{kutta}}{\partial u_j \partial u_i}
&=-\frac{n_x^2}{9}\\
\frac{\partial ^2\mathcal{I}_{kutta}}{\partial v_i \partial u_i}
&=\frac{n_x n_y}{9}\\
\frac{\partial ^2\mathcal{I}_{kutta}}{\partial v_j \partial u_i}
&=-\frac{n_x n_y}{9}
\end{align}
\hrulefill\\ 
\begin{align}
\frac{\partial\mathcal{I}_{kutta}}{\partial v_i}&=\frac{n_y}{3}(\mathbf u_1- \mathbf u_2)\cdot \mathbf n\\
\frac{\partial ^2\mathcal{I}_{kutta}}{\partial u_i \partial v_i}
&=\frac{n_x n_y}{9}\\
\frac{\partial ^2\mathcal{I}_{kutta}}{\partial u_j \partial v_i}
&=-\frac{n_x n_y}{9}\\
\frac{\partial ^2\mathcal{I}_{kutta}}{\partial v_i \partial v_i}
&=\frac{n_y^2}{9}\\
\frac{\partial ^2\mathcal{I}_{kutta}}{\partial v_j \partial j_i}
&=-\frac{n_y^2}{9}
\end{align}
We get similar expressions for $\frac{\partial\mathcal{I}_{kutta}}{\partial u_j}$ and 
$\frac{\partial\mathcal{I}_{kutta}}{\partial v_j}$.
\subsection{Results}
\begin{figure}[htbp]
\begin{minipage}[t]{.5\linewidth}
\centering
\includegraphics[width=2.7in]{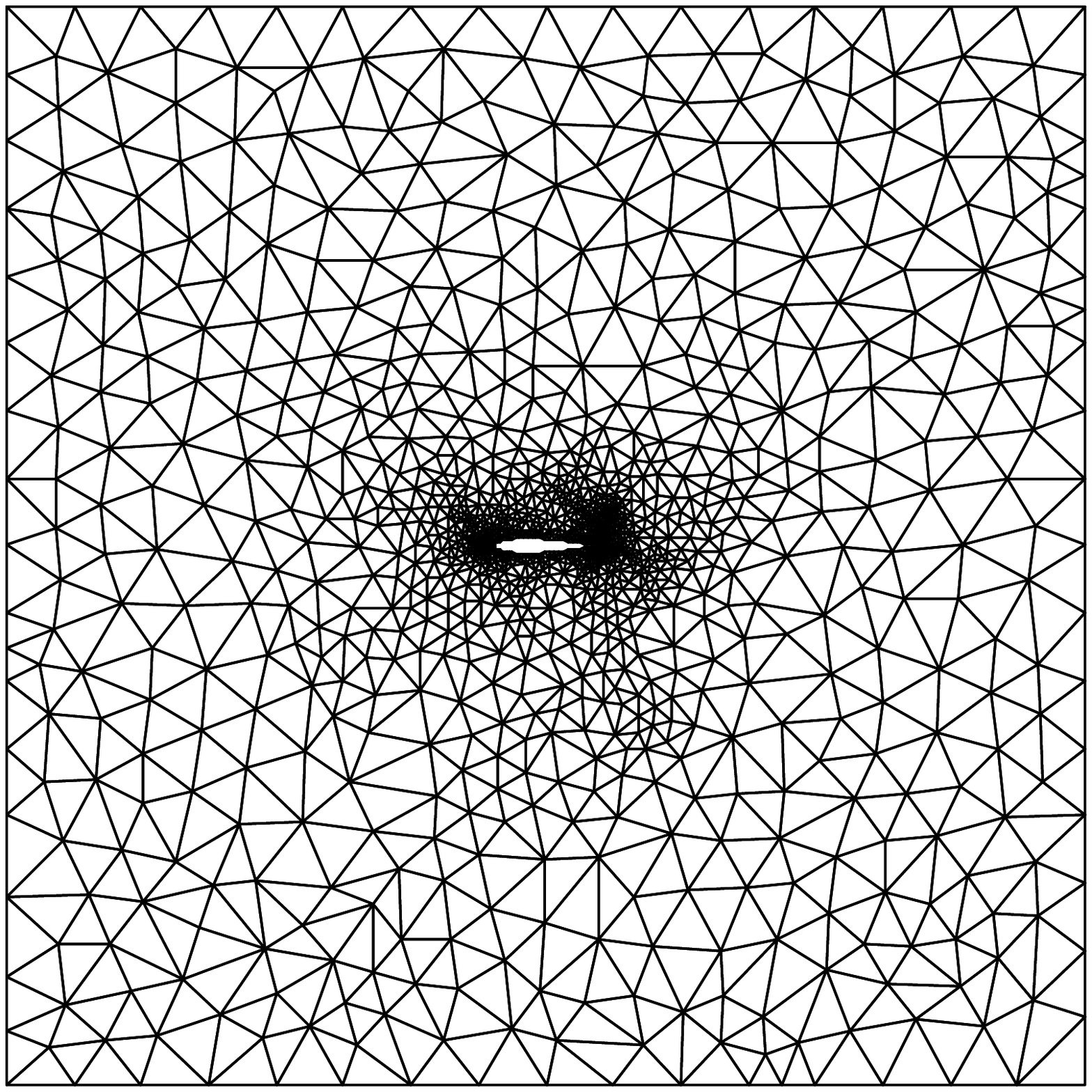}
         \caption{Unstructured Grid About Airfoil\label{airgrid1}}
\end{minipage}%
\begin{minipage}[t]{.5\linewidth}
\centering
\includegraphics[width=2.7in]{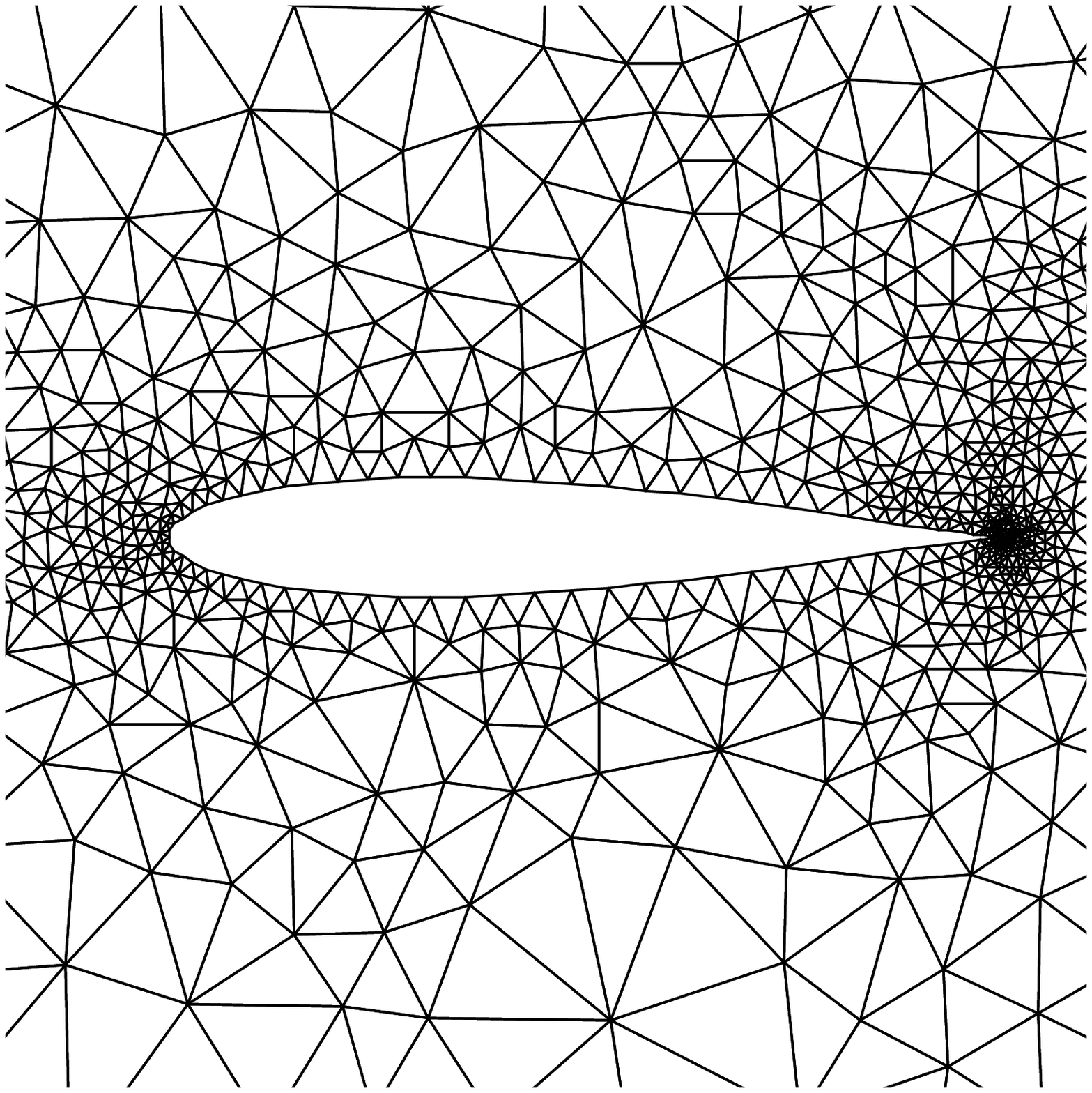}
         \caption{Close-up: 70 surface points\label{airgrid2}}
\end{minipage}
\end{figure}

\begin{figure}[htbp]
        \centering
                \includegraphics[totalheight=3.7in]{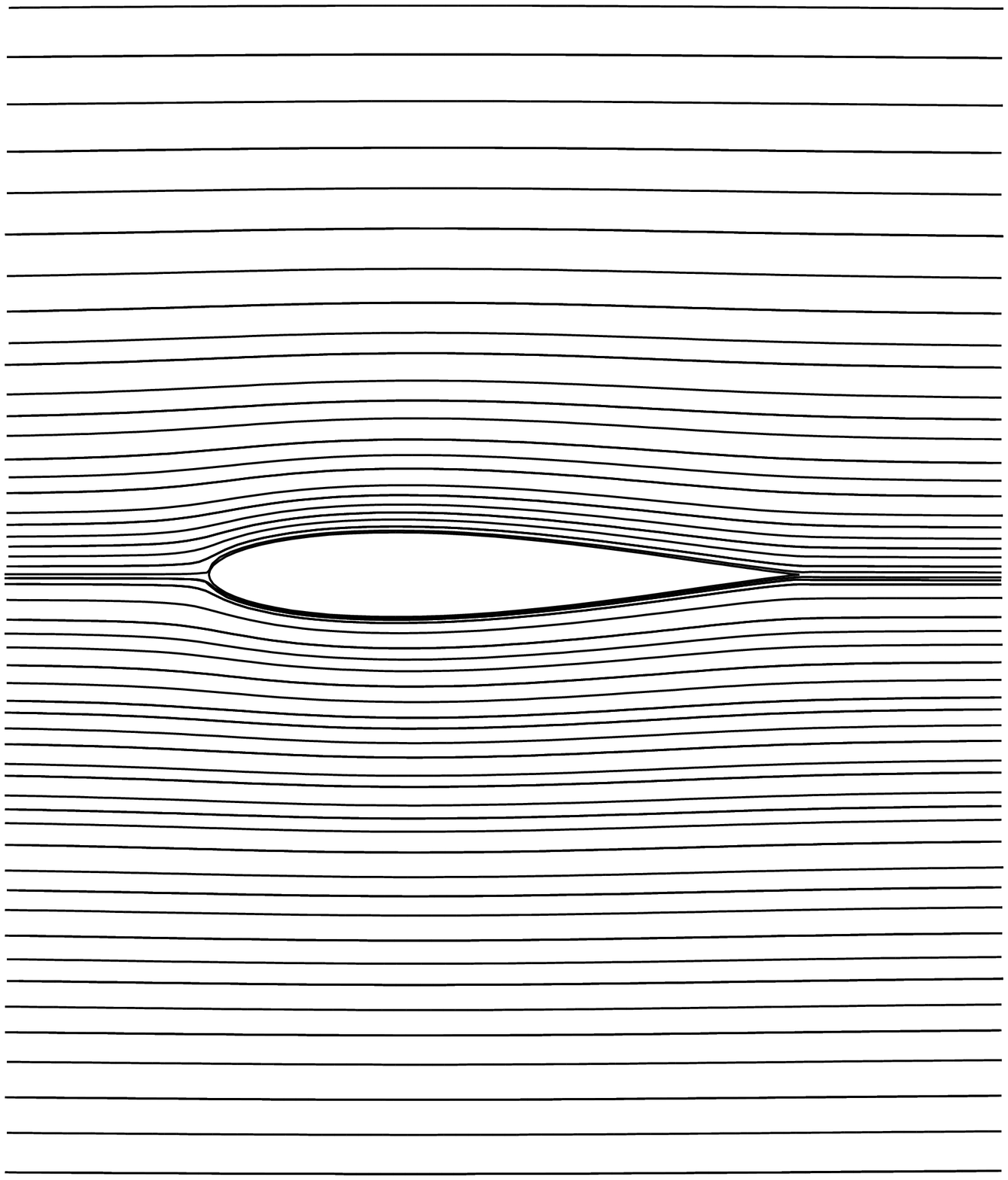}
\caption{Streamlines ($\alpha=0^{\circ}$, $t=.12$, $\tau=10^{\circ}$) \label{airs1}}
\end{figure}

\begin{figure}[htbp]
\begin{minipage}[t]{.5\linewidth}
\centering
\includegraphics[width=2.7in]{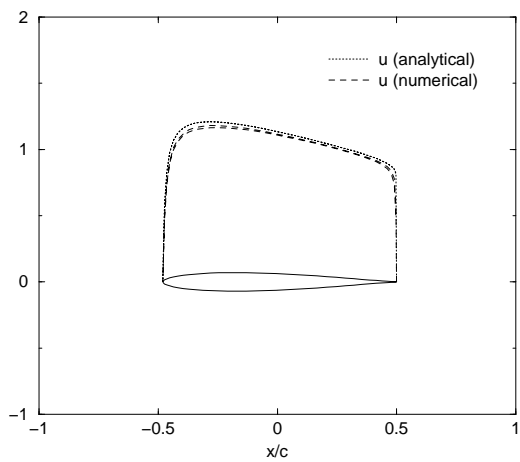}
         \caption{u-profile about airfoil ($\alpha=0$) \label{airp1}}
\end{minipage}%
\begin{minipage}[t]{.5\linewidth}
\centering
\includegraphics[width=2.7in]{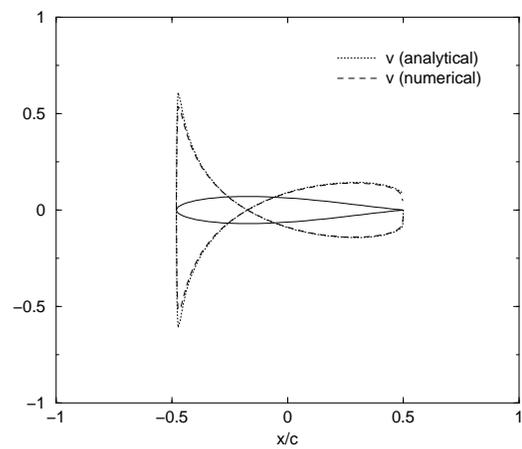}
         \caption{v-profile about airfoil ($\alpha=0$) \label{airp2}}
\end{minipage}
\end{figure}

\begin{figure}[htbp]
        \centering
                \includegraphics[totalheight=3.7in]{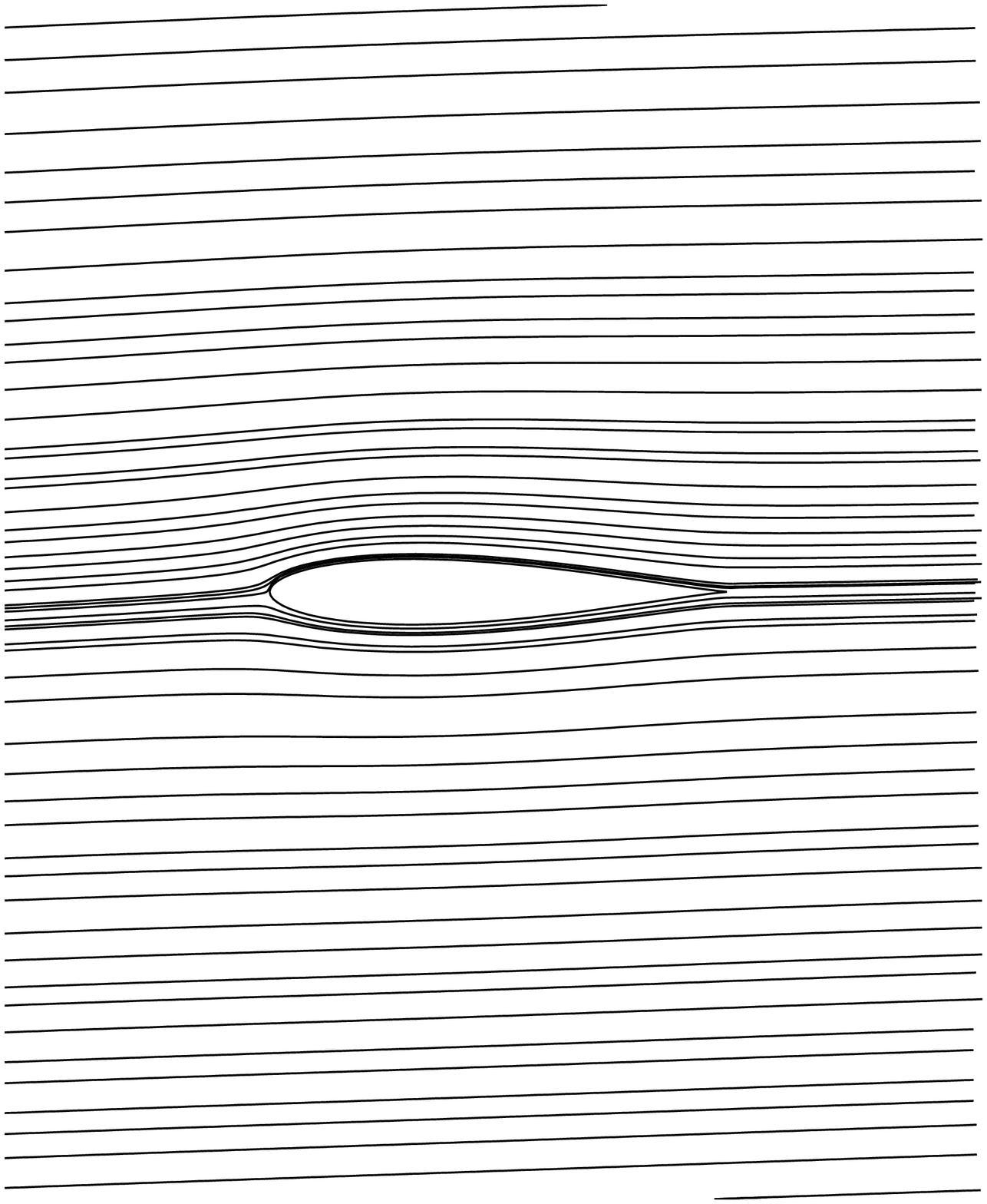}
\caption{Streamlines ($\alpha=2^{\circ}$, $t=.12$, $\tau=10^{\circ}$) \label{airs2}}
\end{figure}

\begin{figure}[htbp]
\begin{minipage}[t]{.5\linewidth}
\centering
\includegraphics[width=2.7in]{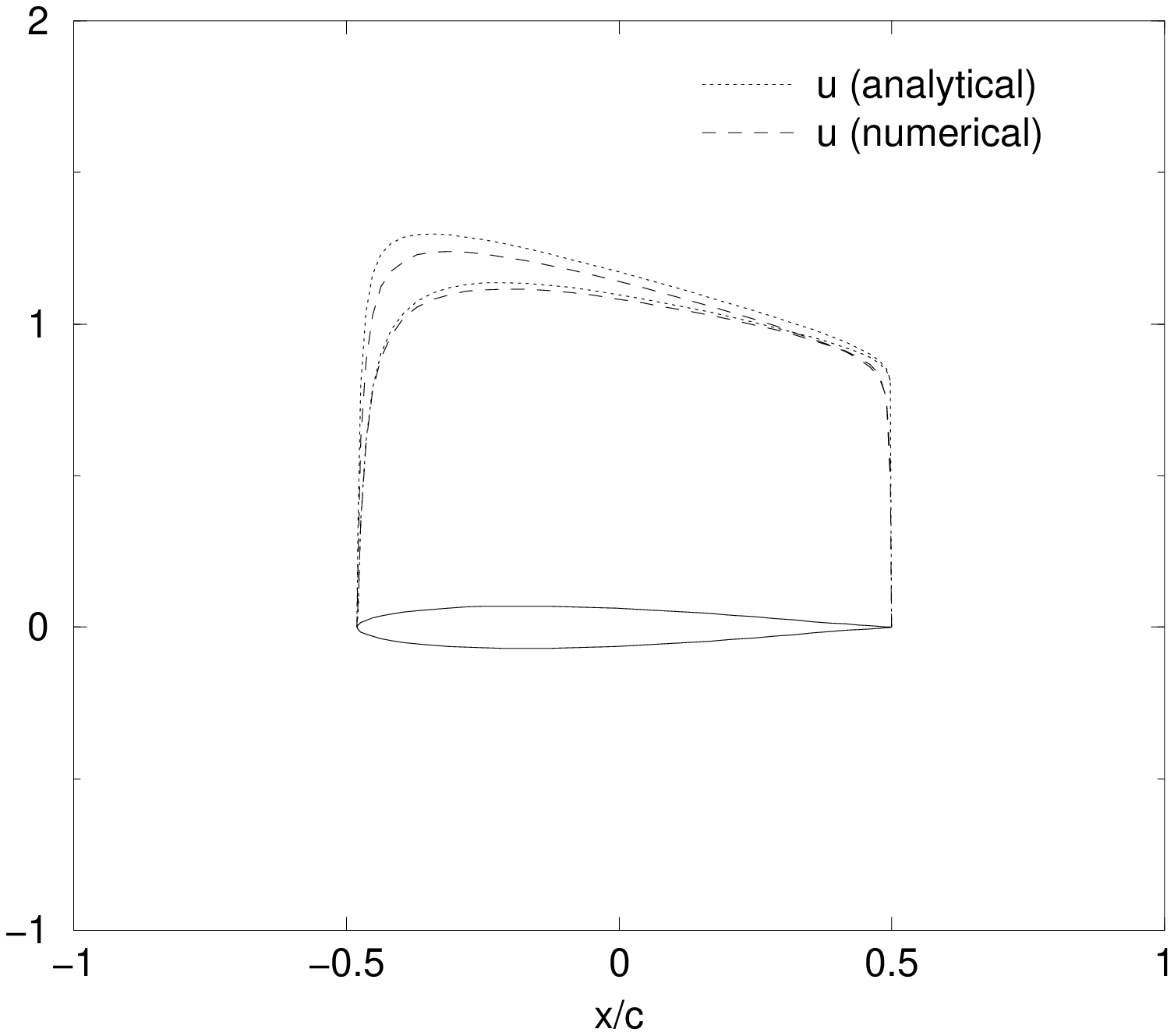}
         \caption{u-profile about airfoil ($\alpha=2^{\circ}$) \label{airp3}}
\end{minipage}%
\begin{minipage}[t]{.5\linewidth}
\centering
\includegraphics[width=2.7in]{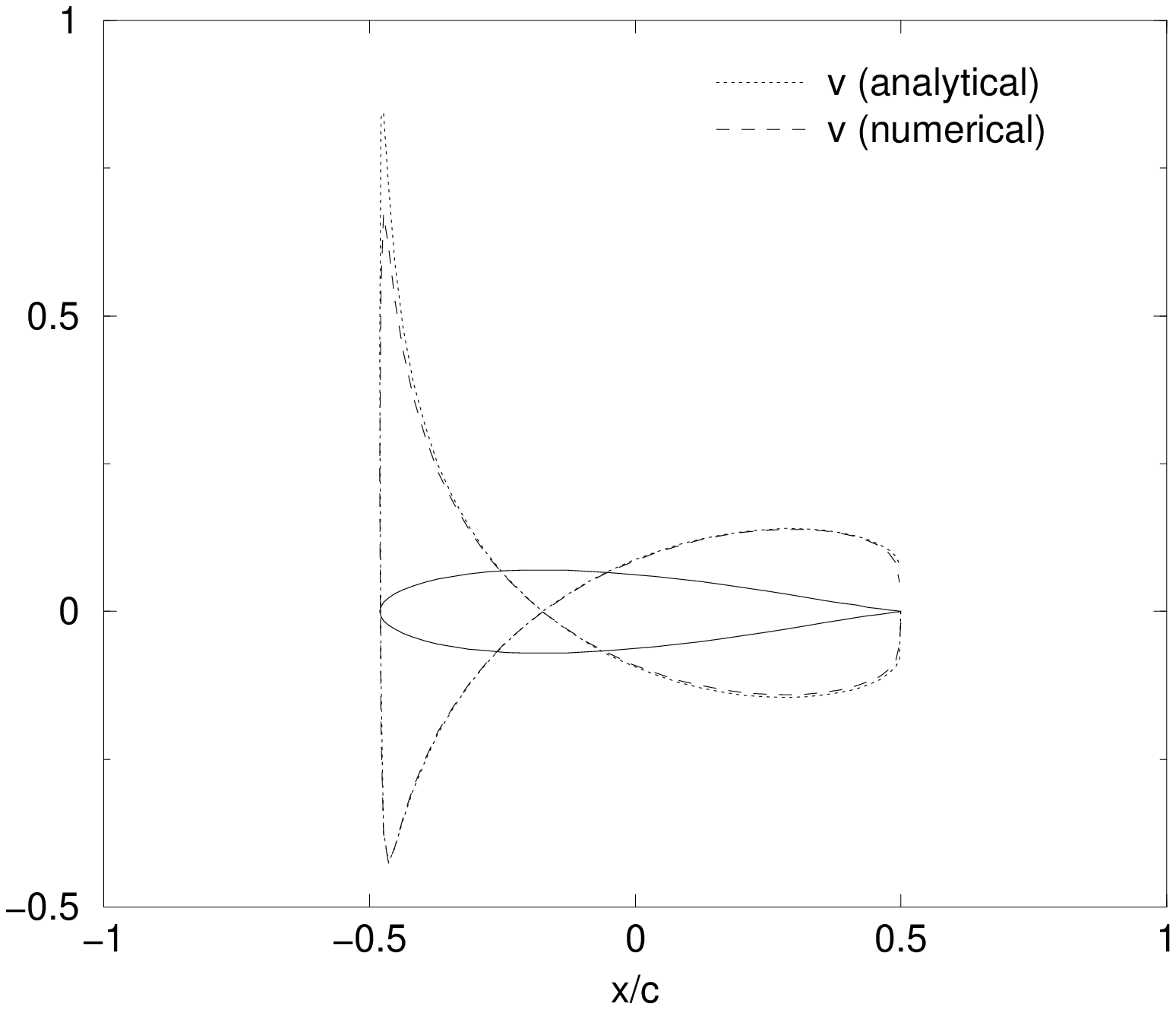}
         \caption{v-profile about airfoil ($\alpha=2^{\circ}$) \label{airp4}}
\end{minipage}
\end{figure}

Figures \ref{airgrid1} and \ref{airgrid2} depict the grid used in the computations.
The airfoil is contained in a $[-5,5]\times[-5,5]$ box.
Figures \ref{airs1}- \ref{airp2} show plots for the case $\alpha=0$.  The results are
in good agreement with the analytical solutions.  
Figures \ref{airs2}- \ref{airp4} show plots for the case $\alpha=2^{\circ}$.  The results
show the method is able to handle flow over an airfoil with circulation however, it is surprising
that the error in $u$ is already becoming appreciable at a very small angle of attack.  The error
in $v$ is mostly concentrated at the leading edge.  Otherwise, it matches up well with the analytical
solution.  Adding more points on the surface of the airfoil and refining the mesh did little to improve
the solution.  In the process of reaching these solutions, different formulations of the Kutta condition were enforced.
Indeed, improving the solution significantly will probably entail using a formulation of the Kutta condition that is more
conducive to the least-squares method. 
\chapter{Boundary-Layer Flow}
\section{Formulation and Derivations}
In this chapter, the least-squares method is applied to a nonlinear set of governing equations in an effort to calculate
boundary-layer flow over a flat plate.  For high Reynolds regime it is known that the Prandtl boundary-layer equations coupled with
potential flow can be
used to approximate the flow as opposed to solving the Navier-Stokes equations.  The equations are parabolic and simpler to solve.
\paragraph{} 
In dimensionless form, the Prandtl boundary-layer equations are:
\begin{eqnarray}
u_x+v_y&=&0\\
uu_x+vu_y-\omega_y-\mathcal{U}\mathcal{U}_x&=&0\\
u_y-\omega&=&0
\end{eqnarray}
The dimensionless variables are related to the dimensional variables by:
\begin{eqnarray*}
x'&=&Lx\\
y'&=&\delta y\\
v'&=&\frac{U_\infty\delta v}{L}\\
\omega'&=&\frac{U_\infty\omega}{\delta}\\
\mathcal{U}'&=&U_\infty\mathcal{U}
\end{eqnarray*}
where $\mathcal{U}'$ is the parallel component of velocity at the outer edge of the boundary layer,
L is the characteristic length of the flate plate and $\delta=\frac{L}{\sqrt{\mathcal{R}e}}$.
\subsection{Functional Minimization}
The functional for the boundary-layer equations is
\begin{eqnarray} 
\mathcal{I}_{_T}&=&\frac{1}{2}\left((u_x+v_y)^2
     +(\bar{u}u_x+\bar{v}u_y-\omega_y-\mathcal{U}\mathcal{U}_x)^2+(u_y-\bar{\omega})^2\right)\Omega_{_T}
\end{eqnarray} 
where $M$ is the number of elements in the domain.  On each element, for each node $i$ we have,
\begin{eqnarray}
\frac{\partial \mathcal{I}_{_T}}{\partial u_i}&=&\left((\bar{u}u_x+\bar{v}u_y-\omega_y
-\bar{\mathcal{U}}\mathcal{U}_x)(\frac{u_x}{3}-\frac{1}{2\Omega_{_T}}(\bar{u}n_{x_i}+\bar{v}n_{y_i}))\right.\label{bdy1}\\ 
& &\left.-\frac{1}{2\Omega_{_T}}\left((u_y-\bar{\omega})n_{y_i}+(u_x + v_y)n_{x_i}\right)\right)\Omega_{_T}\nonumber\\
\frac{\partial^2 \mathcal{I}_{_T}}{\partial u_j\partial u_i}&=&\left(\frac{\vec{n}_i\cdot\vec{n}_j}{4\Omega^2_{_T}}-(\bar{u}u_x+\bar{v}u_y-\omega_y
-\bar{\mathcal{U}}\mathcal{U}_x)(\frac{n_{x_j}+n_{x_i}}{6\Omega_{_T}})\right.\\
& &+\left.(\frac{u_x}{3}-\frac{1}{2\Omega_{_T}}(\bar{u}n_{x_i}+\bar{v}n_{y_i}))(\frac{u_x}{3}-\frac{1}{2\Omega_{_T}}(\bar{u}n_{x_j}+\bar{v}n_{y_j}))\right)\Omega_{_T}\nonumber\\
\frac{\partial^2 \mathcal{I}}{\partial v_j\partial u_i}&=&\left(\frac{n_{x_i}n_{y_j}}{4\Omega^2_{_T}}+(\bar{u}u_x+\bar{v}u_y-\omega_y
-\bar{\mathcal{U}}\mathcal{U}_x)\frac{n_{y_i}}{6\Omega_{_T}}\right.\nonumber\\
& &\left.+\frac{u_y}{3}(\frac{u_x}{3}-\frac{1}{2\Omega_{_T}}(\bar{u}n_{x_i}+\bar{v}n_{y_i}))\right)\Omega_{_T}\\
\frac{\partial^2 \mathcal{I}_{_T}}{\partial\omega_j\partial u_i}&=&\left((\frac{u_x}{3}-\frac{\bar{u}n_{x_i}+\bar{v}n_{y_i}}{2\Omega_{_T}})\frac{n_{y_j}}{2\Omega_{_T}}+\frac{n_{y_i}}{6\Omega_{_T}}\right)\Omega_{_T} \label{bdy12}
\end{eqnarray}
\hrulefill\\ 
\begin{align}
\frac{\partial \mathcal{I}_{_T}}{\partial v_i}&=\left(-\frac{1}{2\Omega_{_T}}(u_x + v_y)n_{y_i}+\frac{u_y}{3}(\bar{u}u_x+\bar{v}u_y-\omega_y
-\bar{\mathcal{U}}\mathcal{U}_x)\right)\Omega_{_T}\label{bdy21}\\
\frac{\partial^2 \mathcal{I}_{_T}}{\partial u_j\partial v_i}&=\left(\frac{n_{x_j}n_{y_i}}{4\Omega^2_{_T}}+(\bar{u}u_x+\bar{v}u_y-\omega_y
-\bar{\mathcal{U}}\mathcal{U}_x)\frac{n_{y_j}}{6\Omega_{_T}}\right.\nonumber\\
&\quad +\left.\frac{u_y}{3}(\frac{u_x}{3}-\frac{\bar{u}n_{x_j}+\bar{v}n_{y_j}}{2\Omega_{_T}})\right)\Omega_{_T}\\
\frac{\partial^2 \mathcal{I}_{_T}}{\partial v_j\partial v_i}&=(\frac{n_{y_i}n_{y_j}}{4\Omega^2_{_T}}+\frac{(u_y)^2}{9})\Omega_{_T}\\
\frac{\partial^2 \mathcal{I}_{_T}}{\partial w_j\partial v_i}&=-n_{y_j}\frac{u_y}{6}\label{bdy22}
\end{align} 
\hrulefill\\
\begin{align}
\frac{\partial \mathcal{I}_{_T}}{\partial\omega_i}&=\left((\bar{u}u_x+\bar{v}\bar{u_y}-\omega_y
-\bar{\mathcal{U}}\mathcal{U}_x)\frac{n_{y_i}}{2\Omega_{_T}}-\frac{u_y-\bar{\omega}}{3}\right)\Omega_{_T} \label{bdy31}\\
\frac{\partial^2 \mathcal{I}_{_T}}{\partial u_j\partial\omega_i}&=-\left((\frac{u_x}{3}+\bar{u}n_{x_j}+\bar{v}n_{y_j})n_{y_i}+\frac{n_{y_j}}{6\Omega_{_T}}\right)\Omega_{_T}\\
\frac{\partial^2 \mathcal{I}_{_T}}{\partial v_j\partial\omega_i}&=n_{y_i}\frac{u_y}{6}\\
\frac{\partial^2 \mathcal{I}_{_T}}{\partial \omega_j\partial\omega_i}&=(\frac{1}{9}+\frac{n_{y_i}n_{y_j}}{4\Omega^2_{_T}})\Omega_{_T}\label{bdy33}
\end{align}
\subsection{Boundary Conditions}
\begin{figure}[htbp]
        \centering
         \includegraphics[totalheight=3in]{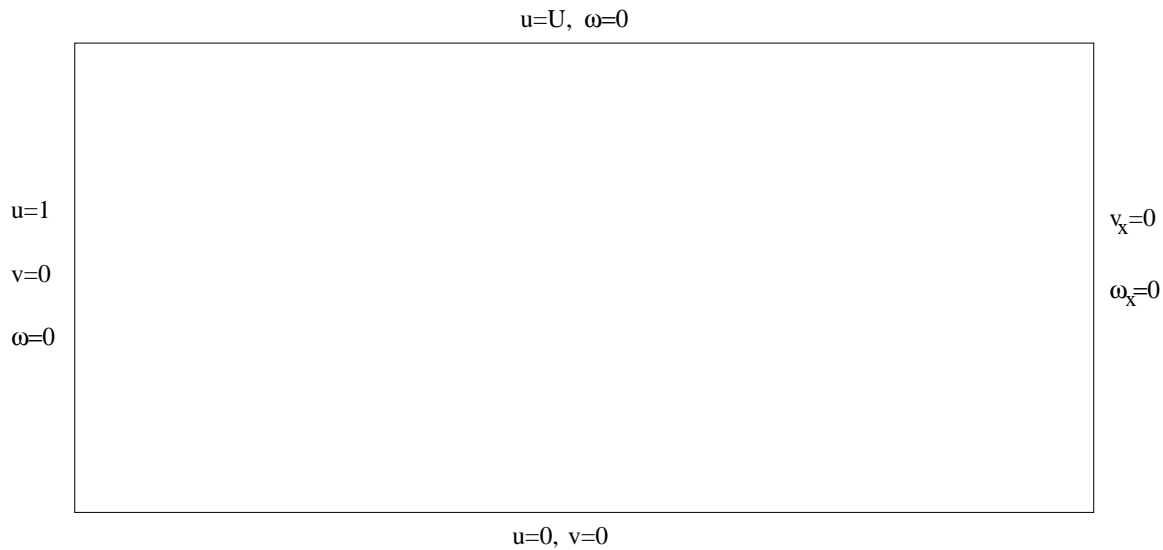}
         \caption{computational domain and boundary conditions}
\end{figure}

\underline{Bottom}\\\\
If the index $i$ corresponds to a boundary node we have $u=0$, $v=0$.  Consequently, instead of using equations (\ref{bdy1})-(\ref{bdy33}) we have,

\begin{align}
\mathcal{I}_i&=\frac{(u_i-1)^2+v_i^2}{2}\\
\frac{\partial \mathcal{I}_i}{\partial u_i}&=u_i\\
\frac{\partial^2 \mathcal{I}_i}{\partial u_i\partial u_i}&=1\\
\frac{\partial \mathcal{I}_i}{\partial v_i}&=v_i\\
\frac{\partial^2 \mathcal{I}_i}{\partial v_i\partial v_i}&=1\\
\end{align}
For $\omega$, we enforce equations (\ref{bdy31})-(\ref{bdy33}).\\\\
\underline{Left     }\\\\
Here we minimize the equation $\mathcal{I}_i=\frac{(u-1)^2+v^2+\omega^2}{2}$.  This gives,
\begin{align}
\frac{\partial \mathcal{I}_i}{\partial u_i}&=u_i-1\\
\frac{\partial^2 \mathcal{I}_i}{\partial u_i\partial u_i}&=1\\
\frac{\partial \mathcal{I}_i}{\partial v_i}&=v_i\\
\frac{\partial^2 \mathcal{I}_i}{\partial v_i\partial v_i}&=1\\
\frac{\partial \mathcal{I}_i}{\partial\omega_i}&=\omega_i\label{om1}\\
\frac{\partial^2 \mathcal{I}_i}{\partial \omega_i\partial\omega_i}&=1\label{om2}
\end{align}
\underline{Top     }\\\\
On the top boundary for $u$ and $\omega$ me minimize $\frac{(u-\mathcal{U})^2+\omega^2}{2}=0$ whereby we get,
\begin{align}
\frac{\partial \mathcal{I}_i}{\partial u_i}&=u_i-\mathcal{U}_i\\
\frac{\partial^2 \mathcal{I}_i}{\partial u_i\partial u_i}&=1
\end{align}
For $\omega$ we again enforce equations (\ref{om1})-(\ref{om2}), and for $v$ we use equations (\ref{bdy21})-(\ref{bdy22}).\\\\
\underline{Right   }\\\\
On the right side for $u$ we simply enforce equations (\ref{bdy1})-(\ref{bdy12}).  For $\omega$ and $v$ we minimize 
\begin{equation*}
\mathcal{I}_i=\frac{1}{2}\left((v_x)^2+(\omega_x)^2\right)\Omega_{_T}
\end{equation*}
which gives us,
\begin{align}
\frac{\partial \mathcal{I}_i}{\partial v_i}&=-\frac{v_xn_{x_i}}{2}\\
\frac{\partial^2 \mathcal{I}_i}{\partial v_j\partial v_i}&=\frac{n_{x_i}n_{x_j}}{4\Omega_{_T}}\\
\frac{\partial \mathcal{I}_i}{\partial v_i}&=\frac{\omega_xn_{x_i}}{4\Omega_{_T}}\\
\frac{\partial^2 \mathcal{I}_i}{\partial \omega_j\partial \omega_i}&=\frac{n_{x_i}n_{x_j}}{4\Omega_{_T}}
\end{align}

\subsection{Additional Computational Aspects}
\begin{figure}[htbp]
        \centering
         \includegraphics[totalheight=4in]{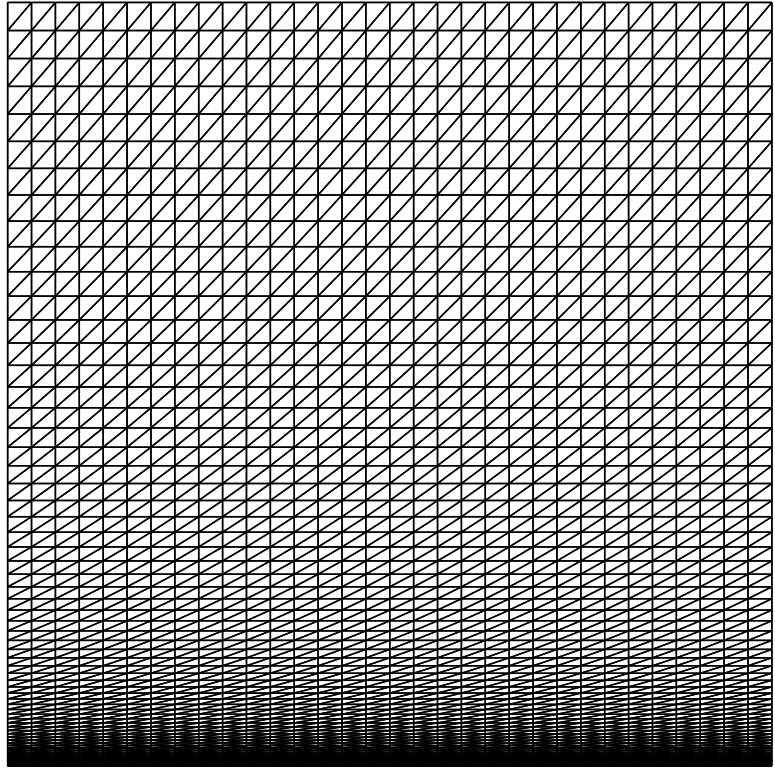}
         \caption{33x65 points}\label{bddygrid}
\end{figure}
The Newton iteration is initialized with uniform flow across the plate combined with the no-slip condition.
That is, $u$ is set to 1 everywhere but on the plate, and $v$ and $\omega$ are set to zero everywhere.
Figure \ref{bddygrid} illustrates the grid used in these computations.  On the horizontal axis $x$ varies from 0 to 1, and on the
vertical axis $y$ varies from 0 to 10.\\\\
For the problem, $\mathcal{U}=1$ and $U_\infty=1$.  Once a solution is computed it is compared to the analytical solution which is found by solving Blasius' equation:
\begin{align}
ff''+2f'''&=0\label{blasius}
\end{align}
where $f=f(\eta)$.  The boundary conditions are just,
 \begin{align}
\eta=0,\: f=0,\: f'=0;\: \eta=\infty,\:f'=1
\end{align} 
For the dimensionless form of the governing equations we have $\eta=\frac{y}{\sqrt{x}}$.  Also,
\begin{align}
u&=f'\\
v&=\frac{1}{2\sqrt{x}}(\eta f'-f)\\
\omega&=\frac{1}{\sqrt{x}}f''
\end{align}
When $y=0$ we define the shear stress $\tau$ by,
\begin{align}
\tau&=0.33206\frac{1}{\sqrt{x}}f''
\end{align}
Equation \ref{blasius} is solved using 4th order Runge-Kutta via the Shooting Method for boundary value problems.
\section{Results}
\begin{figure}[htbp]
        \centering
         \includegraphics[totalheight=4in]{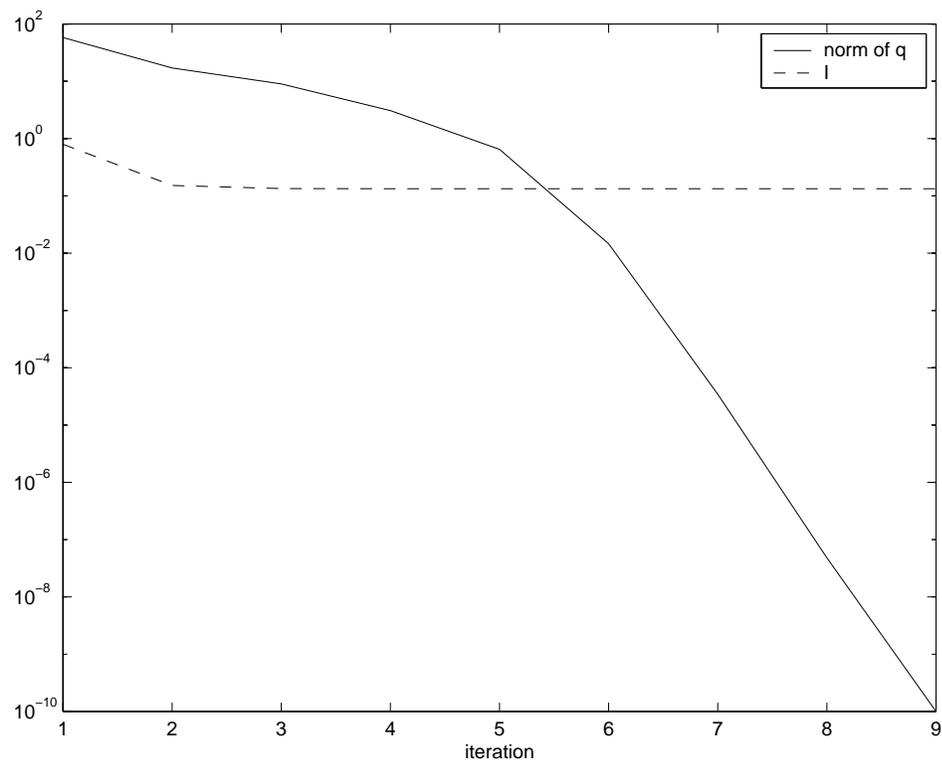}
         \caption{Semilog plot of number of iterations versus $\|\vec{q}\|,\:\mathcal{I}$}\label{res5}
\end{figure}

\begin{figure}[htbp]
        \centering
         \includegraphics[totalheight=4in]{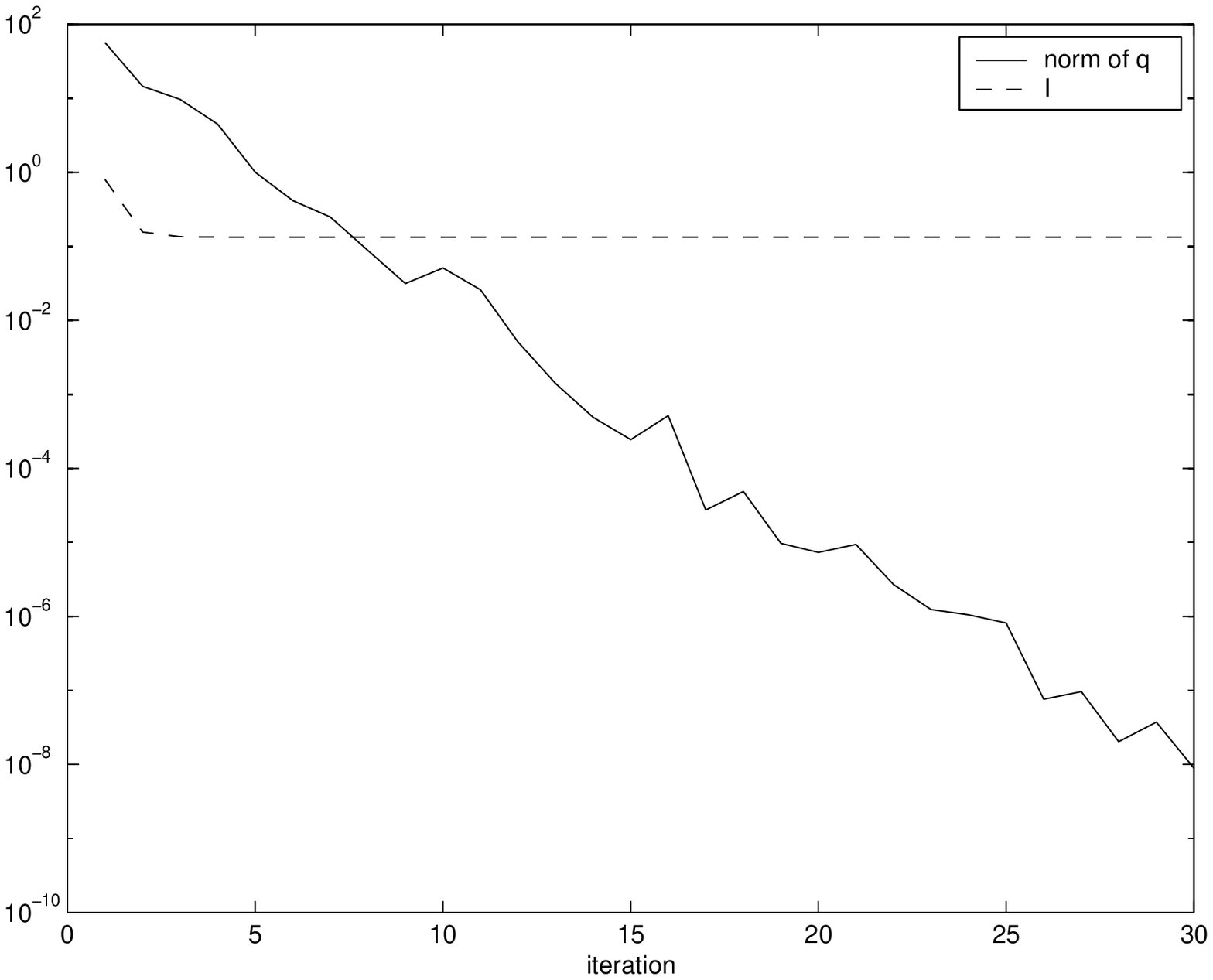}
         \caption{Semilog plot of number of iterations versus $\|\vec{q}\|,\:\mathcal{I}$}\label{res6}
\end{figure}

\begin{figure}[htbp]
        \centering
         \includegraphics[totalheight=3.7in]{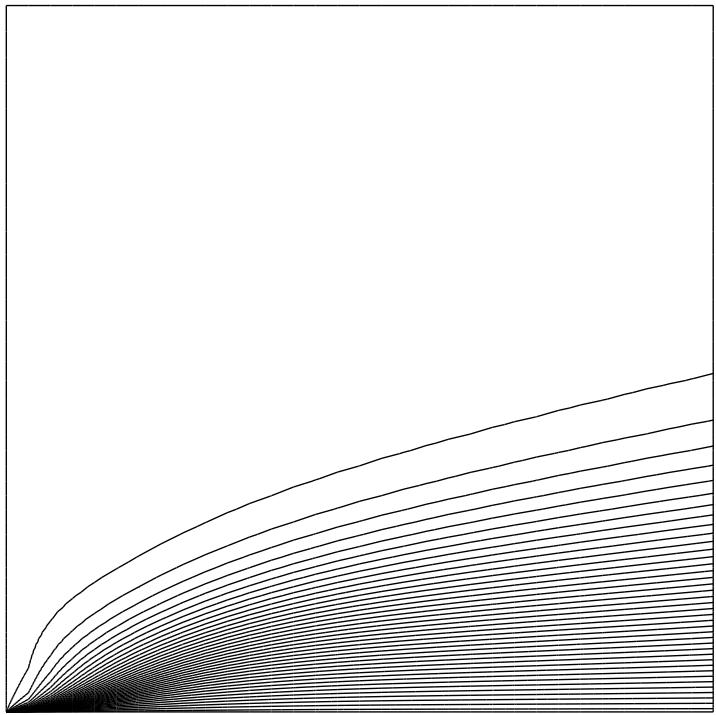}
         \caption{contour plot of $u$}\label{cont1}
\end{figure}

\begin{figure}[htbp]
        \centering
         \includegraphics[totalheight=3.7in]{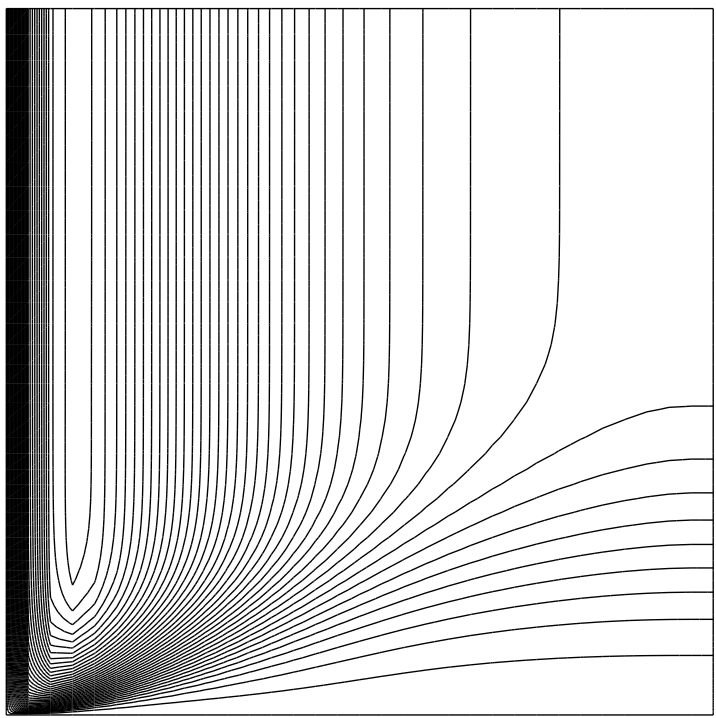}
         \caption{contour plot of $v$}\label{cont2}
\end{figure}

\begin{figure}[htbp]
        \centering
         \includegraphics[totalheight=4in]{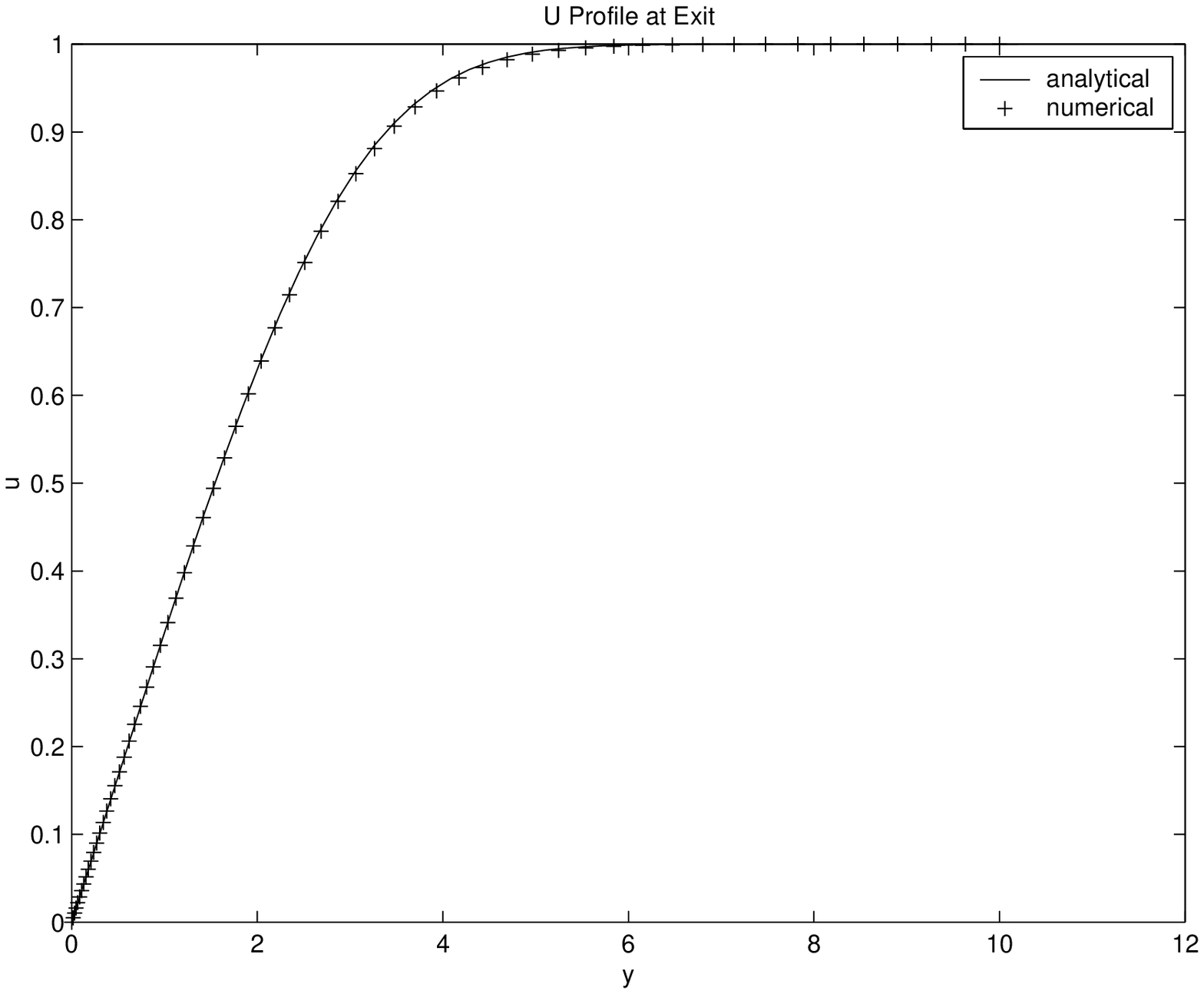}
         \caption{}\label{bp1}
\end{figure}

\begin{figure}[htbp]
        \centering
         \includegraphics[totalheight=4in]{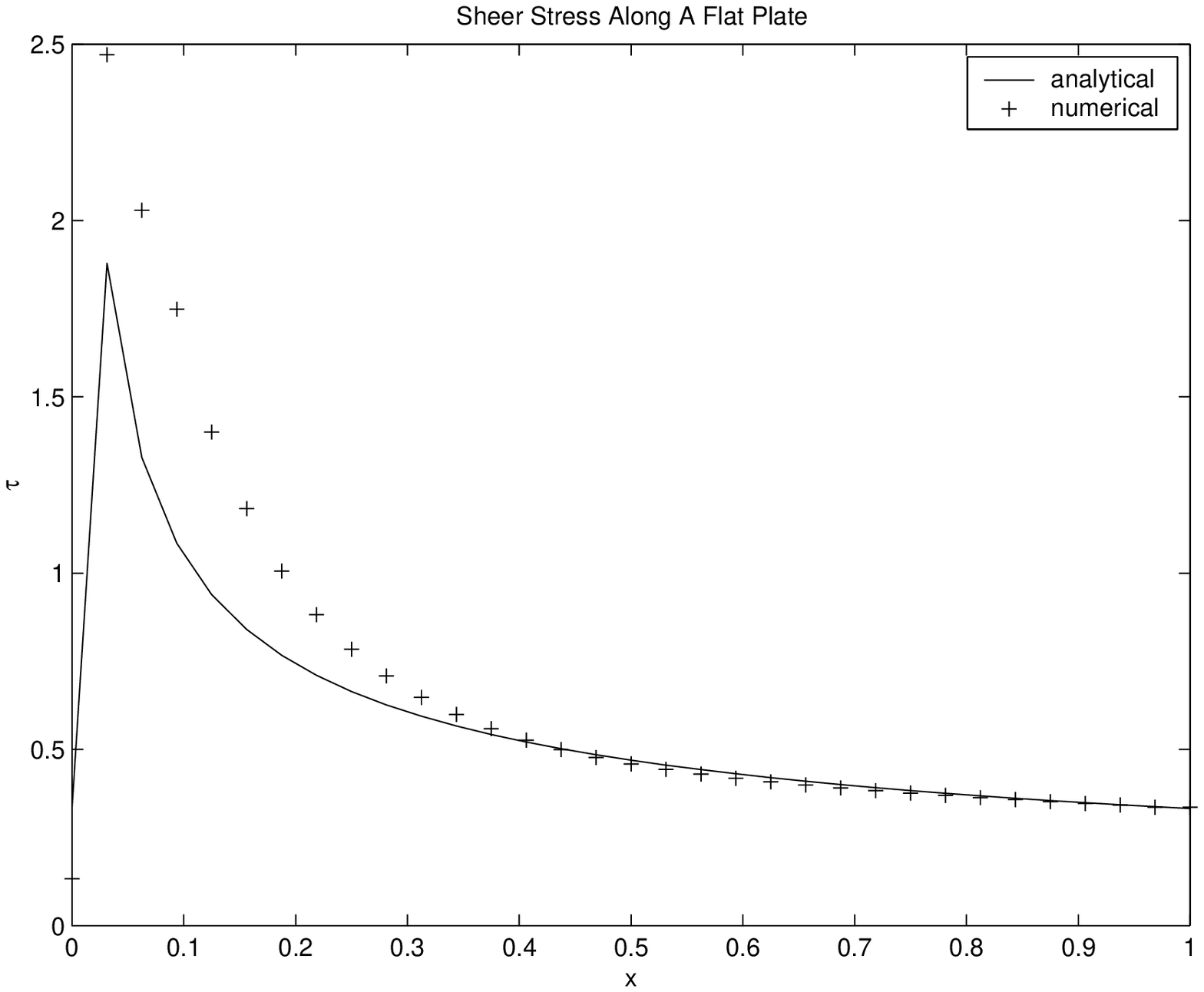}
         \caption{}\label{bp2}
\end{figure}

From figure \ref{res5}, we can conclude that the solution obtained is converged to within the prescribed tolerance.  
Although more iterations could be shown, further improvements in the corrective term $\vec{q}$ results in
negligible improvements in the functional $\mathcal{I}$.  Figure \ref{res6} is also from a converged solution that
is exactly the same as that corresponding to figure \ref{res5}.  But note that in figure \ref{res6} the correction does not
converge quadratically which should be the case for Newton's method.  The order of the method has deteriorated because
the residual has been compromised initially in order to reduce CPU time.  Whereas the correction
corresponding to figure \ref{res5} was always calculated with a residual on the order of $10^{-8}$ or less in norm, and the inner iterations
on GMRES reset after 500 iterations, those corresponding to figure \ref{res6} were reset at 100 and terminated at 400 iterations
regardless of the norm of the residual.  This initially results in a residual on the order of $10^{-3}$ but quickly
dies down and is much less than $10^{-8}$ when the solution is converged.  When executed on Linux machine equipped with a 550 Mhz dual Xeon Pentium, the result is a CPU time of 353.06 seconds.  This is a drastic improvement over a CPU time of 1570.08 seconds if we require a highly accurate solution of Newton's Method at every iteration.   
Figures \ref{cont1} and \ref{cont2} are contour 
plots of $u$ and $v$ respectively.  Figures \ref{bp1} through \ref{bp2} show there is very good agreement between the
analytical and numerical results.  

\begin{figure}[htbp]
        \centering
         \includegraphics[totalheight=4in]{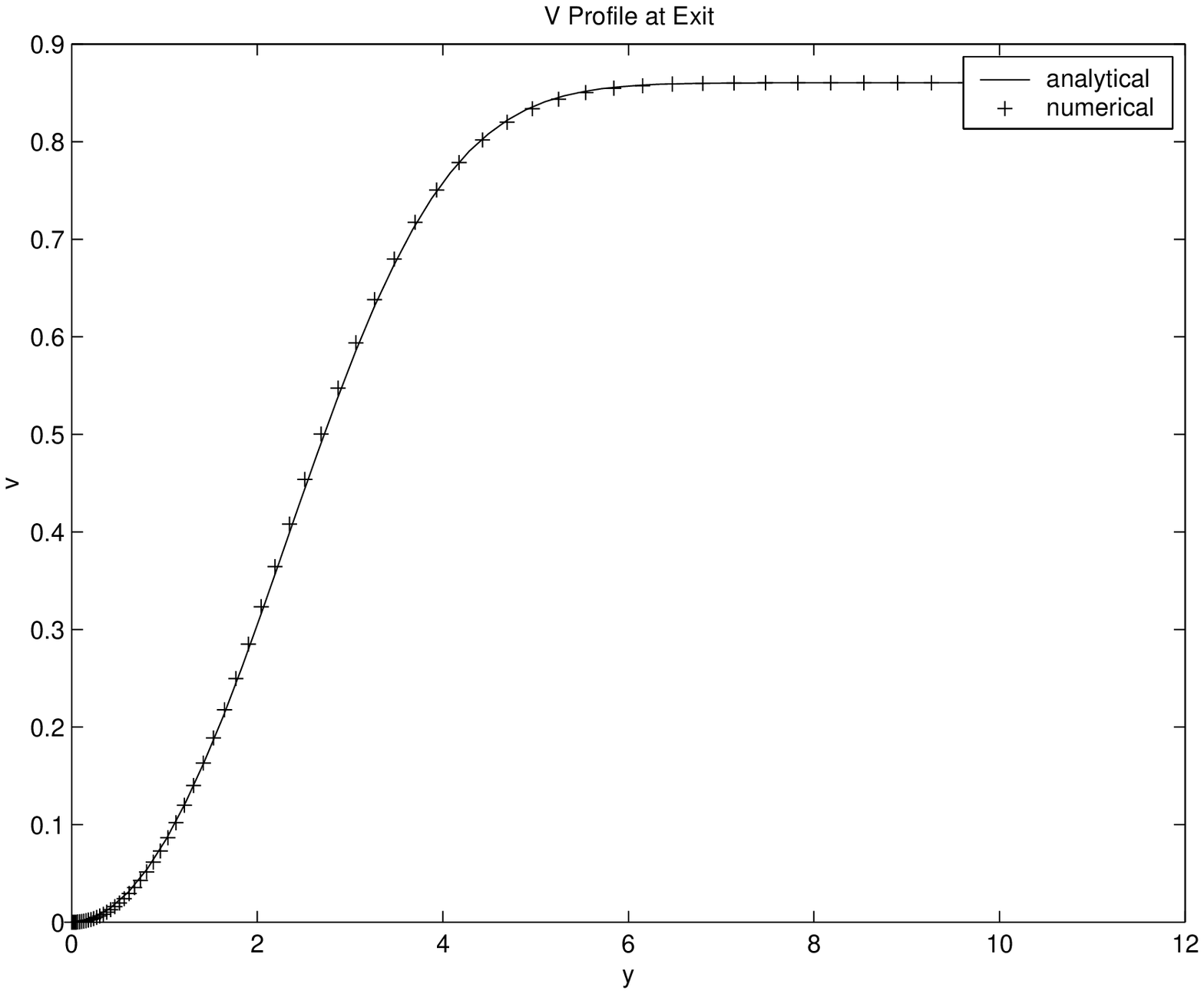}
         \caption{}
\end{figure}
\begin{figure}[htbp]
        \centering
         \includegraphics[totalheight=4in]{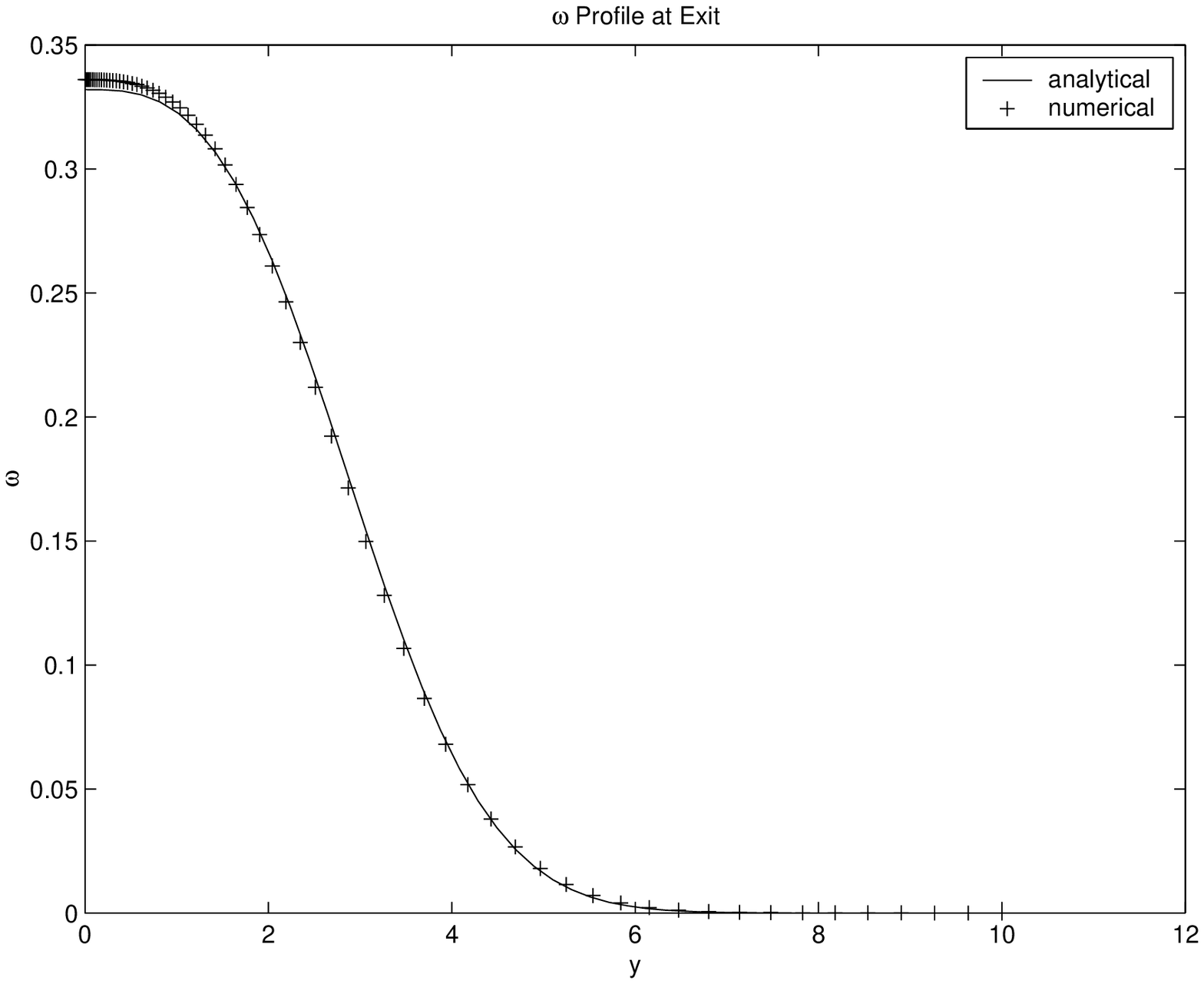}
         \caption{}
\end{figure}

\section{Conclusions And Future Projects}
In this dissertation two numerical methods were derived.  The accuracy of the methods was analyzed and the methods were applied
to problems of incompressible flow.
\paragraph{}
The first method studied was a cell-centered method.  Though the order of the local truncation error for the method in general is at most first order, the order of accuracy of the solutions obtained with the method can be higher than first order.  It was shown that when the method was applied to the solution of Laplace's equation, even in the test case that the method was not pointwise consistent a reasonable solution could be obtained.  This can in part be traced to the fact that the truncation error never diverged in this case but approached a plateau.  The method appears to be far more flexible than one would expect since it is derived by restricting the variables
to be located at the circumcenter of the triangular elements.  Relaxing this condition when necessary by placing the variables at the
centriod of the triangular elements still tends to work reasonably well.  One drawback to the method is that for potential flow, the velocities can only be calculated normal to the edges of the elements in a triangulation.  Also, enforcing Dirichlet or mixed boundary
conditions can be difficult.  Nevertheless, the method is a robust alternative to solving potential flow problems and could easily
be extended to handle advection-diffusion type equations.
\paragraph{}
The second method studied was a least-squares finite volume method.  The method was shown to be in general second order accurate.
For some of the problems solved the method worked very well.  However, for problems requiring a tangency boundary condition, there appears to be a need for further investigation into the best implementation of the boundary conditions.   In particular, it can be troublesome for the method to resolve stagnation points along such a boundary.  It is also possible that solutions can be improved if as
suggested in \cite{hub:least} the vertices of a triangulation are also considered variable and included in the minimization of the corresponding functional for a system of equations.  This has the potential to allow for instance, characteristic grid alignment for hyperbolic systems of
equations.  Also, it seems reasonable to try constructing a least-squares method where the variables are at the mid-edge instead of the vertices.
\paragraph{}
Ultimately, solutions to the Navier-Stokes equation with a least-squares finite volume method is desired.  The methodology should be extendible to three-dimensional problems as well.

\nocite{*}
\bibliographystyle{unsrt}
\bibliography{mybib}

\begin{thebibliography}{10}

\bibitem{chattot2}
Jean-Jacques Chattot.
\newblock M\'{e}thode variationnelle pour les probl\'{e}mes hyperboliques et
  mixtes du premier ordre.
\newblock In R.~Glowinski and J.~L. Lions, editors, {\em Computing Methods in
  Applied Sciences and Engineering}. North-Holland Publishing Company, 1980.

\bibitem{hub:least}
M.~E. Hubbard.
\newblock {Least Squares Minimisation and Steepest Descent Methods for the
  Scalar Advection Equation and a Cauchy-Riemann System on an Adaptive Grid}.
\newblock Technical report, Department of Mathematics, The University of
  Reading, P. O. Box 220, Whiteknights, Reading, Berkshire, RG6 6AX, U. K.,
  September 1997.
\newblock
  "{\urlBiBTeX{http://www.rdg.ac.uk/AcaDepts/sm/wsm1/research/na-reports-1997/%
na_6-97.ps}}".

\bibitem{herbin}
Rapha\'{e}le Rubin.
\newblock An error estimate for a finite volume scheme for a diffusion
  convection problem on a triangular mesh.
\newblock {\em Numerical Methods for P.D.E.}, 11(2):165--174, 1995.

\bibitem{chattot}
Jean-Jacques Chattot.
\newblock A test of the formal accuracy of numerical schemes on unstructured
  meshes.
\newblock {\em Computational Fluid Dynamics Journal}, 9, April 2000.
\newblock 8th International Symposium on Computational Fluid Dynamics, Bremen,
  Sep. 5-10, 1999.

\bibitem{lazarov}
Raytcho~D. Lazarov and Ilya~.D. Mishev.
\newblock Finite volume methods for reaction-diffusion problems.
\newblock Technical Report ISC-96-07-MATH, Texas A\&M University, College
  Station, TX, 1996.
\newblock "{\urlBiBTeX{http://www.isc.tamu.edu/iscpubs/9607.ps}}".

\bibitem{vaz}
Guilherme Vaz.
\newblock {2D Boundary Element Application - The Morino Method}.
\newblock Technical report, Instituto Superior T\'{e}cnico, Av. Rovisco Pais,
  1049-001 Lisboa, Portugal, South America, March 2000.
\newblock "{\urlBiBTeX{http://www.math.ist.utl.pt/~gvaz/Files/Report2.pdf}}".

\bibitem{shewchuk}
Jonathan Shewchuk.
\newblock {Triangle: Engineering a 2D Quality Mesh Generator and Delaunay
  Triangulator}.
\newblock In {\em First Workshop on Applied Computational Geometry}, pages
  124--133, Philadelphia, Pennsylvania, May 1996. ACM.
\newblock "{\urlBiBTeX{http://www-2.cs.cmu.edu/~quake/triangle.html}}".

\bibitem{alcrudo}
Francisco Alcrudo and Jean-Jacques Chattot.
\newblock {Mid-Edge Box Schemes for Hyperbolic Partial Differential Equations}.
\newblock In Hermes, editor, {\em Finite Volumes for Complex Applications},
  pages 301--308, 1996.

\bibitem{baines:mov}
M.~J. Baines.
\newblock {Moving Finite Element and Moving Least Squares Approximation of
  Time-dependent and Steady PDEs in Multidimensions}.
\newblock Technical report, Department of Mathematics, The University of
  Reading, P. O. Box 220, Whiteknights, Reading, Berkshire, RG6 6AX, U. K.,
  December 2000.
\newblock
  "{\urlBiBTeX{http://www.rdg.ac.uk/AcaDepts/sm/wsm1/dsgs/publs/mjb1.pdf}}".

\bibitem{darcha}
Darryl Whitlow and Jean-Jacques Chattot.
\newblock {A Least-Squares Finite Volume Method For Viscous Laminar Flow}.
\newblock Technical report, Mathematics Department, University of California,
  Davis, March 2001.

\bibitem{chattot3}
Jean-Jacques Chattot.
\newblock Box schemes for first order partial differential equations.
\newblock In Wagdi~G. Habashi and Mohamed~M. Hafez, editors, {\em Computational
  Fluid Dynamics Techniques}, pages 307--331. Gordon \& Breach Publishing
  Group, 1995.

\bibitem{white}
Frank White.
\newblock {\em Fluid Mechanics}.
\newblock McGraw-Hill, Inc., second edition, 1986.

\bibitem{baines:sol}
M.~J. Baines.
\newblock {The Solution of Steady PDEs on Adjustable Meshes in Multidimensions
  Using Local Descent Methods}.
\newblock Technical report, Department of Mathematics, The University of
  Reading, P. O. Box 220, Whiteknights, Reading, Berkshire, RG6 6AX, U. K.,
  2001.
\newblock
  "{\urlBiBTeX{http://www.rdg.ac.uk/AcaDepts/sm/wsm1/dsgs/publs/mjb2.pdf}}".

\bibitem{bertin}
John~J. Bertin and Micheal~L. Smith.
\newblock {\em Aerodynamics For Engineers}.
\newblock Prentice-Hall, Englewood Cliffs, New Jersey, second edition, 1989.

\bibitem{burden}
Richard Burden and J.~Douglas Faires.
\newblock {\em Numerical Analysis}.
\newblock Brooks/Cole Publishing Company, sixth edition, 1997.

\bibitem{chorin}
Alexandre~J. Chorin and Jerrold~E. Marsden.
\newblock {\em A Mathematical Introduction to Fluid Mechanics}.
\newblock Springer-Verlag, third edition, 1993.

\bibitem{darryl}
Darryl Whitlow.
\newblock {A Finite Volume Method For The Solution of Potential Flow}.
\newblock Technical report, Mathematics Department, University of California,
  Davis, June 1997.

\bibitem{hermann}
Hermann Schlichting.
\newblock {\em Boundary-Layer Theory}.
\newblock McGraw-Hill Classic Textbook Reissue. McGraw-Hill, Inc., seventh
  edition, 1979.

\end{thebibliography}
\appendix
\chapter{Order of Accuracy For the Cell-Centered Method}
\begin{proposition}
On a grid composed of equilateral triangles
\begin{equation}
e_{_\Delta}=\overline{\Delta\psi}=\frac{\int_{\Omega_1}\Delta\psi\,dx dy}{\Omega_1}=\mathcal{O}(h)
\end{equation} 
\end{proposition}
Proof:\\ 
Recall (refer to the figure below),
\begin{eqnarray}0 &=&\int_{\Omega_1}\Delta\psi\,dx dy\nonumber\\
&=&\frac{\psi_4\,l_{31}}{d_{14}}+\frac{\psi_3\,l_{23}}{d_{13}} +\frac{\psi_2\,l_{12}}{d_{12}}-
\psi_1(\frac{l_{31}}{d_{14}}+\frac{l_{23}}{d_{13}}+\frac{l_{12}}{d_{12}})\nonumber
\end{eqnarray}
\begin{figure}[htbp]
        \centering
         \includegraphics[totalheight=3in]{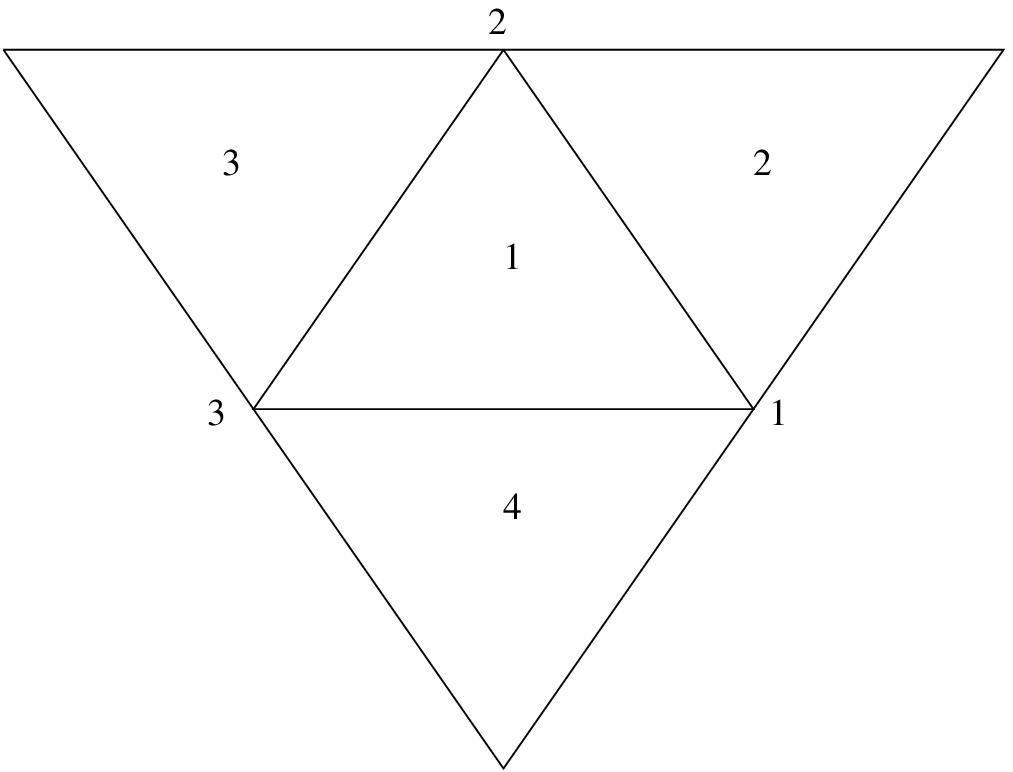}
\end{figure}
For a set of equilateral triangles this reduces to:
\begin{equation}
0=\sqrt{3}(\psi_4+\psi_3+\psi_2-3\psi_1)
\end{equation}
Let $l$ equal the length of a side of a triangle.  Define $h=\delta x =\frac{l}{2}$ and choose
\begin{equation}\label{iden}
k=\delta y=\frac{l}{\sqrt{3}}=\frac{2 \delta x}{\sqrt{3}}
\end{equation}
Then $\Omega_1=h^2\sqrt{3}$.  With respect to the figure above, let the origin be the circumcenter of element numbered one and let $\Psi$ equal the solution to Laplace's equation over the domain of the triangulation.  Then for the local truncation error we have,
\begin{align*}
\Omega_1 e_{_\Delta}&=\sqrt{3}(\Psi_4+\Psi_3+\Psi_2-3\Psi_1)\\
h^2 e_{_\Delta}&=(\Psi_4+\Psi_3+\Psi_2-3\Psi_1)\\
 &=\Psi (x,y-k)+\Psi (x-h,y+\frac{k}{2})+\Psi (x+h,y+\frac{k}{2})-3\Psi(x,y)\\
 &=\left(\Psi(x,y)-k\Psi_y+\frac{k^2}{2}\Psi_{yy}-\frac{k^3}{6}\Psi_{yyy}+\ldots \right)\\
 & \quad + \left(\Psi(x,y) -h\Psi_x-\frac{k}{2}\Psi_y+\frac{1}{2}(h^2\Psi_{xx}-hk\Psi_{xy}\right.
         +\frac{k^2}{4}\Psi_{yy})\\
 &\quad +\left.\frac{1}{6}(-h^3\Psi_{xxx}+\frac{3h^2k}{2}\Psi_{xxy}
        -\frac{3hk^2}{4} \Psi_{xyy}+\frac{k^3}{8}\Psi_{yyy})+\ldots\right)\\
 & \quad + \left(\Psi(x,y) +h\Psi_x+\frac{k}{2}\Psi_y+\frac{1}{2}(h^2\Psi_{xx}+hk\Psi_{xy}
         +\frac{k^2}{4}\Psi_{yy})\right.\\
 &\quad +\frac{1}{6}(h^3\Psi_{xxx}+\frac{3h^2k}{2}\Psi_{xxy}
        +\left.\frac{3hk^2}{4} \Psi_{xyy}+\frac{k^3}{8}\Psi_{yyy})+\ldots\right)\\
 &\quad -3\Psi(x,y)\\
 &=\frac{k^2}{2}\Psi_{yy}-\frac{k^3}{8}\Psi_{yyy}+h^2\Psi_{xx}+\frac{k^2}{4}\Psi_{yy}
   +\frac{h^2k}{2}\Psi_{xxy}+\ldots\\
\end{align*}
\begin{align*}
h^2 e_{_\Delta} &=h^2\Psi_{xx}+\frac{3k^2}{4}\Psi_{yy}+\frac{h^2k}{2}\Psi_{xxy}-\frac{k^3}{8}\Psi_{yyy}+\ldots\\
 &=h^2\Psi_{xx}+\frac{3}{4}\left(\frac{4h^2}{3}\right)\Psi_{yy}+\frac{h^3}{\sqrt{3}}\left(\Psi_{xxy}-\frac{1}{3}\Psi_{yyy}
\right)+\ldots
 &&\text{by \eqref{iden}}
\end{align*}
Division by $h^2$ completes the proof:
\begin{equation}
\boxed{e_{_\Delta}=\frac{h}{\sqrt{3}}\left(\Psi_{xxy}-\frac{1}{3}\Psi_{yyy}\right)+\mathcal{O}(h^2)} 
\end{equation} 
From here it is easy to see why the method is exact for test case 1 in chapter 2 since for
$\Psi=(x-1)^3 -3(x-1)y^2$ we have $\Psi_{xxy}=\Psi_{yyy}=0$.
\chapter{On The Decoupling Of The Jacobian For The Cauchy-Riemann Equations}
For the least-squares method studied in this work, here it is shown that the resulting Jacobian used
in the minimization process of the Cauchy-Riemann equations is very simple if the grid is composed entirely
of congruent right triangles.  
This can result in a reduction in computer time if Dirichlet boundary conditions are enforced 
for a perspective problem.  The resulting Jacobian will have a structure that can be exploited by a particular solver for large problems. It is due to the fact that the equations
in the Jacobian corresponding to the minimization of any particular variable contains no entries due to contributions from the other
variables.  
This is clarified by the following proposition:   
\begin{proposition}
Let $u,v$ be piecewise linear functions on a triangulation of a domain $\Omega$.
If $u,v$ satisfy the aforementioned discrete functional for the Cauchy-Riemann equations on $\Omega$
and the triangulation is composed entirely of congruent right triangles then at all inner nodes of the domain we have, 
\begin{equation}
\sum_T\frac{\partial ^2 \mathcal{I}_{_T}}{\partial v_j\partial u_i}=0\label{cross}
\end{equation}
where $T$ indexes all elements abutting the node $i$.
\end{proposition}
Proof:\\
Recall that for the Cauchy-Riemann equations we have,
\begin{equation}\label{cross2}
\frac{\partial ^2\mathcal{I}_{_T}}{\partial v_j \partial u_i}
=\frac{1}{4\Omega_{_T}}(n_{x_i}n_{y_j}-n_{y_i}n_{x_j})=\frac{1}{4\Omega_{_T}}(\mathbf n_i \times \mathbf n_j)\cdot \mathbf e_3
\end{equation}
where $\mathbf e_1=\langle 1,0,0 \rangle$, $\mathbf e_2=\langle 0,1,0 \rangle$ and $\mathbf e_3=\langle 0,0,1 \rangle$.  
Without loss of generality, assume that the triangles have sides of unit length and suppress the term $\frac{1}{4\Omega_{_T}}$.  From equation
(\ref{cross2}) the case $i=j$ is trivial.  Therefore, due to symmetry, we need only show
equation (\ref{cross}) holds at one node where $i \neq j$ for configurations of the figures below. 
\begin{figure}[H]
\begin{minipage}[t]{.35\linewidth}
\centering
\includegraphics[width=1.7in]{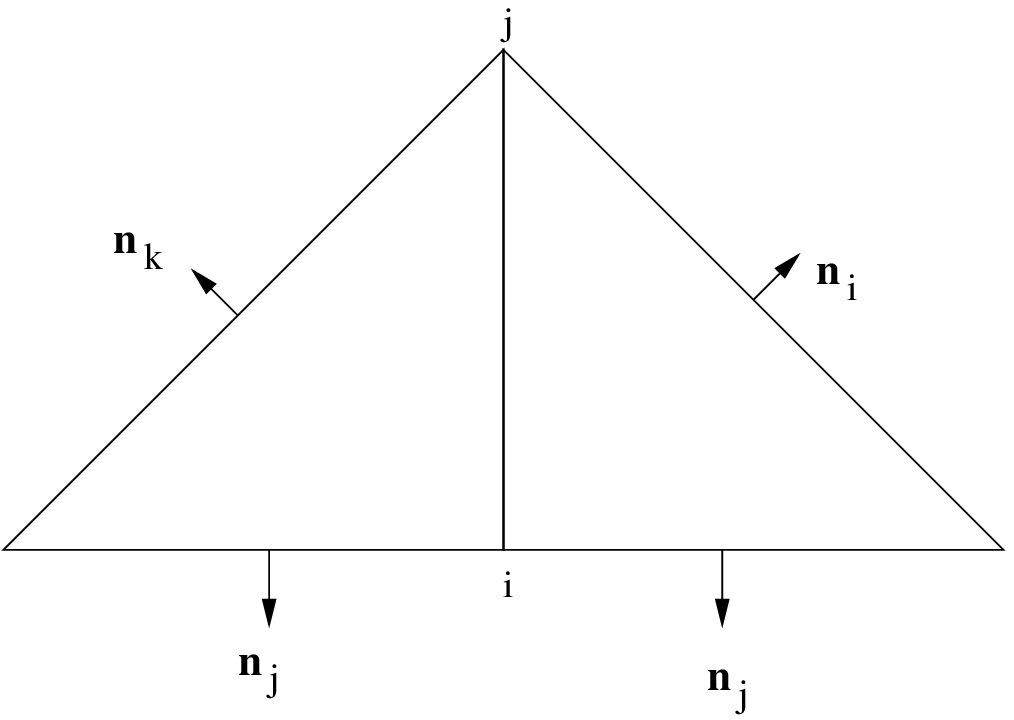}
         \caption{$1^{st}$ config.\label{config1}}
\end{minipage}%
\begin{minipage}[t]{.35\linewidth}
\centering
\includegraphics[width=1.7in]{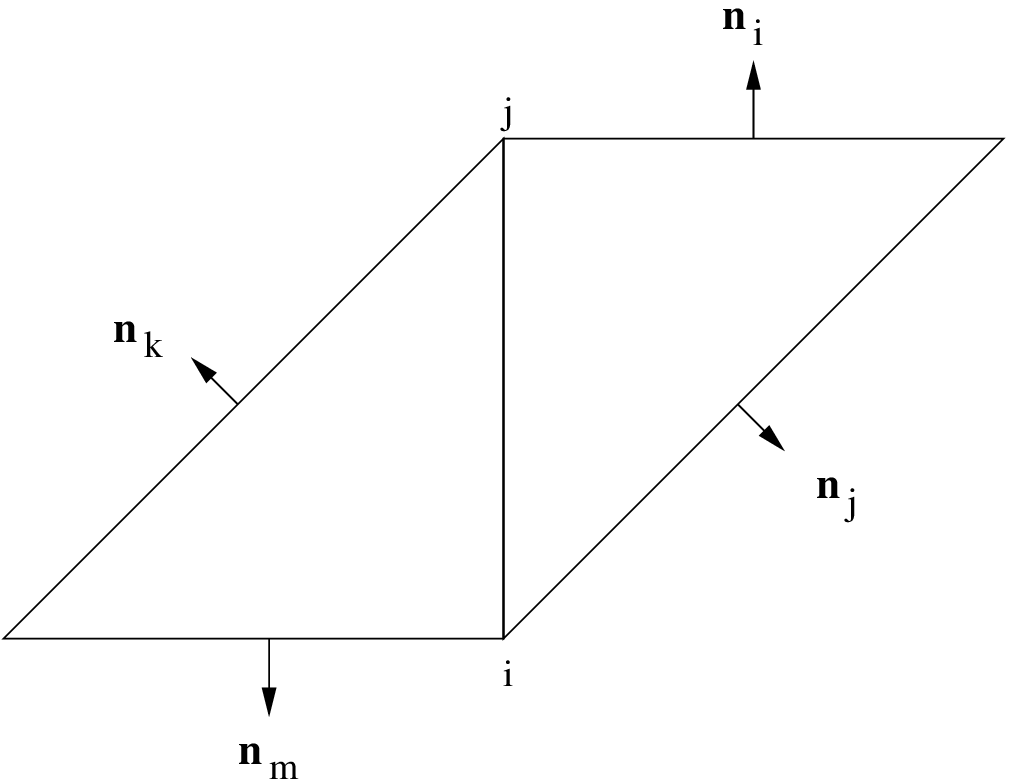}
         \caption{$2^{nd}$ config.\label{config2}}
\end{minipage}%
\begin{minipage}[t]{.35\linewidth}
        \centering
         \includegraphics[height=1.3in, width=1.3in]{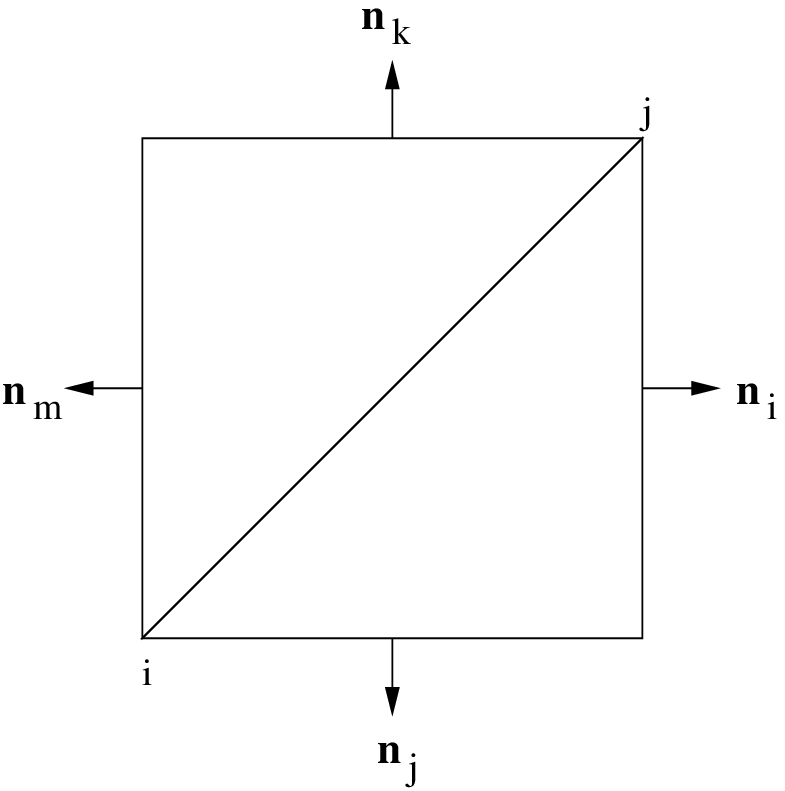}
         \caption{$3^{rd}$ config.\label{config3}}
\end{minipage}
\end{figure}
\underline{$1^{st}$ Configuration}\\               
From equations (\ref{cross}) and (\ref{cross2}) we have,
\begin{align}
\sum_T\frac{\partial ^2 \mathcal{I}_{_T}}{\partial v_j\partial u_i}&=(\mathbf n_i \times \mathbf n_j)\cdot \mathbf e_3+
(\mathbf n_k \times \mathbf n_j)\cdot \mathbf e_3\\
&=((\mathbf e_1-\mathbf n_j)\times \mathbf n_j )\cdot \mathbf e_3+((-\mathbf e_1-\mathbf n_j)\times \mathbf n_j )\cdot \mathbf e_3\\
&= -2(\mathbf n_j \times \mathbf n_j)\cdot \mathbf e_3\\
&=0
\end{align}
\hrulefill\\ 
\underline{$2^{nd}$ Configuration}\\ 
\begin{align}
\sum_T\frac{\partial ^2 \mathcal{I}_{_T}}{\partial v_j\partial u_i}&=(\mathbf n_i \times \mathbf n_j)\cdot \mathbf e_3+
(\mathbf n_k \times \mathbf n_m)\cdot \mathbf e_3\\
&=(\mathbf n_i \times (\mathbf e_1 - \mathbf e_2))\cdot \mathbf e_3+((\mathbf e_2 - \mathbf e_1)\times(-\mathbf n_i))\cdot \mathbf e_3\\
&=(\mathbf n_i \times (\mathbf e_1 - \mathbf e_2))\cdot \mathbf e_3-(\mathbf n_i \times (\mathbf e_1 - \mathbf e_2))\cdot \mathbf e_3\\
&=0
\end{align}
\hrulefill\\ 
\underline{$3^{rd}$ Configuration}\\ 
\begin{align}
\sum_T\frac{\partial ^2 \mathcal{I}_{_T}}{\partial v_j\partial u_i}&=(\mathbf n_i \times \mathbf n_j)\cdot \mathbf e_3+
(\mathbf n_k \times \mathbf n_m)\cdot \mathbf e_3\\
&=(\mathbf n_i \times \mathbf n_j)\cdot \mathbf e_3+((-\mathbf n_j)\times (-\mathbf n_i))\cdot \mathbf e_3\\
&=(\mathbf n_i \times \mathbf n_j)\cdot \mathbf e_3-(\mathbf n_i \times \mathbf n_j)\cdot \mathbf e_3\\
&=0
\end{align}
This completes the proof.  Incidentally, this result holds if the triangulation is composed of congruent equilateral triangles.  Referring to the figure below,
\begin{figure}[H]
        \centering
         \includegraphics[totalheight=3in]{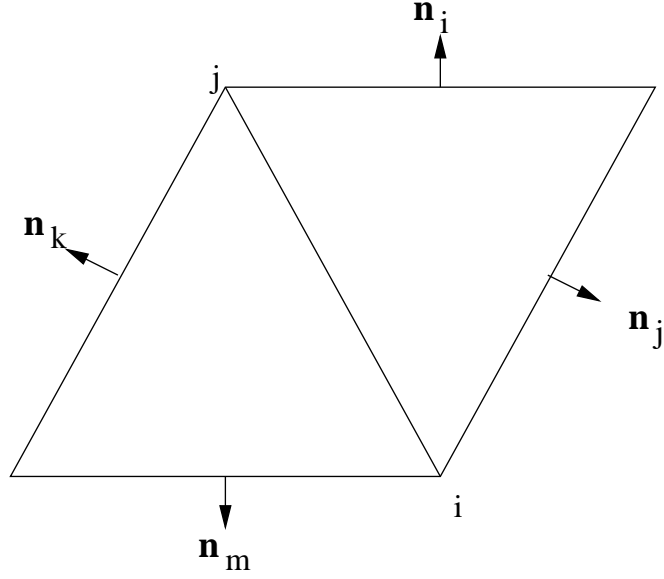}\caption{Configuration of equilateral triangles} 
\end{figure}
Then just as in the previous configuration we have,
\begin{align}
\sum_T\frac{\partial ^2 \mathcal{I}_{_T}}{\partial v_j\partial u_i}&=(\mathbf n_i \times \mathbf n_j)\cdot \mathbf e_3+
(\mathbf n_k \times \mathbf n_m)\cdot \mathbf e_3\\
&=(\mathbf n_i \times \mathbf n_j)\cdot \mathbf e_3+((-\mathbf n_j)\times (-\mathbf n_i))\cdot \mathbf e_3\\
&=(\mathbf n_i \times \mathbf n_j)\cdot \mathbf e_3-(\mathbf n_i \times \mathbf n_j)\cdot \mathbf e_3\\
&=0
\end{align}
This property can be used to trace errors in the coding of the numerical method.
\end{document}